\documentclass[12pt]{article}

\usepackage{amscd,amsfonts,amsmath,amssymb,latexsym,amsthm}
\usepackage[mathscr]{eucal}

\input{xy}
\xyoption{all}

\theoremstyle{plain}
\newtheorem{thm}{Theorem}[section]
\newtheorem{prop}[thm]{Proposition}
\newtheorem{lemma}[thm]{Lemma}
\newtheorem{cor}[thm]{Corollary}
\newtheorem*{tthm}{Technical Theorem}
\newtheorem*{ktt}{Kasparov's Technical Theorem}

\theoremstyle{definition}
\newtheorem{definition}[thm]{Definition}
\newtheorem{example}[thm]{Example}

\newtheorem{goodrmk}[thm]{Remark}
\newtheorem{goodclaim}[thm]{Claim}
\newtheorem{overview}[thm]{Overview}

\theoremstyle{remark}
\newtheorem*{rmk}{Remark}

\def\singlespace{%
\vskip\parskip%
\vskip\baselineskip%
\def\baselinestretch{1}%
\ifx\@currsize\normalsize\@normalsize\else\@currsize\fi%
\vskip-\parskip%
\vskip-\baselineskip%
}

\renewcommand{\qedsymbol}{$\blacksquare$}
\newcommand{\lam}{\lbrack\!\lbrack}
\newcommand{\ram}{\rbrack\!\rbrack\,}
\newcommand{\asy}{\dasharrow}

\newcommand{\C}{\mathbb{C}}
\newcommand{\CI}{C_0((-1,1))}

\newcommand{\pcli}[1]{\mathcal{C}_{+#1}}
\newcommand{\cc}{C_c^{\infty}}
\newcommand{\cs}{$C^*$}
\newcommand{\cscateg}{$\mathbf{C^*}${\bf-alg}}
\newcommand{\mcs}[1]{C^*(#1)}

\newcommand{\der}{\partial}

\newcommand{\eps}{\varepsilon}
\newcommand{\gcs}[1]{$#1$-$C^*$}
\newcommand{\gcscateg}{$\mathbf{G}${\bf-}$\mathbf{C^*}${\bf-alg} }
\newcommand{\hm}{{\mathcal{E}}_{\bullet}}
\newcommand{\hs}{\mathcal{H}}
\newcommand{\ie}{{\em i.e.\;}}
\newcommand{\inp}[2]{\langle #1, #2\rangle}
\newcommand{\lra}{\longrightarrow}
\newcommand{\mic}{\vspace{.2cm}}
\newcommand{\mare}{\vspace{.5cm}}

\newcommand{\ra}{\rightarrow}
\newcommand{\R}{\mathbb{R}}
\newcommand{\res}[2]{(\text{\rm Res}_{#1})_*\left(\, #2 \,\right)} %
\newcommand{\SSS}{\mathbb{S}}
\newcommand{\s}{\mathcal{S}}
\newcommand{\sh}{$*$-}
\newcommand{\stack}[1]{\xrightarrow{#1}}
\newcommand{\ui}{t\in [1,\infty)}

\newcommand{\Z}{\mathbb{Z}}

\newcommand{\bdd}{\mathscr{B}}                
\newcommand{\cpct}{\mathscr{K}}               
\newcommand{\loc}{\mathscr{J}}                
\newcommand{\cval}{\mathscr{C}}               
\newcommand{\mlt}{\mathscr{M}}                
\newcommand{\uni}{\mathscr{U}}                

\newcommand{\tpr}{\otimes}
\newcommand{\tprm}{\otimes_{max}}

\newcommand{\itpr}[1]{\otimes_{\scriptscriptstyle{#1}}}

\newcommand{\Kprod}[3]{\,{#1}\,\sharp_{\scriptscriptstyle{#2}}\,{#3}\,}   
\newcommand{\Mh}{M^{\scriptscriptstyle{\frac12}}}

\newcommand{\Nh}{N^{\scriptscriptstyle{\frac12}}}

\newcommand{\sg}[1]{\sigma_{\scriptscriptstyle{#1}}}
\newcommand{\kke}{\Theta}
\newcommand{\kee}{\Xi}
\newcommand{\keclass}[1]{\boldsymbol{[\!\![\!\![}\,#1\,\boldsymbol{]\!\!]\!\!]}}
\newcommand{\kkclass}[1]{\boldsymbol{#1}}
\newcommand{\eclass}[1]{\lam #1 \ram}

\newcommand{\fm}[1]{{\mathcal{E}}_{#1}}         
\newcommand{\ofm}[1]{\widehat{{\mathcal{E}}_{#1}}} 
\newcommand{\fr}[1]{\varphi_{#1}}               
\newcommand{\fop}[1]{F_{#1}}                    
\newcommand{\ofop}[1]{\widehat{F_{#1}}}            
\newcommand{\fv}[1]{#1}                         
\newcommand{\cnn}[1]{\underline{#1}}            

\newcommand{\et}{\mathcal{E}_t}

\newcommand{\E}{\mathcal{E}}
\newcommand{\ec}{\mathcal{E}}
\newcommand{\ecu}{\mathcal{E}_{1}}
\newcommand{\ecd}{\mathcal{E}_{2}}
\newcommand{\cfm}[2]{continuous field of ($#1$,$#2$)-modules}
\newcommand{\akm}[2]{asymptotic Kasparov ($#1$,$#2$)-module}

\newcommand{\ga}{\alpha}
\newcommand{\gb}{\beta}
\newcommand{\gd}{\delta}
\newcommand{\gD}{\Delta}

\newcommand{\gf}{\phi}
\newcommand{\gvf}{\varphi}
\newcommand{\gG}{\Gamma}
\newcommand{\gi}{\iota}
\newcommand{\gl}{\lambda}
\newcommand{\go}{\omega}
\newcommand{\gs}{\sigma}
\newcommand{\gS}{\Sigma}

\newcommand{\gth}{\theta}
\newcommand{\gz}{\zeta}
\begin{document}

\title{On an Intermediate Bivariant Theory for $C^*$-algebras, I}
\author{
Dorin Dumitra\c{s}cu
\thanks{Dartmouth College, Department of Mathematics, 6188 Bradley Hall,
Hanover, NH 03755, USA, {\tt dumitras@hilbert.dartmouth.edu}}
}

\author{
Dorin Dumitra\c{s}cu}

\maketitle

%
%
%
%

\begin{abstract}
We construct a new bivariant theory, 
that we call $KE$-theory, which is intermediate between
the $KK$-theory of G.~G.~Kasparov, and the $E$-theory of 
A.~Connes and N.~Higson. 
For each pair of separable graded $C^*$-algebras $A$ and $B$, acted upon
by a locally compact $\gs$-compact group $G$, we define an abelian
group $KE_G(A,B)$. We show that there is an associative product
$KE_G(A,D) \otimes KE_G(D,B) \ra KE_G(A,B)$. 
Various functoriality properties of the $KE$-theory groups and of 
the product are presented. 
The new theory has a simpler product than $KK$-theory and there
are natural transformations
$KK_G \rightarrow KE_G$ and $KE_G \rightarrow E_G$. 
The complete description of these maps will form the substance of a second 
paper.
\end{abstract}

\mare\noindent

\mic\noindent

\tableofcontents


\setcounter{section}{-1}
\section{Introduction}\label{C:intro}

In the 1960's, the algebraic topology of manifolds produced 
one of the profound theorems of $\text{XX}^{\text{th}}$ 
century mathematics: the Index Theorem of Atiyah and Singer 
(\cite{AtSi63}, \cite{Pals}). The conceptual proof given in 
\cite{AtSiI} is based on a cohomology theory invented by 
M.~Atiyah and F.~Hirzebruch \cite{AtHr} --- namely $K$-theory. 
Using hints coming from various generalizations 
of the index theorem, Atiyah \cite{Atiy69} also proposed a way of 
defining cycles of the {\em dual} theory --- namely $K$-homology.
The only thing left open by Atiyah was the definition of the 
equivalence relation that would make these cycles into a group. 
This issue was resolved by G.~G.~Kasparov \cite{Kas75}. 
He succeeded in creating (see also \cite{Kas81}) a {\em bivariant} 
theory --- named $KK$-theory --- which associates to any two 
\cs -algebras $A$ and $B$ a group $KK(A,B)$. His theory generalizes
both $K$-theory for compact manifolds
(obtained when $A$ is $\C$, and $B$ is the continuous 
functions on the manifold) and $K$-homology (obtained when $A$ is
the continuous functions on the manifold, and $B$ is $\C$). For a
very well written account of the origins of $KK$-theory see the 
introductory sections of \cite{Hg87a} and \cite{Hg90a}.

Besides a wealth of functorial properties, 
the key feature of $KK$-theory
is the existence for any separable \cs -algebras 
$A$, $B$, and $D$ of an associative product map
$
KK(A,D) \tpr KK(D,B) \longrightarrow KK(A,B).
$
Following an approach 
indicated by J.~Cuntz (\cite{Cu83}, \cite{Cu84}), N.~Higson 
\cite{Hg87a} gave the following description of $KK$-theory:
it is the universal category with homotopy invariance,
stability, and split-exactness. This category has separable 
\cs -algebras as objects, elements of $KK$-groups as morphisms, and the
Kasparov product as the composition of morphisms.

In a subsequent paper \cite{Hg90b} Higson described the universal 
category with homotopy invariance, stability, and {\em exactness}. The 
resulting new theory --- named $E$-theory --- has become important 
in \cs -algebra theory after A.~Connes and N.~Higson \cite{CoHg90} 
described it concretely in terms of asymptotic morphisms.
(An asymptotic morphism between two \cs -algebras is a family of maps
between the two, indexed by $[1,\infty)$, which satisfies the conditions
of a \sh homomorphism in the limit at $\infty$.) The
description of $KK$-theory and $E$-theory using category theory
implies, in a rather abstract and algebraic way, the existence of a
map $KK(A,B)\ra E(A,B)$, for any two \cs -algebras $A$ and $B$.
This map is an isomorphism when $A$ is nuclear \cite{Sk88}. 
Similar descriptions of the universality property for the
{\em equivariant} theories are also known: 
for the equivariant $KK$-theory \cite{Kas88} 
under the action of a group see \cite{Thms98},
for the equivariant $E$-theory 
under the action of a group see \cite{GHT},
and for both theories under the action of a groupoid see \cite{Pop}.

Equivariant $KK$-theory and $E$-theory have become essential tools
in \cs -algebra theory because of their use in solving 
topological/geometrical problems,
notably cases of the Novikov conjecture \cite{Kas88}, \cite{Rsn84}, 
and the Baum-Connes conjecture \cite{BC82}, \cite{BCH}.

In this paper a new theory is constructed, that we call
{\em $KE$-theory}, which is intermediate between $KK$-theory and 
$E$-theory. It applies to \cs -algebras that are separable, graded, 
and admit an action of a locally compact $\gs$-compact Hausdorff group.
For such a group $G$, and for any two such \gcs{G}-algebras $A$ and $B$,
the resulting abelian group is denoted by $KE_G(A,B)$. 
The new theory recovers in a rather direct way the ordinary $K$-theory 
of ungraded \cs -algebras.
The $KE$-theory groups 
satisfy some of the good functorial properties of the other two 
bivariant theories, and there exists an associative product 
$KE_G(A,D)\tpr KE_G(D,B)\ra KE_G(A,B)$. We have also proved the
existence of two natural transformations, 
$\kke:KK_G(A,B)\ra KE_G(A,B)$ and $\kee:KE_G(A,B)\ra E_G(A,B)$, 
which preserve the product structures. Their composition 
$\kee\,\circ\,\kke$ provides an {\em explicit} construction of the map 
$KK\ra E$, abstractly known to exist because of the universality 
properties of the two theories (as we mentioned above). 
The complete analysis of the maps $\kke$ and $\kee$ will
form the substance of a subsequent paper.
The idea of constructing a theory intermediate between 
$KK$-theory and $E$-theory was suggested 
by V.~Lafforgue (private communication to N.~Higson).

Intermediate theories between $KK$-theory and $E$-theory 
appear also in the work of
J.~Cuntz \cite{Cu97}, \cite{Cu98}. Our construction is
different in initial motivation, concrete realization, and final goal:
we wanted to produce a solid framework for another proof to the
Baum-Connes conjecture for a-T-menable groups \cite{HgKas97},
\cite{HgKas01}. Details for this application will be given elsewhere.


\mic
The paper is structured as follows.
In {\bf Section 1} we briefly review the essential
definitions, theorems and constructions related to $KK$-theory. 
We also use it to set up notation. 
{\bf Section 2} constructs the new $KE$-theory. 
In subsection 2.1 we introduce and study its `cycles',
which we call {\it asymptotic Kasparov modules}.
They are appropriate families of pairs, indexed by $[1,\infty )$.
Each pair consists of a Hilbert module and
an operator on it, that are put together in a field satisfying
conditions that resemble those appearing in $KK$-theory.
An example of such cycle, motivated by the $K$-homology class of the 
Dirac operator on a spin manifold, consists of a \cs -algebra $A$,
a Hilbert space $\hs$ (constant family), a \sh -homomorphism 
$\fr{}: A \ra \bdd(\hs)$, and a family $\{ \fop{t} \}_{\ui}$
of bounded linear operators on $\hs$ satisfying:

\noindent
\makebox[2.2cm][l]{\bfseries (aKm1)}
$ F_t = F_t^* $, for all $t$; \\
\makebox[2.2cm][l]{\bfseries (aKm2)}
$ \| \, [F_t, \fr{}(a)] \,\| \xrightarrow {t\ra \infty} 0$,
for all $a\in A$; \\
\makebox[2.2cm][l]{\bfseries (aKm3)}
$ \fr{}(a)\, (F_t^2-1) \, \fr{}(a)^* \geq 0$,
\index{aaKm1@(aKm1-3)}
modulo compact operators and operators which converge
\makebox[2cm][l]{}$ $ in norm to zero.

\noindent
Such a family $\{ (\hs, \fop{t}) \}_{\ui}$ is an \akm{A}{\C}. 
Axiom (aKm2) encodes the pseudo-locality of first order elliptic
differential operators, and axiom (aKm3) is supposed to encode
the Fredholm property of elliptic operators on smooth manifolds.
The definition can be also adapted to include a group action,
and in subsection 2.2 we define, for a locally compact group $G$ and two
graded separable \gcs{G}-algebras $A$ and $B$, the group
$KE_G(A, B)$ of homotopy equivalence classes of asymptotic
Kasparov $G$-$(A,B)$-modules. Various functoriality properties
of these groups are proved in the remaining part of the section.

In {\bf Section \ref{S:KEprod}} the product in $KE$-theory is constructed
using the notions of `two-dimensional' connection
and quasi-central approximate unit.
Let $G$ be a locally compact group, 
and $A_1$, $A_2$, $B_1$, $B_2$, $D$ be \gcs{G}-algebras.
As in $KK$-theory, in its most general form, the product is a map
$$
KE_G(  A_1, B_1\tpr D ) \tpr KE_G( D\tpr A_2, B_2 ) \ra
  KE_G( A_1\tpr A_2, B_1\tpr B_2 ), \;\;\; 
    (x, y) \mapsto \Kprod{x}{D}{y}.
$$
 Insight about
the product in the new theory can be obtained by looking at the
particular case when $B_1 = B_2 = D = \C$, which corresponds to
the external product in $K$-homology. Consider two
asymptotic Kasparov modules as described above:
$\{ (\hs_1, \fop{1,t}) \}_{t}\in KE(A_1, \C)$, and
$\{ (\hs_2, \fop{2,t}) \}_{t}\in KE(A_2, \C)$. Their product is
$$\{ (\hs_1 \tpr \hs_2, 
    \fop{1,t}\tpr 1+ 1\tpr\fop{2,t}) \}_{t}\in KE(A_1\tpr A_2, \C).$$
The reader familiar with $KK$-theory will notice that no Kasparov
Technical Theorem was used in our construction. The general case is
more involved, but we hope that it is still simpler 
than in $KK$-theory. Our method is summarized in Overview 
\ref{overview}.
In subsection 3.4 we analyze the algebra behind the product. We show that
the product is associative and its various compatibilities with the
functoriality of $KE$-groups are worked out. The stability of
$KE$-theory is an easy consequence of the corresponding property
of $KK$-theory.
Subsection \ref{S:TT} plays the role of an appendix to this section.
It contains the proof of Theorem \ref{T:TT}
used to construct the product.

In the short {\bf Section 4} we present the axioms of asymptotic
Kasparov modules from the perspective of two concrete examples.
In subsection \ref{ss:ktheory} we show that the $KE$-theory groups
recover the ordinary $K$-theory for trivially graded \cs -algebras,
and in subsection \ref{ss:keGamma} we compute $KE_{\gG}(\C,\C)$,
for a discrete group $\gG$.

%
%
%
%
In {\bf Section \ref{S:kke}} we define the two natural transformations
$\kke: KK_G \ra KE_G$ and $\kee: KE_G \ra E_G$, whose composition 
gives an explicit characterization of the map from 
$KK$-theory into $E$-theory. We also briefly discuss what we were able 
to prove related to the exactness property: $KE$-theory has Puppe exact
sequences and is split-exact. Full details will appear elsewhere.

\mare
All definitions and results (theorems, propositions, corollaries, lemmas, 
important remarks and examples) are numbered in each section in order of
appearance. This convention dictates that a reference of the type
{\it m.n} sends to the result {\it n} from section {\it m}.
All equations and diagrams are numbered after the section.
The end of a proof is marked by \qedsymbol. 

\mare
\noindent
{\em Acknowledgement.}
Most of the material in this paper is part of the author's Ph.D. thesis
defended at Penn State University in 2001. I want to thank my advisor,
professor Nigel Higson, for suggesting the subject of this research, 
and I am grateful for his guidance, constant encouragement and partial
financial support during my graduate study years.
It is also a pleasure to recall here many bivariantly stimulating 
discussions with Radu Popescu.


\section{Preliminaries: review and notation}\label{C:review}

The focus of this section is to briefly review the essential
definitions, theorems and constructions related to $KK$-theory. 
We also use it to set up notation. Among the covered
topics we mention: tensor products of \cs -algebras, Hilbert modules
and tensor products of Hilbert modules, group actions, approximate 
units, and Kasparov's Technical Theorem. 
A short overview of $KK$-theory, including
its product, is given in subsection \ref{S:kk}. 

The standing assumption for the entire paper is: we
shall work in the category 
\nolinebreak$\mathbf{C^*}$\nolinebreak{\bf -\nolinebreak alg}, 
\index{$C*alg$@\cscateg}
\linebreak whose objects are the
{\it separable} and ($\Z_2$-){\it graded} \cs -algebras
\cite[14.1]{Blck}, and whose
morphisms are \sh homomorphisms that preserve the grading.

\subsection{\cs -algebras, Hilbert modules and tensor products}
\label{S:tpr} 

Given a separable graded \cs -algebra $A$, the {\em commutator}
\index{commutator}
of two elements $a,b\in A$ is: 
$[a,b] = a b - (-1)^{\der a\,\der b} b a$. %
\index{$AA.comm$@$[\,,]$}
\index{$C$@$\C$}
The \cs -algebra of complex numbers, $\C$, is trivially graded. 
As a general rule, given a locally compact space $X$, the
\cs -algebra $C_0(X)$ 
\index{$Co(X)$@$C_0(X)$}
of complex valued continuous functions on $X$ vanishing at infinity,
will be trivially graded. 

All the tensor products are graded \cite[14.4]{Blck}. 
The minimal \cs -algebra tensor product is denoted by $\tpr$, %
\index{tensor product!of \cs -algebras}
\index{$AA.tpra$@$\tprm$, $\tpr$, tensor product of \cs -algebras}
the maximal one by $\tprm$. For two \cs -algebras $A_1$ and $A_2$,
there is a transposition isomorphism $A_1\tpr A_2 \simeq A_2\tpr A_1$,
given on elementary tensors by $a_1\tpr a_2\mapsto(-1)^{\der a_1\,
\der a_2} a_2\tpr a_1$.
We recall also two of the identities that hold true with graded
tensor products:
$(a_1\tpr a_2) (b_1\tpr b_2)=
   (-1)^{\der a_2\, \der b_1} a_1 b_1 \tpr a_2 b_2$, and
$(a_1\tpr a_2) ^* = (-1)^{\der a_1 \, \der a_2}a_1^* \tpr a_2^*$,
for all $a_1, b_1 \in A_1$, and $a_2, b_2 \in A_2$.

Let $L=[1,\infty)$, %
\index{$L$@$L$, $LL$|textbf}
and $LL=[1,\infty)\times [1,\infty)$. 
For any \cs -algebra $B$
and any locally compact space $X$, the \cs -algebra $B(X)$ of $B$-valued 
continuous functions on $X$ vanishing at infinity is
$B(X) = C_0(X,B) = C_0(X)\tpr B$.
We further simplify and write:
\index{$BL$@$BL$, $BLL$, $B[0,1]$}
$BL = C_0(L,B)$, $BLL = C_0(LL,B)$, and $B[0,1] = C([0,1],B)$.

Given a Hilbert $B$-module $\fm{}$ \cite[Ch.1]{Lan}, %
\index{$E$@$\fm{}$}
the \cs -algebra of {\em adjointable operators on $\fm{}$} 
(see \cite{Kas81}, \cite{Kas80}, \cite{Lan}) is denoted by $\bdd(\fm{})$. %
\index{$Bdd(E)$@$\bdd(\fm{})$}
The closed ideal of {\em `compact' operators on $\fm{}$}
is denoted by $\cpct(\fm{})$. %
\index{$K(E)$@$\cpct(\fm{})$}
It is generated by the rank-one operators $\gth_{\xi,\eta}(\gz) =
\xi \inp{\eta}{\gz}$, for $\xi, \eta, \gz \in \fm{}$. %
\index{$xxheta.fnrk$@$\gth_{\xi,\eta}$}

Let $\ecu$ and $\ecd$ be graded Hilbert modules over $B_1$ and $B_2$,
respectively. The completion $\ecu\tpr\ecd$ of the algebraic tensor
product $\ecu\odot\ecd$ with respect to the $B_1\tpr B_2$-valued 
semi-inner product 
$\inp{\xi_1\tpr \eta_1}{\xi_2\tpr \eta_2} =
     (-1)^{\der \eta_1 (\der \xi_1+\der \xi_2)}
       \inp{\xi_1}{\xi_2} \tpr \inp{\eta_1}{\eta_2}$ 
is a Hilbert $B_1\tpr B_2$-module, called the 
{\em external tensor product of $\ecu$ and $\ecd$}.
\index{$AA.tprm$@$\tpr$, $\itpr{D}$, $\itpr{\gvf}$,
tensor product of Hilbert modules}
\index{tensor product!of Hilbert modules}
If $\gvf:B_1\ra \bdd(\ecd)$ is a \sh -homomorphism, we can also 
construct the {\em internal tensor product $\ecu\itpr{B_1}\ecd$
of $\ecu$ and $\ecd$}. (The notation $\ecu\itpr{\gvf}\ecd$
will also be used.)
It is the Hilbert $B_2$-module obtained as 
completion of the algebraic tensor product $\ecu\odot_{B_1}\ecd$ 
with respect to the $B_2$-valued semi-inner product 
$\inp{\xi_1\tpr \eta_1}{\xi_2\tpr \eta_2} =
  \inp{\eta_1}{\gvf(\inp{\xi_1}{\xi_2})(\eta_2)}$.
In both cases the grading is 
$\der(\xi\tpr\eta)=\der \xi + \der \eta$.
For details we refer the reader to \cite{Kas80}, \cite[Ch.4]{Lan}, 
or \cite[Secs.13,14]{Blck}.

Given two Hilbert modules $\ecu$ and $\ecd$, there is an embedding
$\bdd(\ecu)\tpr\bdd(\ecd)\ra\bdd(\ecu\tpr\ecd)$, given by
$(F_1\tpr F_2)(\xi\tpr\eta) = (-1)^{\der \xi \, \der F_2}
    F_1(\xi)\tpr F_2(\eta)$.
Its restriction to compact operators gives an isomorphism 
$\cpct(\ecu)\tpr\cpct(\ecd)\simeq \cpct(\ecu\tpr\ecd)$.
In the case of internal tensor product of Hilbert modules, we only
have a natural graded \sh -homomorphism $\bdd(\ecu)\ra
\bdd(\ecu\itpr{B_1}\ecd)$, $F\mapsto F\itpr{B_1} 1$,
$(F\itpr{B_1} 1)(\xi\tpr\eta) = F(\xi)\tpr \eta$.

Given a Hilbert $B$-module $\fm{}$ and a space $X$, $\fm{}(X)$ %
is the Hilbert $B(X)$-module $C_0(X)\tpr \fm{}$ (external tensor
product of Hilbert modules). We shall use the notation:
$\fm{}L = C_0(L)\tpr \fm{} = \{ \fm{}\}_t$ 
 = constant family with `fiber' $\fm{}$ indexed by $[1,\infty)$, %
\index{$EL$@$\fm{}L$, $\fm{}LL$|textbf}
$\fm{}LL = C_0(LL)\tpr \fm{}$ = constant family with `fiber' $\fm{}$
indexed by $[1,\infty)\times [1,\infty)$.

The {\em multiplier algebra}
\index{$M$@$\mlt$}
\index{multiplier algebra}
$\mlt(A)$ of a \cs -algebra $A$ is the largest \cs -algebra in which
$A$ embeds as an essential ideal (\cite[Sec.12]{Blck}, \cite[Ch.2]{Lan}). 
We recall the following two facts about multipliers:
$ \mlt ( \cpct (\fm{}) ) \simeq \bdd ( \fm{} )$,
for any Hilbert $B$-module $ \fm{} $
(\cite[13.4.1]{Blck}, \cite[2.4]{Lan}), and
$ \mlt \big( C_0([1,\infty),\cpct (\fm{})) \big) \simeq 
     C_b \big( [1,\infty),\bdd_{\scriptstyle{str}} (\fm{}) \big)$,
where $\bdd_{\scriptstyle{str}} (\fm{})$ denotes the strict topology
(\cite[3.4]{APT}).

\subsection{Group actions}
As reference for this section see \cite[Sec.1]{Kas88}.
Besides being separable and graded, the \cs -algebras that we consider
have an additional structure: the action of a group by automorphisms.
A standing assumption for the entire paper is the following:
{\em all groups are supposed to be locally
compact, $\gs$-compact and Hausdorff.}
Given such a group $G$ and a \cs -algebra
$A$, an {\it action of $G$ on $A$} is a group homomorphism
$G\ra \text{Aut}(A)$, where $\text{Aut}(A)$ is the group of 
automorphisms of $A$, with no topology on it. An element $a\in A$ is
called {\it $G$-continuous} 
\index{G-continuous@$G$-continuous element}
if the map
$G\ra A$, $g\mapsto g(a)$ is continuous.
We denote by \gcscateg  the category 
\index{$GC*alg$@\gcscateg}
with objects the separable graded \cs -algebras equipped with
$G$-action compatible with the grading
and having all the elements $G$-continuous, 
and with morphisms the equivariant \sh homomorphisms.
The objects of \gcscateg are called {\em \gcs{G}-algebras}.
The action of any group $G$ on $\C$ is trivial.

\begin{definition}\label{D:action}
Given a group $G$, a \gcs{G}-algebra $B$, and a Hilbert $B$-module
$\fm{}$, an {\it action of $G$ on $\fm{}$}, %
\index{action of a group, $G$-action!on a Hilbert module}
or a {\it $G$-action},  is an action of $G$ by grading
preserving linear automorphisms such that: (i) $G\times\fm{}\ra\fm{}$,
$(g,\xi)\mapsto g(\xi)$, is continuous in the norm topology of $\fm{}$;
(ii) $g(\xi b) = g(\xi) g(b)$; and (iii) $\inp{g(\xi)}{g(\eta)} = 
g(\inp{\xi}{\eta})$, for all $\xi, \eta \in\fm{}, b\in B, g\in G$.
We call such a Hilbert module $\fm{}$ a {\em $G$-$B$-module}.
\end{definition}

Given an action of $G$ on $\fm{}$, there is an induced action of $G$ on 
$\bdd (\fm{})$ as follows: $g (T) (\xi)= g(T(g^{-1} \xi))$, for all
$g\in G$, $T\in \bdd (\fm{})$, and $\xi\in\fm{}$. 
In this way, for any \gcs{G} -algebra $B$, there is a canonical
induced action on $\mlt(B)$.
Let $\fm{1}$ be a Hilbert $D$-module, with a $G$-action, and
$\fm{2}$ be a $G$-$(D,B)$-module. The action of $G$ on the internal
tensor product $\fm{1}\itpr{D}\fm{2}$ is given by
$g(\xi\itpr{D}\eta)=g(\xi)\itpr{D}g(\eta)$, for all $\xi\in\fm{1}$,
$\eta\in\fm{2}$. This implies, for $T\in\bdd(\fm{1})$, that
$g(\, T\itpr{D} 1\,) = g(T)\itpr{D} 1$.

The {\em standard Hilbert $G$-space} is%
\index{$H_G$@$\hs_G$}
\index{standard Hilbert $G$-space, $\hs_G$}
$\hs_G = L^2(G) \oplus L^2(G) \oplus \dots$, with infinitely many
summands, graded alternately even and odd, and equipped with the
left regular representation of $G$. Let $\cpct = \cpct (\hs_G )$%
\index{$K$@$\cpct$, $\cpct (\hs_G )$}
be the compact operators on $\hs_G$.
For any \gcs{G}-algebra $B$, the 
{\em standard Hilbert $G$-$B$-module} is 
\index{$H_B$@$\hs_B$}
\index{standard Hilbert $G$-$B$-module, $\hs_B$}
$\hs_B = l^2 \tpr L^2(G) \tpr (B\tpr B^{\text{op}})$.

\subsection{Quasi-central approximate units, and
Kasparov's Technical Theorem}\label{S:au}

Recall that a \cs -algebra is called {\em $\gs$-unital}
\index{s-unital@$\gs$-unital \cs -algebra}
if it has a countable approximate unit.

\begin{definition}\label{Def:qiqcau}
Let $G$ be a group. Consider an inclusion $I\subset B\subset A$,
where $A$ is a \gcs{G}-algebra, $B$ is a $\gs$-unital \gcs{G}-subalgebra
of $A$, and $I$ is a $\gs$-unital $G$-ideal of $A$. 
A {\em quasi-invariant quasi-central approximate unit 
for $I$ in $B$} %
\index{quasi-invariant quasi-central approximate unit}
(abbreviated q.i.q.c.a.u.) is a
continuous family $\{ u_t \}_{\ui}$ of positive, increasing, 
even elements of $I$ satisfying: \\
\makebox[1.5cm][l]{(a.u.)}
$ \|\, x u_t - x \,\| \xrightarrow{t\ra \infty} 0$, 
for all $x\in I$;\\
\makebox[1.5cm][l]{(q.c.)}
$ \|\, y u_t - u_t y \,\| \xrightarrow{t\ra \infty} 0$,
for all $y\in B$; and \\
\makebox[1.5cm][l]{(q.i.)}
$ \|\, g( u_t ) - u_t \,\| \xrightarrow{t\ra \infty} 0$, 
uniformly on compact subsets of $G$.
\end{definition}

\begin{prop}\label{P:qiqcau}
Let $G$ be a group.
A quasi-invariant quasi-central approximate unit exists for 
any closed $G$-invariant ideal $I$
of a \gcs{G}-algebra $A$.
\end{prop}

For a proof see \cite[Lemma 1.4]{Kas88}, or \cite[5.3]{GHT}.
Without a group action, the  existence of quasi-central 
approximate units is proved in \cite[3.12.14]{Pdrs}, or 
\cite[Thm.1]{Arvs}. 
Such approximate units exist for any $I\vartriangleleft A$, 
but in this paper we need a {\em countable} approximate unit
$\{u_n\}_n$ (which by interpolation gives the family $\{ u_t \}_{t}$),
and this justifies the presence of the 
separable subalgebra $B$. It is
usually clear from the context what $B$ is (the biggest subalgebra
that one needs in each particular application!), and we shall 
usually omit mention of it.

\mare

The following result in pure \cs -algebra theory is due to Kasparov.
His initial proof \cite[Sec.3]{Kas81} was complicated. 
N.~Higson \cite{Hg87b} gave an elegant proof in the non-equivariant
case based on the notion of quasi-central approximate unit. The 
statement that follows is Theorem 1.5 of \cite{Kas88}.
We shall often abbreviate the result as KTT.

\begin{ktt}
\index{Kasparov's Technical Theorem}
Let $J$ be a $\gs$-unital \gcs{G}-algebra. Assume that $E_1$ and $E_2$
are subalgebras of $\mlt(J)$, $E_1$ with $G$-action and having all
the elements $G$-continuous, such that:
(i) $\; E_1$, $E_2$ are $\gs$-unital, (ii) $\; E_1\,E_2 \subset J$.
Assume also that $\gD$ is a subset of $\mlt(J)$, separable in the
norm topology, consisting of $G$-continuous elements, and satisfying:
(iii) $\; [\gD,E_1]\subset E_1$.
Further assume that $\phi:G\ra \mlt(J)$ is a bounded function, such 
that:
(iv) $\; E_1\, \phi(G) \subset J$, $ \phi(G)\,E_1\subset J$, and
(v) $\; g\mapsto a \phi(g)$,  $g\mapsto \phi(g) a$ 
are norm continuous on $G$, for any $a\in E_1+J$.
Then there exist $G$-continuous positive even elements 
$M, N \in \mlt(J)$ with the properties:
(1) $\; M+N=1$;
(2) $\; M\,E_1\subset J$, $N\, E_2\subset J$;
(3) $\; [M,\gD]\subset J$;
(4) $\; \big(g(M)-M\big)\in J$, for all $g\in G$;
(5) $\; N\,\phi(G)\subset J$, $\phi(G)\,N \subset J$; and
(6) $\; g\mapsto N\phi(g)$, $g\mapsto \phi(g)N$ are norm 
continuous on $G$.
\end{ktt}

\subsection{$KK$-theory}\label{S:kk}

Taking into account the fact that some constructions in $KE$-theory 
are motivated by $KK$-theory constructions, and in order to have the
paper self-contained, we present in this subsection a quick review
of the theory of Gennadi Kasparov.
The $KK$-theory groups were introduced and studied in \cite{Kas75},
\cite{Kas81}, the equivariant ones under the action of a group
in \cite{Kas95}, \cite{Kas88}, under the action of a Hopf algebra
in \cite{BaaSk}, and under the action of a groupoid
in \cite{leG99}. We follow the equivariant presentation of Kasparov 
\cite{Kas88}.

\begin{definition}\label{D:Kasmod}
Consider a group $G$, and two
graded separable \gcs{G}-algebras $A$ and $B$.
A {\it Kasparov $G$-$(A,B)$-module}%
\index{Kasparov module}
is a pair $(\E, F)$, where $\E$ is a
Hilbert $B$-module, admitting a $G$-action and an action of $A$
via a \sh homomorphism $\gvf:A\ra\bdd(\E)$,
and $F\in\bdd(\E)$ is an odd $G$-continuous operator
such that for every $a\in A$ and $g\in G$
\begin{equation}\label{E:axkk}
(F-F^*)\gvf(a), \, [F,\gvf(a)], \, (F^2-1)\gvf(a), 
\text{ and } (g(F)-F)\gvf(a) \,
   \text{ all belong to } \cpct(\E).
\end{equation}
The set of all Kasparov $G$-$(A,B)$-modules
will be denoted by $kk_G(A,B)$.%
\index{$kk_G$@$kk_G(A,B)$}
A Kasparov $G$-$(A,B)$-module $( \E, F )$ is said
to be {\it degenerate} %
\index{Kasparov module!degenerate}
if for all $a\in A$ and $g\in G$:
$
(F-F^*)\gvf(a)=0, \, [F,\gvf(a)]=0, \, (F^2-1)\gvf(a)=0, 
\text{ and } (g(F)-F)\gvf(a)=0.
$
Whenever there is no risk of confusion, we shall write
$a$ instead of $\gvf(a)$.
\end{definition}

\begin{definition}
An element $(\E, F)$ of $kk_G(A, B[0,1])$ gives
by `evaluation at $s$' a family
$\{\, ( \E_s, F_s ) \in kk_G(A,B) \,|\, s\in [0,1] \,\}$,
with $\E_s=\E\itpr{\text{ev}_s} B$,
$F_s=F\itpr{\text{ev}_s} 1$. Such an element $(\E, F)$
and the family that it generates are called a {\it homotopy}%
\index{homotopy!of Kasparov modules}
between $(\E_0,F_0)$ and $(\E_1,F_1)$.
\end{definition}

\begin{definition}
The set $KK_G(A,B)$ %
\index{$KK$theory groups@$KK_G(A,B)$}
\index{KK-theory@$KK$-theory!groups}
is defined as the quotient of $kk_G(A,B)$ by
the equivalence relation generated by homotopy. 
Given an element $x=(\ec,F)\in kk_G(A,B)$, its class in $KK_G(A,B)$
will be denoted by $\kkclass{x}$.
\index{$x$@$\kkclass{x}$, $KK$-theory class}
The {\it addition} %
\index{group operation (addition)!in $KK$-theory}
of two Kasparov $G$-$(A,B)$-modules is given by the obvious notion of
direct sum.
\end{definition}

Under the above defined addition $KK_G(A,B)$ becomes an abelian group.
The following elements play an important role in the theory:
$1=1_{\C}\in KK_G(\C,\C)$,
\index{$1$@$1$, $1_{\C}$}
the class of the Kasparov module $(\C,0)$, and $1_A\in KK_G(A,A)$,
\index{$1A$@$1_A$}
the class of the Kasparov module $(A,0)$.
Given $A$, $B$, and $D$, there is a map
$\gs_D: KK_G(A,B) \ra KK_G(A\tpr D, B\tpr D), \,
  (\E,F)\mapsto (\E\tpr D, F\tpr 1)$.
\index{$xxsigma$@$\gs_D$}

\begin{definition}\label{D:cnn}
\index{connection!in $KK$-theory}
\index{$F$@$\cnn{F}$}
(\cite[Def.A.1]{CoSk}, \cite[Def.8]{Sk84})
Assume that the following elements are given:
a Hilbert $D$-module $\ecu$,
a Hilbert $(D, B)$-module $\ecd$,
and $F_2 \in\bdd( \ecd )$. Let $\ec = \ecu\itpr{D}\ecd$.
An operator $\cnn{F}\in \bdd(\ec)$ is called an
{\it $F_2$-connection for $\ecu$} if it has 
the same degree as $F_2$ and if it satisfies for every $\xi\in\ecu$:
\begin{equation}
\big(\, T_{\xi}\,F_2 - (-1)^{\der \xi\cdot \der F_2} \cnn{F}\,T_{\xi}
     \,\big)\in \cpct(\ecd,\ec), \text{ and }
\big(\, F_2\,T_{\xi}^* - (-1)^{\der \xi\cdot \der F_2}T_{\xi}^*\,\cnn{F}
     \,\big)\in \cpct(\ec,\ecd).
\end{equation}
Here $T_{\xi}\in \bdd (\ecd,\ec)$ is defined by
$T_{\xi}(\eta)=\xi \tpr \eta$, for $\eta\in\ecd$.
\end{definition}

The properties of connections are listed in \cite[Prop.9]{Sk84}. 
Here is the one that interests us most (the equivariant part is
contained in \cite[Lemma 2.6]{Kas88}):

\begin{lemma}
Consider the notation of the previous definition.
If $F_2$ satisfies, for all $d\in D$,
$[F_2, d] \in\cpct(\ecd)$,  then an $F_2$-connection $\cnn{F}$
exists for any countably generated $\ecu$.
If $d\,F_2$ and $F_2\,d$ are $G$-continuous for any $d\in D$, then
$\cnn{F}\,(K\itpr{D}1)$ and
$(K\itpr{D}1)\,\cnn{F}$ are $G$-continuous, for any
$K\in\cpct(\ecu)$.
\end{lemma}

\begin{definition}\label{D:axprod}
(\cite[Thm.A.3]{CoSk}, \cite[Def.10]{Sk84})
Let $A$, $B$, $D$ be \gcs{G}-algebras, $x=(\ecu,F_1)\in kk_G(A,D)$,
$y=(\ecd,F_2)\in kk_G(D,B)$, $\ec=\ecu\itpr{D}\ecd$. Denote by
$\Kprod{F_1}{D}{F_2}$ the set of operators $F\in\bdd(\ec)$ satisfying:
\begin{enumerate}
 \item $(\ec,F)\in kk_G(A,B)$;
 \item $F$ is an $F_2$-connection for $\ecu$; and
 \item $a\,[F_1\itpr{D}1,F]\,a^* \geq 0$, modulo $\cpct(\ec)$, for all
$a\in A$.
\end{enumerate}
For any $F\in \Kprod{F_1}{D}{F_2}$, the pair $z=(\ec,F)$ will be called
{\em the product of $x$ and $y$}.
\index{product!in $KK$-theory}
\index{$AA.prod$@$\Kprod{x}{D}{y}$}
We shall use the notation 
$\kkclass{z}=\Kprod{\kkclass{x}}{D}{\kkclass{y}}$. 
(The same notation $\sharp$ will also be used to designate 
the product 
in the new $KE$-theory. We hope that it will be clear from the 
context to what theory a certain product belongs to.)
\end{definition}

\begin{thm}\label{T:xkk}
Let $G$ be a group, and
$A$, $B$, and $D$ be separable graded \gcs{G}-algebras.
The product $\sharp_D$ exists, is unique up to homotopy, and defines
a bilinear pairing:
\begin{equation}
KK_G(A,D)\tpr KK_G(D,B)\stack{\quad\sharp_D\quad} KK_G(A,B), \;
   (\kkclass{x},\kkclass{y}) \mapsto 
      \Kprod{\kkclass{x}}{D}{\kkclass{y}}.
\end{equation}
\end{thm}

\begin{proof}
(\cite[Thm.2.11]{Kas88}, \cite[Thm.12]{Sk84})
As in the definition above, let $x=(\ecu,F_1)\in kk_G(A,D)$,
$y=(\ecd,F_2)\in kk_G(D,B)$, $\ec=\ecu\itpr{D}\ecd$.
Let $\cnn{F}$ be an $F_2$-connection for $\ecu$. Apply KTT for:
$J=\cpct(\ec)$;
$E_1=\cpct(\ecu)\itpr{D} 1 + \cpct(\ec)$;
$E_2$ = the algebra span of $[F_1\itpr{D} 1, \cnn{F}]$, $(\cnn{F} -
\cnn{F}^*)$, $(\cnn{F}^2-1)$, and $[\cnn{F},A]$;
$\gD$ = the vector space span of $F_1\itpr{D} 1$, $\cnn{F}$, and $A$;
$\phi:G\ra \bdd(\ec)$, $\phi(g)=g(\cnn{F})-\cnn{F}$.
With the elements $M$ and $N$ so obtained, define:
\begin{equation}\label{E:xkk}
F = \Mh (F_1\itpr{D} 1) + \Nh \cnn{F} .
\end{equation}
Then $F$ satisfies the conditions of Definition \ref{D:axprod}, and
consequently the class $(\ec,F)$ represents the product 
$\Kprod{\kkclass{x}}{D}{\kkclass{y}}$ of $\kkclass{x}$ and 
$\kkclass{y}$.
\end{proof}

\mic
\begin{example}[External product in $KK$-theory]\label{Ex:exprKK}%
\index{external product!in $KK$-theory|(}
Let $A_1$, $A_2$, $B_1$, $B_2$ be \gcs{G}-algebras.
The {\em external product} is the map:
\begin{equation}
KK_G(  A_1, B_1 ) \tpr KK_G( A_2, B_2 ) \stack{\quad\sharp_{\C}\quad}
  KK_G( A_1\tpr A_2, B_1\tpr B_2 ).
\end{equation}
Let $x=(\ecu,F_1)\in kk_G(A_1, B_1)$,
$y=(\ecd,F_2)\in kk_G(A_2, B_2)$, $\ec=\ecu\tpr\ecd$ (external product
of Hilbert modules). We shall still apply Kasparov's Technical Theorem,
as in the proof of the theorem above, but things are simpler because
as $F_2$-connection we can choose $\cnn{F}=1\tpr F_2$. Let:
$J=\cpct(\ec)=\cpct(\ecu)\tpr\cpct(\ecd)$;
$E_1=\cpct(\ecu)\tpr A_2 + J$;
$E_2=A_1\tpr\cpct(\fm{2})$;
$\gD$ = the vector space span of $F_1\tpr 1$, $1\tpr F_2$, and 
$A_1\tpr A_2$; and $\gvf\equiv 0$.
With the elements given by KTT, define:
\begin{equation}
F = \Mh (F_1\tpr 1) + \Nh (1\tpr F_2) .
\end{equation}
Then $F$ satisfies the conditions of Definition \ref{D:axprod}, and
consequently the class $(\ec,F)$ represents the product
$\Kprod{\kkclass{x}}{\C}{\kkclass{y}}$ of $\kkclass{x}$ and 
$\kkclass{y}$.
\index{external product!in $KK$-theory|)}
\end{example}

We make the remark that even in the external product case one
cannot in general obtain the product without the `partition of unity'
provided by KTT. See nevertheless Section 10.7 of \cite{HgRoe} for 
the example of elliptic differential operators where the `ideal'
formula $F = F_1\tpr 1 + 1\tpr F_2$ holds true. The search
for a theory in which such a simple product always exists and which
is well suited to deal with elliptic operators provided leads towards
the new $KE$-theory. See Section \ref{S:KEprod} and especially 
subsection \ref{sS:ExtProd}.

\section{$KE$-theory: definitions and functorial properties}\label{C:KEdef}

In this section we introduce the new bivariant theory.
Its cycles are appropriate families of pairs,
indexed by $[1,\infty )$. Each pair consists of a Hilbert module and
an operator on it, and they are put together in a field satisfying
conditions that resemble those appearing in $KK$-theory.
Various functoriality properties of the theory are discussed.

\mare
\subsection{Asymptotic Kasparov modules}

\begin{definition}\label{D:gab} 
Consider a group $G$, and two \gcs{G}-algebras $A$ and $B$. 
A {\it continuous field of $G$-$(A,B)$-modules}%
\index{continuous field of modules|textbf}
is a countably generated $G$-$(A,BL)$-module, that is 
a Hilbert $BL$-module $\fm{}$, admitting
a $G$-action and  
a left action of $A$ through an equivariant   
\sh -homomorphism $\fr{} :A \ra \bdd (\fm{}\,)$.
(We recall the notation:
$L=[1,\infty)$, $BL=C_0(L) \tpr B = C_0(L, B)$.)
We omit $G$ in the non-equivariant case.
\end{definition}

A continuous field $\fm{}$ of $G$-$(A,B)$-modules may be thought of as a
family $\{\et\}_{\ui}$ of Hilbert $B$-modules,
each acted on the left by $A$ and $G$, satisfying certain
continuity conditions for the left and right actions.
Indeed, for any $\ui$, let $ev_t : BL\ra B$, 
$ev_t(f\tpr b)=f(t) b$, be
the evaluation \sh homomorphism at $t$. We obtain the Hilbert
$G$-$(A,B)$-module $\et = \fm{}\itpr{ev_{t}} B$, with inner product
$\langle \fv{\xi} \tpr b, \fv{\xi}' \tpr b' \rangle_t =
    b^* \, ev_{t}(\langle \fv{\xi}, \fv{\xi}' \rangle) \, b'$.
The $A$-action an each $\et$ is $\varphi_t: A \ra \bdd(\et)$,
$\varphi_t (a) = \fr{}(a) \itpr{ev_{t}} 1$.
Whenever there is no risk of confusion, we shall write
$a$ instead of $\fr{}(a)$, and $a_t$ instead of $\fr{t}(a)$.
It is also the case that an  operator $\fop{}\in \bdd(\fm{})$ gives a
family $\{F_t\}_{\ui} = \{\fop{}\itpr{ev_{t}} 1\}_{\ui}$.
When $\fm{}=\hm L$, for a fixed Hilbert $B$-module $\hm$, the 
function $L\ra \bdd(\hm)$, $t\mapsto \fop{t}$, is `bounded and
$*$-strong continuous' \cite[3.16]{Hg90a}, \ie the family 
$\{\fop{t}\}_t$ is norm bounded, and for each $\xi\in \hm$
the functions $t\mapsto \fop{t}(\xi)$ and
$t\mapsto \fop{t}^*(\xi)$ are norm continuous.
Indeed, we have: $\bdd(\hm L) = \bdd(C_0(L)\tpr \hm)
= \mlt( \cpct(C_0(L)\tpr \hm) ) = \mlt( C_0(L,\cpct(\hm)) )
= C_b (L,\bdd_{\scriptstyle{str}} (\hm) )$, and strict continuity
is $*$-strong continuity.
On $\bdd(\hm)$ $*$-strong continuity is weaker than norm continuity.

For the remaining part of this subsection we assume no group action.
Given any Hilbert $BL$-module $\fm{}$, besides the 
{\it adjointable operators} $\bdd(\fm{})$ on $\fm{}$ and the
{\it compact operators} $\cpct(\fm{})$,%
\index{$Bdd(E)$@$\bdd(\fm{})$}
\index{$K(E)$@$\cpct(\fm{})$}
two other ideals will play an important role in our presentation:

\begin{definition}
The closed ideal of 
{\it locally compact-valued families of operators} is%
\index{$C(E)$@$\cval(\fm{})$|textbf}
\index{compact-valued families of operators|textbf}
\begin{equation}\label{E:cval}
\cval( \fm{} ) 
  = \{ \fop{}\in\bdd(\fm{}) \; \big\vert \;
     F\; f \in\cpct( \fm{} ), \text{ for all } f\in C_0(L) \}. 
\end{equation}
(Here $C_0(L)$ is viewed as a sub-\cs -algebra of $\bdd(\fm{})$ as
follows: let $\{b_n\}_n$ be an approximate unit for $B$, then
for $\xi\in\fm{}$, let $f(\xi)= \lim_{n\ra \infty} \xi (f\tpr b_n)$.)
The closed ideal of {\it vanishing families of operators} is %
\index{$J(E)$@$\loc(\fm{})$|textbf}
\index{vanishing families of operators|textbf}
\begin{equation}\label{E:loc}
 \loc ( \fm{} ) = \{ \fop{}\in\bdd(\fm{}) \; \big\vert \;
               \lim_{t\ra \infty} \| F_t \| = 0 \}.
\end{equation} 
\end{definition}

\begin{lemma}
$\cpct(\fm{}) = \cval( \fm{} ) \cap \loc ( \fm{} )$.
\end{lemma}
\begin{proof}
The inclusion $\cpct(\fm{})\subseteq \cval(\fm{})\cap\loc(\fm{})$
is clear. Let $F\in \cval(\fm{})\cap\loc(\fm{})$. From the fact that
$F\in \loc(\fm{})$ it follows that for every positive integer $n$
there exists $t_n$ such that $\| \fop{t} \| < 2^{-n}$, for all
$t> t_n$. Consider a partition of unity for $L$,
$\{\chi_0,\chi_1,\cdots, \chi_n,\cdots\}$,  subordinated to the
cover $\{[1,t_1+2^{-1})\}\cup\{(t_n,t_{n+1}+2^{-n-1})|n=1,2,\dots\}$.
Then $\fop{}=\fop{} \cdot 1=\sum_{n=0}^{\infty} \fop{}\cdot \chi_n
\in \cpct(\fm{})$, due to the fact that each term $\fop{}\cdot \chi_n$
of the sum is compact ($F\in \cval(\fm{})$), and 
of norm less than $2^{-n}$ (for $n\geq 1$).
\end{proof}

\begin{lemma}
If $\fm{}=\hm L$ is a constant family of Hilbert $B$-modules, then
any $\fop{}\in\cval(\fm{})$ generates a norm-continuous family 
of operators $\{F_t\}_t$ in $\cpct(\hm)$.
\end{lemma}
\begin{proof}
We first notice that the elements of $\cpct(\fm{})$ 
generate norm-continuous
families of operators. This is because any $\xi \in \hm L$ is 
a norm-continuous section vanishing at infinity in the constant field
of Hilbert modules $\{\hm\}_t$. Consequently the generators 
$\gth_{\xi,\eta}$, $\xi, \eta \in \hm L$, 
of $\cpct(\fm{})$ are norm-continuous.
Now, given $\fop{}\in\cval(\fm{})$, the continuity of the family
$\{F_t\}_t$ that it generates is a local property. For any $t_0$,
choose$f\in C_c(L)$, $f\equiv 1$ in a neighborhood of $t_0$. The 
definition of $\cval(\fm{})$ says that $F f\in\cpct(\fm{})$, 
and consequently $F f$ is 
a norm-continuous family. This gives the norm-continuity of 
$\{\fop{t}\}_t$ at $t_0$.
\end{proof}

\begin{goodrmk}
$\cval(\fm{})$ does not coincide with $\{ \,\fop{}\in\bdd(\fm{}) \,|\,
\fop{t}\in\cpct(\fm{t}), \text{ for all } t \,\}$.
Indeed, it is not difficult to construct a $*$-strongly continuous
family $\{ P_t\}_{\ui}$ of rank-one projections on an infinite 
dimensional Hilbert space which is not norm continuous.
\end{goodrmk}

We summarize the relations between these various ideals in
the following diagram:
\begin{equation}\label{E:ideals}
\xymatrix@1{
                       & \bdd(\fm{}) & \\
\cval( \fm{} ) \ar[ur] &             & \loc ( \fm{} ) \ar[ul] \\
                       & \cpct(\fm{}) \ar[ul] \ar[ur] & 
}
\end{equation}

\mare
\begin{definition}\label{D:akm}
Let $A$ and $B$ be graded separable \cs -algebras
(with no group action).
An {\it asymptotic Kasparov $(A,B)$-module} is a pair %
\index{asymptotic Kasparov module!non-equivariant|textbf}
$( \fm{}\,, \fop{} )$, where $\fm{}$ is a \cfm{A}{B}, 
and $\fop{} \in \bdd (\fm{}\,)$ is odd and satisfies
for any $a\in A$: \\
{\bfseries (aKm1)} \hspace{.6cm}%
$ ( \fop{} - \fop{}^* ) \fr{} (a) 
      \in \loc( \fm{} )$;\\
{\bfseries (aKm2)} \hspace{.6cm}%
\index{aaKm1@(aKm1-3)|textbf}
$ [\fop{},\fr{} (a)]\in \loc( \fm{} )$; and \\
{\bfseries (aKm3)} \hspace{.6cm}%
$ \fr{}(a) \,( \fop{}^2-1 ) \,\fr{}(a)^* \geq 0$, 
modulo $\cval( \fm{} ) + \loc( \fm{} )$.
\par\noindent
The set of all asymptotic Kasparov $(A,B)$-modules will be denoted
by $ke(A,B)$. 
\end{definition}

\begin{goodrmk}\label{R:kk-ke}
Compare these axioms with the ones that a Kasparov module $(\ec,F)$
must satisfy (Definition \ref{D:Kasmod}, (\ref{E:axkk})). 
It is worth noticing that the
third axiom of a Kasparov module, 
$(F^2-1)\gvf(a) \in \cpct(\ec)$, can be replaced 
(at least when $\| F \| \leq 1$) 
by
$\gvf(a) (F^2-1) \gvf(a)^*\geq 0$, modulo $\cpct(\ec)$, 
which looks more like our (aKm3).
\end{goodrmk}

\begin{goodrmk}
We introduce the following notation: given two operators $T, T'\in
\bdd(\fm{})$, then $T\sim T'$ if $(T-T')\in \loc(\fm{})$.
\index{$0$@$\sim$|textbf}
With this convention
(aKm1) reads $( \fop{} - \fop{}^* ) \fr{} (a)\sim 0$, and
(aKm2) reads $[\fop{},\fr{} (a)]\sim 0$, for all $a\in A$.
\end{goodrmk}

\begin{goodrmk}
In terms of families we can rephrase the conditions of Definition 
\ref{D:akm}
as follows: $\{ \et \}_{\ui}$ is a family of Hilbert $(A,B)$-modules,
$\{ F_t\}_{\ui}$ is a bounded $*$-strong continuous family of odd 
operators, meaning that for each continuous section $\xi=\{\xi_t\}_t$
the maps $t\mapsto \fop{t}(\xi_t)$ and 
$t\mapsto \fop{t}^*(\xi_t)$ are continuous sections of the field
$\{ \et \}_{\ui}$, and for each $a\in A$

\noindent
{\bfseries (aKm1$'$)} \hspace{.6cm}%
$ \|\, (F_t - F_t^*) a_t \,\| \xrightarrow{t\ra \infty} 0$; \\
{\bfseries (aKm2$'$)} \hspace{.6cm}%
\index{aaKm1a@(aKm1-3$'$)|textbf}
$ \| \, [F_t, a_t] \,\| \xrightarrow {t\ra \infty} 0$; and \\
{\bfseries (aKm3$'$)} \hspace{.6cm}%
$ a_t\, (F_t^2-1) \, a_t^* \geq K_t^a$,
for a family $ \{K_t^a\}_t \in \cval( \fm{} ) + \loc( \fm{} )$
depending on $a$.

\noindent
(Here $a_t$ denotes $\fr{t}(a) = \fr{}(a)\itpr{ev_{t}} 1$.)
\end{goodrmk}

\begin{example}\label{Ex:sh}
\index{KE-theory@$KE$-theory!class!of a \sh homomorphism}
Given a \sh homomorphism $\psi:A \ra B$,
we form the \akm{A}{B} $( \fm{}, \fop{} )$, where $\fm{}=BL$
and $\fop{}=0$. The representation of
$A$ is $\fr{} = 1\tpr\psi$.
Axioms (aKm1) and (aKm2) are trivially satisfied,
and (aKm3) follows from the fact that
$\fr{}(a) \, \fr{}(a)^* \in \cval (\fm{})$
(the family $\{K_t^a\}_t$ is constant,
$ K_t^a = - \psi(a) \, \psi(a)^* \in B \simeq \cpct(B)$,
for all $\ui$ and all $a\in A$).
More generally, given a \sh homomorphism $\psi:A \ra \cpct(\hs) \tpr B$,
with $\hs$ a countable generated Hilbert space, we form the \akm{A}{B}
$(\hs_B L, 0)$, with constant action of $A$ on `fibers' as above.
In this situation 
$ K_t^a = - \psi(a) \, \psi(a)^* \in \cpct(\hs) \tpr B
    \simeq \cpct ( \mathcal{H}_B )$. 
This simple but fundamental example implies the following principle:
if a Kasparov module $(\ec,F)$ with $F=0$ exists, then 
an asymptotic Kasparov module can be constructed from it, namely
$(\ec L, 0)$.
\end{example}

\begin{example}[The $K$-homology class of the Dirac operator]
\label{Ex:Drc}
\index{KE-theory@$KE$-theory!class!of a Dirac operator|(}
Let $M^{2n}$ be an even-\-dimensional, complete,
$\text{spin}^c$-manifold, with spinor bundle $\SSS=\SSS_M$, and Dirac
operator $D=D_M$. ($D$ is essentially self-adjoint, and
whenever functional calculus is used 
$D$ actually denotes the closure $\overline{D}=D^*$.)
The {\it fundamental asymptotic Kasparov $(C_0(M),\C)$-module} 
is constructed as follows:
$\fm{}=\{ L^2(M,\SSS) \}_{\ui}$, constant family;
the action of $C_0(M)$ is the same on each `fiber',
by multiplication operators $\fr{t}(f)=M_f$; and
$\fop{}=\{\chi (\frac1t D)\}_{\ui}$, where $\chi$ is a
{\it normalizing function} (\ie $\chi:\R\ra [-1,1]$ is odd, smooth,
and $\lim_{x\ra \pm\infty} \chi(x)=\pm 1$; for example one could take
$\chi(x) = x/{(1+x^2)^{1/2}}).$ 
Let us show that this is an asymptotic Kasparov module.
(For a thorough exposition of elliptic operators on manifolds see
\cite[Chaps.10,11]{HgRoe}. This reference also explains the 
terminology that we use in this example.)

\noindent
$\bullet$ $\fop{} \in \bdd(\fm{}\,).$
Indeed, this is implied by the norm continuity of
$t\mapsto \chi (\frac1t D).$ 

\noindent
$\bullet$ {\it $\fop{}$ satisfies} (aKm1).
As noted above, when we write $D$ we actually mean 
$\overline{D}=D^*$, which is
self-adjoint, and the functional calculus gives $F=F^*$.

\noindent
$\bullet$ {\it $\fop{}$ satisfies} (aKm2).
Let $f\in \cc(M)$. Then
$\left[\, \tfrac1t D, f\,\right] = \tfrac1t \, (D f- fD)
    = \tfrac1t \, \nabla f \xrightarrow{t\ra\infty} 0$, in norm
($\nabla f$ represents Clifford multiplication by the vector 
field $\nabla f$). This gives (\cite{Hg91}):
\begin{equation}
\begin{aligned}\label{E:akm2D}
   \left[ (\tfrac1t D \pm i )^{-1}, f \right] &=
    (\tfrac1t D \pm i )^{-1} f -
     f (\tfrac1t D \pm i )^{-1} \\
   &= (\tfrac1t D \pm i )^{-1}
    \{  f (\tfrac1t D \pm i ) -
     (\tfrac1t D \pm i ) f \}
      (\tfrac1t D \pm i )^{-1} \\
   &=  (\tfrac1t D \pm i )^{-1}
    \,(\tfrac1t \, \nabla f )\,
     (\tfrac1t D \pm i )^{-1} \\
   &\xrightarrow{t\ra\infty} 0.
\end{aligned}
\end{equation}
It follows that we obtain norm convergence 
$[\, \phi (\tfrac1t D), f \, ] \xrightarrow{t\ra\infty} 0$,
for all $\phi\in C_0(\R), f\in C_0(M)$. The significance is that
the asymptotic Kasparov module that we construct will not depend
on the normalizing  function, any two such having difference in
$C_0(\R)$. Moreover it suffices now to prove (aKm2) for {\em one}
particular normalizing function $\chi_0$. We choose it such that
its distributional Fourier transform $\widehat{\chi}_0$ 
is compactly supported, and $s\mapsto s\,\widehat{\chi}_0(s)$ is
smooth (or in $L^1(\R)$). (Such functions exist: see
\cite[10.9.3]{HgRoe}.) We know, basically from Stone's 
theorem, that:
\begin{equation}\label{E:Dirac1}
\langle \, \chi_0(D) u,v \,\rangle =
    \int_{\R} \langle \,e^{isD} u,v \,\rangle \;
      \widehat{\chi_0} (s)\; ds ,
      \text{ for all } u,v \in \cc (M,\SSS).
\end{equation}
Consider for the moment a function $f\in C^{\infty}(M)$ 
which takes values in $S^1\subset \C$ (\ie $M_f$ is an unitary operator),
and such that $\nabla f$ is also a bounded operator. 
We have:
\begin{equation}\label{E:Dirac2}
\begin{aligned}
\,[\chi_0(\tfrac1t D), f] 
 &= \chi_0(\tfrac1t D) \, f - f\, \chi_0(\tfrac1t D) 
   = f\, \big( f^{-1}\, \chi_0(\tfrac1t D) \, f - \chi_0(\tfrac1t D)
         \big) \\
 &= f\, \big( \chi_0(\tfrac1t f^{-1} D f) - \chi_0(\tfrac1t D) \big).
\end{aligned}
\end{equation}
Putting together (\ref{E:Dirac1}) and (\ref{E:Dirac2}), we obtain:
\begin{equation}\label{E:Dirac2.5}
\langle \, [\chi_0(\tfrac1t D), f] u,v \,\rangle =
  \int_{\R} \langle \,\big( e^{is t^{-1} f^{-1} D f} -
     e^{is t^{-1} D}) u , \overline{f} v \,\rangle \;
      \widehat{\chi_0} (s)\; ds.
\end{equation}
By our first computation of this paragraph, $f^{-1} D f-D=
f^{-1} [D,f]=f^{-1} \nabla f$ is a bounded operator. In accordance 
with \cite[Lemma 10.3.6]{HgRoe}, applied to $T_1=\tfrac1t f^{-1} D f$
and $T_2=\tfrac1t D$, we have:
\begin{equation}\label{E:Dirac3}
\| e^{i s T_1} - e^{i s T_2} \| \leq |s| \;\|T_1-T_2\|, \;
    \text{ for all } s\in \R.
\end{equation}
Because of (\ref{E:Dirac3}), the inner product in the integral 
of (\ref{E:Dirac2.5})
equals $|s|$ times a smooth function which is pointwise bounded by
$\tfrac1t \,\|\nabla f\| \cdot \|u\| \cdot \|v\|$. The required
norm asymptotic commutation now follows:
$$
\| \,[\chi_0(\tfrac1t D), f]  \,\|
      \leq \tfrac1t \,\|\nabla f\| \int_{\R} |s\,\widehat{\chi_0} (s)|
        \; ds.
$$
The computation made in the last part of the argument above is
\cite[Prop. 10.3.7]{HgRoe}.

Finally, any arbitrary non-zero $f\in \cc(M)$  can be written as a 
linear combination of functions on $M$ which are $S^1$-valued. Indeed,
$f=\text{Re}(f) + i \,\text{Im}(f)$, and for a  real valued 
$f\neq 0$ one writes:
$$
f = (\| f\|/2) \, \Big( \big(
f/\| f\|+ i \sqrt{1-f^2/\| f\|^2} \big) + \big(
f/\| f\|- i \sqrt{1-f^2/\| f\|^2} \big) \Big).
$$
We are through (due to the density of $\cc(M)$ in $C_0(M)$).

\noindent
$\bullet$ {\it $\fop{}$ satisfies} (aKm3).
The standard theory of elliptic first order differential operators
shows that $f\, \left(\chi^2(\tfrac1t D) -1\right)$ is compact
for $f\in C_0(M)$.
It follows that $ f\, \left(\fop{}^2 -1\right)\,\overline{f} = 0$,
modulo $\cval(\fm{})$. (The norm continuity of $t\mapsto F_t$ was used 
again here.)
\end{example}
\index{KE-theory@$KE$-theory!class!of a Dirac operator|)}

\begin{goodrmk}
Given an asymptotic Kasparov $(A,B)$-module
$(\fm{},\fop{})$ then $(\fm{},(\fop{}+\fop{}^*)/2)$ is another 
such object. Indeed, the only axiom which is not clear is (aKm3).
It reduces to showing that 
$(\fop{}+\fop{}^*)^2/4 \geq (\fop{}^2 + (\fop{}^*)^2)/2$, which in
turn is equivalent to the obvious 
$(\fop{}-\fop{}^*)\,(\fop{}-\fop{}^*)^* \geq 0$.
\end{goodrmk}

\subsection{The $KE$-theory groups}

In this subsection we define the new bivariant theory and we study
some of its functorial properties. {\em A group 
(locally compact, $\gs$-compact, Hausdorff) is assumed to act 
continuously on all the objects under study.}
We start with an extension of our previous Definition \ref{D:akm} to the
equivariant context.

\begin{definition}\label{D:akmg}
Consider a group $G$, and two
graded separable \gcs{G}-algebras $A$ and $B$.
An {\it asymptotic Kasparov $G$-$(A,B)$-module}%
\index{asymptotic Kasparov module|textbf}
is a pair $(\fm{}, \fop{})$, where $\fm{}$ is a
continuous field of $G$-$(A,B)$-modules (see Definition \ref{D:gab}),
and $\fop{}\in\bdd(\fm{})$ is an odd $G$-continuous operator that
satisfies 
{\bf (aKm1)}, {\bf (aKm2)}, {\bf (aKm3)} %
\index{aaKm1@(aKm1-3)|textbf}
of Definition \ref{D:akm}, and the extra condition:

\noindent {\bfseries (aKm4)} \hspace{.6cm}%
\index{aaKm4@(aKm4)|textbf}
$ (\, g (\fop{}) - \fop{} \,)\, \gvf(a)\, \in \loc(\fm{})$, 
for all $g\in G$, $a\in A$.
\par\noindent
In terms of families this last condition reads:
\par\noindent {\bfseries (aKm4$'$)} \hspace{.6cm}%
\index{aaKm4a@(aKm4$'$)|textbf}
$\left\| \; (g (F_t) - F_t)\, a_t \; \right\| 
  \xrightarrow{t\ra\infty} 0$, for all $g\in G$, $a\in A$,
and with $a_t=\fr{t}(a)$.
\newline
The set of all asymptotic Kasparov $G$-$(A,B)$-modules
is denoted by $ke_G(A,B)$.%
\index{$ke_G$@$ke_G(A,B)$|textbf}
\end{definition}

\begin{example}\label{Ex:splitting_m}
Consider an equivariantly split exact sequence of \gcs{G}-algebras
$$
\xymatrix@1{
 0 \ar[r] & B \ar[r]_j & D \ar[r]_p & A \ar@/_/[l]_s \ar[r] & 0
},
$$
meaning that all the \sh homomorphisms are equivariant and that
$p\circ s = id_A$. Let $\go:D\ra \mlt(B)=\bdd(B)$ be the
canonical extension of the inclusion $B\ra \mlt(B)$ (the construction
of the extension given in the proof of \cite[Prop.2.1]{Lan}
is equivariant). Let $\{u_t\}_t$ be a quasi-invariant quasi-central 
approximate unit for $B\subset \go(D) \subset \mlt(B)$ (recall 
Defini\-tion \ref{Def:qiqcau}). We associate to the above extension 
the asymptotic Kasparov $G$-$(D,B)$-module 
$
\left\{ 
  \left( 
     B\oplus B^{op},
     \left( 
         \begin{smallmatrix} 
           0 & 1-u_t \\ 1-u_t & 0
         \end{smallmatrix} 
     \right)
  \right) 
\right\}_t
$,
where the action of $D$ is constant on fibers
$\gvf_t : D \ra \bdd (B\oplus B^{op})$,
$\gvf_t(d)=\left( \begin{smallmatrix}
\go(d) & 0 \\ 0 & (\go\circ s\circ p)(d)
\end{smallmatrix} \right)$. Its class in $KE_G(D,B)$ is the 
{\em splitting morphism} of the exact sequence 
(see \cite[17.1.2b]{Blck} and \cite[Sec.5]{CoHg89}).
\end{example}

\begin{example}[The Bott element]
\label{Ex:Bott}
Let $V$ be a separable Euclidean space, and denote by $\mathcal{A}(V)$ the
non-commutative \cs -algebra used by Higson-Kasparov-Trout in their proof 
of Bott periodicity (\cite[Def.3.3]{HKT}, \cite[Def.4.1]{HgKas97}).
One considers $C_0(\R)$ graded by even and
odd functions. For a finite dimensional affine subspace $V_a$ of $V$,
denote by $V_A^0$ its linear support, and by
$\mathcal{A}(V_a)=C_0(\R) \tpr C_0(V_a,\text{Cliff}(V_a^0))$.
The \cs -algebra $\mathcal{A}(V)$ is defined as the direct limit over the 
directed set of all finite dimensional affine subspaces $V_a\subset V$
of $\mathcal{A}(V_a)$: $\mathcal{A}(V) = \varinjlim \mathcal{A}(V_a)$.
Then let $\gb:C_0(\R) \ra \mathcal{A}(V)$ be the \sh homomorphism 
given by the inclusion $(0)\subset V$, and use it to construct a family of 
\sh homomorphisms $\{ \gb_t \}_{\ui} : C_0(\R) \ra \mathcal{A}(V)$,
$\gb_t(f)=\gb(f_t)$, where $f_t(x)=f(t^{-1}x)$.
For each $t$ extend $\gb_t$ to a \sh homomorphism
$\overline{\gb_t}: C_b(\R)=\mlt(C_0(\R))\ra \mlt(\mathcal{A}(V))$. Consider
$\gl(x)=x/(1+x^2)^{1/2}$ and define 
$F_t=\overline{\gb_t}(\gl) \in \mlt(\mathcal{A}(V))$.
Further assume that a group $G$ acts isometrically and by affine
transformations on $V$.
We associate the asymptotic Kasparov $G$-$(\C,\mathcal{A}(V))$-module
$\{ (\mathcal{A}(V),F_t) \}_t$, where the action of $\C$ is constant
on fibers $\fr{t}:\C \ra \bdd(\mathcal{A}(V))$, $\fr{t}(1)=1$. We notice
that, for each $t$, $F_t$ is odd and self-adjoint (because $\gl$ has these
properties), and that $\{ F_t \}_t$ is actually a norm continuous family of
operators. This shows that (aKm1) is satisfied, (aKm2) is trivial, and
(aKm4) follows from the asymptotic equivariance of $\{\gb_t\}_t$
(\cite[Def.4.3]{HgKas97}). Finally, to see that (aKm3) holds true, note that
$F_t^2-1=-\gb_t(1/(1+x^2)) \in \mathcal{A}(V)=\cpct(\mathcal{A}(V))$.
Consequntly $F_t^2-1=0$, modulo $\cval(\mathcal{A}(V)L)$.
\end{example}

\begin{definition}
An element $(\fm{}, \fop{})$ of $ke_G(A, B[0,1])$ gives,
by `evaluation at $s$', a family 
$\{\, ( \fm{s}, \fop{s} ) \in ke_G(A,B) \,|\, s\in [0,1] \,\}$,
with $\fm{s}=\fm{}\itpr{\text{ev}_s} BL$, 
$\fop{s}=\fop{}\itpr{\text{ev}_s} 1$. Such an element $(\fm{}, \fop{})$
and the family that it generates are called a {\it homotopy}%
\index{homotopy!of asymptotic Kasparov modules|textbf}
between $(\fm{0}, \fop{0})$ and $(\fm{1}, \fop{1})$.
An {\it operator homotopy} is a homotopy %
\index{homotopy!of asymptotic Kasparov modules!operatorial}
$\{\, (\fm{}, \fop{s}) \,|\, s\in [0,1] \,\}$, with $s\mapsto \fop{s}$
being norm continuous. Note that $\fm{}$, and the action of $A$ on it,
are constant throughout an operator homotopy.
\end{definition}

\begin{example}\label{Ex:htpy}
Each $(\fm{0}, \fop{0}) = \{ ( \fm{0,t}, \fop{0,t} ) \}_t 
\in ke_G(A,B)$ is homotopic to any of its `translates'
$\{ ( \fm{0,t+N}, \fop{0,t+N} ) \}_t $.
It can also be `stretched' by a homotopy to 
$(\fm{1}, \fop{1}) = \{ ( \fm{0,h(t)}, \fop{0,h(t)} ) \}_t $,
for any increasing bijective function $h : [1,\infty)\ra [1,\infty)$. 
\end{example}

\begin{definition}
An asymptotic Kasparov $G$-$(A,B)$-module $( \fm{}, \fop{} )$ is said
to be {\it degenerate} %
\index{asymptotic Kasparov module!degenerate|textbf}
if for all $a\in A$ and $g\in G$:
$\fop{}=\fop{}^*$, $[\fop{},\fr{} (a)]=0$, 
$( g (\fop{}) - \fop{} )\gvf(a)=0$, and
$ \fr{}(a) \,(\fop{}^2-1) \,\fr{}(a)^* \geq 0$, modulo $\loc( \fm{} )$.
\end{definition}

\begin{goodrmk}
We want to comment on the definition of degenerate elements. The
first three conditions are identical with the ones for a degenerate 
Kasparov module (Definition \ref{D:Kasmod}), but in (aKm3) we require 
positivity modulo $\loc( \fm{} )$. In this way, for example, the generator
of $KE(\C,\C)$ will be described by $\cval(\hs L)/\loc(\hs L)$, 
which corresponds to the Fredholm index as invariant. This result
is required by the dimension axiom that any homology theory has to satisfy.
\end{goodrmk}

\begin{lemma}
If $(\fm{},\fop{})$ is degenerate, then it is homotopic to the 
$0$-module $(0,0)$.
\end{lemma}
\begin{proof}
The pair $(C_0([0,1))\tpr\fm{}, 1\tpr\fop{})$, with $A$ acting as
$1\tpr\fr{}$, is a degenerate 
asymptotic Kasparov $(A,BI)$-module, which gives a homotopy
between $(\fm{},\fop{})$ and $(0,0)$.
\end{proof}

\begin{definition}
Given $(\fm{},\fop{})$ and $(\fm{},\fop{}')$ in $ke_G(A,B)$,
we say that $\fop{}'$ is a {\em `small perturbation'}
\index{small perturbation@`small perturbation'}
of $\fop{}$ if $(\fop{}-\fop{}')\,\fr{}(a)\in
\loc(\fm{})$, for all $a\in A$.
\end{definition}

\begin{lemma}\label{L:smallp}
Consider $(\fm{},\fop{})$ in $ke_G(A,B)$,
and $\fop{}'$ a `small perturbation' of $\fop{}$. Then
$(\fm{},\fop{})$ and $(\fm{},\fop{}')$ are operatorially homotopic.
\end{lemma}

\begin{proof}
Indeed, the straight line segment between $\fop{}$ and $\fop{}'$ is
an operator homotopy: $\boldsymbol{\fop{}}=
\{ s \fop{} + (1-s) \fop{}' \}_{s\in [0,1]}$. We note that it is
the same proof as in $KK$-theory for `compact perturbations'
\cite[Def.17.2.4]{Blck}.
\end{proof}

\begin{cor}\label{C:self-adj}
Any $(\fm{},\fop{})\in ke_G(A,B)$ is homotopic to
$(\fm{},(\fop{}+\fop{}^*)/2)$.
\end{cor}
\begin{proof}
$(\fop{}+\fop{}^*)/2$ is a `small perturbation' of $\fop{}$.
\end{proof}

From the corollary above it follows that (aKm1) can be strengthened:
in Definitions \ref{D:akm} and \ref{D:akmg} we could consider only 
self-adjoint operators $\fop{}$. Other changes are possible too.

A less trivial example of homotopy
is provided by the next result (compare with
\cite[Lemma 11]{Sk84}). Despite the simplicity of its proof, it
will be very useful when we shall analyze in depth the product
in $KE$-theory.

\begin{lemma}\label{L:htpy}
Let $\fm{}$ be a continuous field of $G$-$(A,B)$-modules. 
Consider two asymptotic Kasparov $(A,B)$-modules 
$(\fm{}, \fop{})$, $(\fm{}, \fop{}') \in ke_G(A,B)$,
such that $ \fr{}(a)\, [F, F'] \, \fr{}(a)^*\geq 0$, 
modulo $\cval( \fm{} ) + \loc( \fm{} )$, for all $a\in A$.
Then $(\fm{}, \fop{})$ and $(\fm{}, \fop{}')$ are (operatorially)
homotopic.
\end{lemma}
\begin{proof}
Put $\fop{s} = \cos (s \,\pi/2) \fop{} + \sin (s\,\pi/2) \fop{}'$,
for $s\in [0, 1]$. Then the family $\{ ( \fm{}, \fop{s} )\}_s$ 
realizes the required homotopy.
\end{proof}

\begin{definition}
The set $KE_G(A,B)$ %
\index{$KE$theory groups@$KE_G(A,B)$|textbf}
\index{KE-theory@$KE$-theory!groups|textbf}
is defined as the quotient of $ke_G(A,B)$ by
the equivalence relation generated by homotopy. (We shall omit $G$ in the
non-equivariant case.)
Given $x=(\fm{},\fop{})\in ke_G(A,B)$, its class in $KE_G(A,B)$
will be denoted by $\keclass{x}$.
\index{$x1ke$@$\keclass{x}$, $KE$-theory class|textbf}
\index{$AA.htpclske$@$\keclass{\;}$|textbf}
The {\it addition} %
\index{asymptotic Kasparov module!addition of modules|textbf}
\index{group operation (addition)!in $KE$-theory|textbf}
of two asymptotic Kasparov $G$-$(A,B)$-modules
$(\fm{1}, \fop{1})$ and $(\fm{2}, \fop{2})$ is defined
by $(\fm{1}, \fop{1}) + (\fm{2}, \fop{2}) =
      (\fm{1}\oplus \fm{2}, \fop{1}\oplus \fop{2}) \,\in ke_G(A,B)$.
\end{definition}

\begin{thm}\label{T:KE}
With the notation of the previous definition, $KE_G(A,B)$ is an abelian
group.
\end{thm}
\begin{proof}
The argument is similar to the one for $KK$-theory --
see \cite[Prop.4]{Sk84}. The inverse of $(\fm{},\fop{})$ is
$(\fm{}^{op},-U\fop{}U^*)$, where $\fm{}^{op}$ is $\fm{}$ with the
opposite grading, $U:\fm{} \ra \fm{}^{op}$ is the identity, and $A$ 
acts on $\fm{}^{op}$ as $a (U \xi) = U((-1)^{\der a} a \xi)$.
\end{proof}

\begin{definition}\label{D:1A}
For any group $G$, $1=1_{\C}\in KE_G(\C, \C)$ is the class of 
the identity \sh homomorphism $\psi=\text{id}:\C\ra\C$,
\ie the class of $(C_0(L), 0)$, with trivial action on $C_0(L)$.
\index{$1$@$1$, $1_{\C}$|textbf}
Note that $1$ has nothing to do with the abelian group structure.
More generally, given a \gcs{G}-algebra $A$, the element $1_A\in
KE_G(A, A)$ is the class of the identity
\index{$1A$@$1_A$|textbf}
\sh homomorphism $\psi=\text{id}: A\ra A$ (as in Example 
\ref{Ex:sh}), \ie the class of $(AL,0)$.  
Given an equivariant \sh -homomorphism $\psi:A\ra B$ or more
generally $\psi:A\ra \cpct\tpr B$, its class in $KE_G(A,B)$ is denoted
by $\keclass{\psi}$.
\end{definition}

\mare
\subsection{Functoriality properties}
We discuss next some of the functoriality properties of the $KE$-groups.
They are similar to the ones that the $KK$-theory groups satisfy.

\index{KE-theory@$KE$-theory!groups!functoriality properties|(}
\mic\noindent 
(a) Given a \sh homomorphism $\psi : A_1 \ra A$, we obtain a map:
$$
\psi^* : ke_G(A,B) \ra ke_G(A_1,B)\,,\: (\fm{},\fop{})\mapsto
    (\psi^*\fm{},\fop{}).
$$
Here $\psi^*\fm{}$ denotes the same Hilbert module $\fm{}$, but
with left action by $A_1$ given by the composition $\fr{}\circ\psi:
A_1\ra \bdd(\fm{})$. We observe that $\psi^*$ respects direct sums, and
homotopy of asymptotic Kasparov modules. Consequently we get a
well-defined map, denoted by the same symbol, at the level of groups:
$\psi^* : KE_G(A,B) \ra KE_G(A_1,B)$.
It is clear that for \sh homomorphisms
$A_2 \stackrel{\go}{\lra} A_1 \stackrel{\psi}{\lra} A$ we have
$(\psi \circ \go)^* = \go^*\circ \psi^*$.

\mic\noindent 
(b) Let $\psi: B\ra B_1$ be a \sh homomorphism. Using 
$1\tpr \psi: BL \ra B_1L$, we obtain a map:
$$ 
\psi_* : ke_G(A,B) \ra ke_G(A,B_1)\,,\: (\fm{},\fop{})\mapsto
    (\fm{}\itpr{1\tpr\psi} B_1L,\fop{}\itpr{1\tpr\psi} 1). 
$$ 
This map also respects direct sums, and homotopy of
asymptotic Kasparov modules, and so gives a well-defined map:
$\psi_* : KE_G(A,B) \ra KE_G(A,B_1)$.

\mic\noindent
(c) For any \gcs{G}-algebra $D$ there is a map:
\begin{equation}\label{E:sigma}
\gs_D: ke_G(A,B) \ra ke_G(A\tpr D, B\tpr D), \;
  (\fm{},\fop{})\mapsto (\fm{}\tpr D, \fop{}\tpr 1).
\end{equation}%
\index{$xxsigma$@$\gs_D$|textbf}
It passes to quotients and gives a map 
$\gs_D: KE_G(A,B) \ra KE_G(A\tpr D, B\tpr D)$.
Indeed, we verify first that the axioms for asymptotic
Kasparov modules are satisfied.
\newline $\bullet$ {\it $\fop{}\tpr 1$ satisfies} (aKm1).
$(\fop{}\tpr 1-(\fop{}\tpr 1)^*)\,(a\tpr d) = 
  (F-F^*) a \tpr d \in \loc(\fm{})\tpr D\subseteq \loc(\fm{}\tpr D).$
\newline $\bullet$ {\it $\fop{}\tpr 1$ satisfies} (aKm2).
$(\fop{}\tpr 1)\,(a\tpr d) - (-1)^{\der a + \der d} 
  (a\tpr d)\,(\fop{}\tpr 1) = [F,a]\tpr d 
    \in \loc(\fm{})\tpr D\subseteq \loc(\fm{}\tpr D).$
\newline $\bullet$ {\it $\fop{}\tpr 1$ satisfies} (aKm3).
$(a\tpr d)\,(\fop{}^2\tpr 1) (a^*\tpr d^*) = 
  a F^2 a^* \tpr d d^* \geq a a^*\tpr d d^*$, modulo
$\cval(\fm{})\tpr D + \loc(\fm{}\tpr D) \subseteq
  \cval(\fm{}\tpr D) + \loc(\fm{}\tpr D).$
The last inclusion follows from the isomorphism
$\cpct(\mathcal{F}\tpr D) \simeq \cpct(\mathcal{F})\tpr D$,
where $\mathcal{F}$ is any Hilbert module.
\newline $\bullet$ {\it $\fop{}\tpr 1$ satisfies} (aKm4).
$( g (\fop{}\tpr 1) - \fop{}\tpr 1)\, (a\tpr d) =
  (g(\fop{}) - \fop{}) a \tpr d \in \loc(\fm{}\tpr D).$

\noindent
Finally, $\gs_D$ sends homotopic Kasparov modules to homotopic 
asymptotic Kasparov modules, and this shows that $\gs_D$ is well
defined at the level of groups.
\index{KE-theory@$KE$-theory!groups!functoriality properties|)}

\begin{prop}[Homotopy invariance]
\index{homotopy!invariance of $KE$-theory groups|textbf}
\index{KE-theory@$KE$-theory!groups!homotopy invariance|textbf}
The bifunctor $KE_G(A,B)$ is homotopy invariant in both variables:
\begin{itemize}
\item[({\it a})] let $\psi_0, \psi_1:A_1 \ra A$ 
be homotopic \sh homomorphisms; then, for any $B$,
$\psi_{0}^*=\psi_{1}^*:KE_G(A_1,B)\ra KE_G(A,B)$;
\item[({\it b})] let $\psi_0, \psi_1:B\ra B_1$ 
be homotopic \sh homomorphisms;
then, for any $A$, $\psi_{0\,*}=\psi_{1\,*}:KE_G(A,B)\ra KE_G(A,B_1)$.
\end{itemize}
\end{prop}
\begin{proof}
Once again we may follow the same proof as in $KK$-theory.

\noindent
(a) Let $\boldsymbol{\psi}:A_1\ra A[0,1]$ be a homotopy between
$\psi_0$ and $\psi_1$. If $(\fm{},\fop{})\in ke_G(A,B)$, then
$\boldsymbol{\psi}^*\big(\, \gs_{C([0,1])}((\fm{},\fop{}))\,\big)
\in ke_G(A_1, B[0,1])$ gives a homotopy between 
$\psi_{0}^*((\fm{},\fop{}))$ and $\psi_{1}^*((\fm{},\fop{}))$.

\noindent
(b) Let $\boldsymbol{\psi}:B\ra B_1[0,1]$ be a homotopy between 
$\psi_0$ and $\psi_1$. Because $ev_0$ and $ev_1$ are essential \sh
homomorphisms, it follows that $\psi_{i\,*}=ev_{i\,*}\circ
\boldsymbol{\psi}_*$, for $i=0,1$. Consequently, given $(\fm{},\fop{})
\in ke_G(A,B)$, $\boldsymbol{\psi}_*((\fm{},\fop{}))$ gives a homotopy
between $\psi_{0\,*}((\fm{},\fop{}))$ and $\psi_{1\,*}((\fm{},\fop{}))$.
\end{proof}

\subsection{Some technical results}
We conclude this section with a technical result
(namely Lemma \ref{L:constru}),
three definitions, and a `diagonalization' process, that will be used in 
the definition of the product in Section \ref{S:prod}.
Recall that any self-adjoint element $x$ of a \cs -algebra can
be written as a difference of two positive elements $x=x_{+}-x_{-}$,
with $x_{+}\, x_{-} = x_{-}\, x_{+} = 0$.
The element $x_{-}$ is called the negative part of $x$.

\begin{lemma}\label{L:constru}
Let $A$ and $B$ be separable \gcs{G}-algebras. Given
$(\fm{}, \fop{})\in ke_G(A,B)$, there exists a self-adjoint element 
$u\in \cval(\fm{})^{(0)}$ satisfying:
\begin{enumerate}
 \item[(i)]   $\; [\, u, \fop{} \,] \in \loc(\fm{})$;
 \item[(ii)]  $\; [\, u, a \,] \in \loc(\fm{})$, for all $a\in A$;
 \item[(iii)] $\; (1-u^2) \left( a\,(\,\fop{}^2 -1\,)\,a^* \right)_{-}
                         \in \loc(\fm{})$, for all $a\in A$; and
 \item[(iv)]  $(\, g(u)-u \,) \in \loc(\fm{})$, for all $g\in G$.
\end{enumerate}
\end{lemma}
\begin{proof}
Consider a dense subset $\{ a_n \}_{n=1}^{\infty}$ in $A$, and
an appropriate (see below) cover of $[1,\infty)$ by closed intervals
$\{I_n\}_{n=0}^{\infty}$, of the form $I_n =
[t_n,t_{n+2}]$, with $t_0=1$, and $\{t_n\}_n$ being a strictly
increasing sequence with $\lim_{n\ra\infty} t_n = \infty$.
Choose a partition of unity $\{ \mu_n \}_{n=0}^{\infty}$ 
in $C_0([1,\infty))$ subordinated to this cover. For each positive integer
$n$, let $r_n: BL\ra B(I_n)$ be the restriction \sh homomorphism,
and use it to define the restriction of $\fm{}$ and $\fop{}$ to $I_n$:
$\fm{}|_{I_n} = ( r_n )_* (\fm{})$, 
$\fop{}|_{I_n} = ( r_n )_* (\fop{})$.

Let $u_{0,0}$ be an arbitrary even self-adjoint element of 
$\cpct(\fm{}|_{I_0})$.
For each $n\geq 1$, apply Proposition \ref{P:qiqcau} to construct a 
quasi-invariant approximate unit 
$\{ u_{n,k} \}_{k=1}^{\infty}$ for  $\cpct( \fm{}|_{I_n} )$,
which is quasi-central for $\fop{}|_{I_n}$, $A|_{I_n}$, and
$\big\{ \left( a\,(\,\fop{}^2 -1\,)\,a^* \right)_{-}|_{I_n}
 \, \big| \, a\in A \big\}$.
There exists an index $k_n$ such that
$\|\, [u_{n,k_n}, \fop{}] \,\| < 1/n$, 
$\|\, [u_{n,k_n}, a_m] \,\| < 1/n$,
$\|\, (1-u_{n,k_n}^2) 
        \left( a_m (\fop{}^2-1) a_m^* \right)_{-} \,\| < 1/n$,
for $m=1, 2, \dots, n$. 
(For the third inequality, recall that (aKm3) implies 
$(\; a_m (\fop{}^2-1) a_m^* \;)_{-} \in \cval(\fm{}) + \loc(\fm{})$, 
with the $\cval(\fm{})$ part restricting to an element of 
$\cpct( \fm{}|_{I_n} )$, and the $\loc(\fm{})$ part having norm
$< 1/2n$ by our initial choice of the partition $\{I_n\}_{n}$.) Define:
$u = \sum_{n=0}^{\infty} \mu_n \; u_{n,k_n} \,\in\, \cval(\fm{}).$
We observe that (i) is satisfied, and that (ii) and (iii) hold true for all
the elements of the dense subset $\{ a_n \}_n $ of $A$. 
A density argument finishes the proof.
To have (iv) satisfied, 
one uses quasi-invariance, and a similar argument
after choosing a dense subset $\{ g_n \}_{n=1}^{\infty}$ of $G$.
\end{proof}

\begin{rmk}
The diagram (\ref{E:ideals}) shows that the operators that appear in (i) and
(ii) of the lemma above actually belong to $\cpct(\fm{})$.
\end{rmk}

\mic
\begin{definition}\label{D:sec}
A {\it section of }%
\index{section|textbf}
 $[1,\infty)\times [1,\infty)$ is any increasing
continuous function $h : [1,\infty) \ra [1,\infty)$, with $h(1)=1$,
$\lim_{t\ra\infty} h(t)=\infty$, differentiable on $[1,\infty)$, except
maybe for a countable set of points where it has finite one-sided 
derivatives. (Note that the differentiability assumption is just
a convenience.)
\end{definition}

\begin{lemma}\label{L:sec}
Given a countable family $\{ h_n \}_n$ of sections of 
$[1,\infty)\times [1,\infty)$, one can find a suitable strictly
increasing sequence of numbers $ \{ 1=x_0, x_1, x_2, \dots, x_n, \dots \}$, 
with $\lim_{n\ra \infty} x_n = \infty$, and a section $h$ 
satisfying the following condition: 
for each $n$, $h \geq h_i$, for $i=1, 2,\dots, n$, over the
closed interval $[x_{n-1}, x_n]$.
\end{lemma}
\begin{proof}
The definition of a section implies the existence, for each $h_i$,
of a sequence $\{ 1=x_0^i, x_1^i, \dots, x_n^i, \dots \}$, 
such that $h_i$ is differentiable on $(x_{n-1}^i, x_n^i)$, for
each positive integer $n$, and has finite one-sided derivatives at the
end points. Let $h(1)=1$, and for each integer $n\geq 1$ define:
$$
\begin{aligned} 
x_n &= \max_{1 \leq i \leq n} \{ x_n^i \},  \\
m_n &= \max_{1\leq i \leq n} \; \sup_{t\in [x_{n-1}, x_n]}
   \{ \;h_i'(t)\; \}, \text{ (one-sided derivatives included)}, \\
H_n &= \max \left\{ 0, h_{n+1}(x_n) - \left( h(x_{n-1}) + 
         m_n \, (x_n - x_{n-1}) \right) \right\} .
\end{aligned}
$$
Define $h$ on $( x_{n-1}, x_n ]$ by:
$$
h(t) = h(x_{n-1}) + \left( m_n + \frac{H_n}{ x_n - x_{n-1} } \right)
          \; (t - x_{n-1}).
$$
\end{proof}

\begin{definition}
Consider a Hilbert $BLL$-module $ \fm{} $. Given a section $h$
of $[1,\infty)\times [1,\infty)$ as in Definition \ref{D:sec},
consider the restriction \sh homomorphism:
$ \text{Res}_h : BLL \ra BL$, $f\mapsto f\vert_{\text{graph}(h)}$. 
(The parameter $t\in L$ in $BL$ is such that
$(t,h(t))\in \text{graph}(h)\subset [1,\infty)\times [1,\infty)$.)
The {\it restriction of $\fm{}$ to the graph of $h$} %
\index{restriction to the graph of a section|textbf}
is the Hilbert $BL$-module
$\fm{h} := \res{h}{\fm{}} = \fm{} \tpr_{\text{Res}_h} BL$. %
\index{$Res*$@$\text{Res}_h$, $(\text{Res}_h)_*$|textbf}
Consider now any operator $\fop{} \in \bdd(\fm{})$. The
{\it restriction of $\fop{}$ to the graph of $h$} is the operator
$\fop{h} := \res{h}{\fop{}} = \fop{} \tpr_{\text{Res}_h} 1
\in \bdd (\fm{h})$.
\end{definition}

\begin{definition}\label{D:J_2dim}
Given a Hilbert $BLL$-module $\ec$, let
$$
\loc(\ec)=\{ \; F\in\bdd(\ec) \; | 
   \lim_{t_1,t_2\ra \infty} \| F_{(t_1,t_2)} \| =0 \; \}.
$$
Here $(t_1,t_2)$ designates a point in $LL=[1,\infty)\times [1,\infty)$, and 
the limit is taken when both $t_1$ and $t_2$ approach infinity.
Note that if $F\in \loc(\ec)$ then $\fop{h}\in\loc(\fm{h})$ for any section
$h$ of $[1,\infty)\times [1,\infty)$.
\end{definition}

\section{$KE$-theory: construction of the product}\label{S:KEprod}
\label{S:prod}

In this section the product is defined and various properties, including
its associativity, are proved.

\subsection{A motivational example}\label{sS:ExtProd}

Let $G$ be a locally compact $\sigma$-compact Hausdorff group, 
and $A_1$, $A_2$, $B_1$, $B_2$, $D$ be separable \gcs{G}-algebras.
The aim is to construct a certain bilinear map
\begin{equation}\label{E:gxmap}
KE_G(  A_1, B_1\tpr D ) \tpr KE_G( D\tpr A_2, B_2 ) \ra
  KE_G( A_1\tpr A_2, B_1\tpr B_2 ).
\end{equation}
This will be the {\it product in $KE$-theory} (compare with the
product in $KK$-theory and in $E$-theory), and its construction is based
on the particular case when $B_1 = A_2 = \C$. The intuition, based on
examples coming from $K$-homology and $K$-theory, is that the product
should have the form:
\begin{equation}\label{E:gxform}
\big(\; (\,\fm{1}, \fop{1}\,)\,,\, (\,\fm{2}, \fop{2}\,) \;\big)
   \mapsto \big(\; \fm{1}\boxtimes \fm{2}  \,,\, 
      \fop{1}\boxtimes 1 + 1\boxtimes\fop{2} \;\big),
\end{equation}
where $\boxtimes$ is a certain `tensor product.' Kasparov \cite{Kas75},
\cite{Kas81} succeeded to overcome the serious technical difficulties 
that arise in making sense of (\ref{E:gxform}). We start our
approach by providing a construction of the product 
(\ref{E:gxmap}) in the case when $D=\C$, known as {\it external product}. 
By doing so, we shall present a case when the formula (\ref{E:gxform}) 
is actually correct. We shall also see the axioms (aKm1) - (aKm4) at work, 
and understand some of the difficulties involved in the general construction.

\begin{example}[External product]\label{Ex:exprKE}%
\index{external product!in $KE$-theory|(}
Consider elements $( \fm{1}, \fop{1} )\in ke_G( A_1, B_1 )$ and 
$( \fm{2}, \fop{2} )$ $\in ke_G( A_2, B_2 )$. Construct the 
$( A_1\tpr A_2, B_1L\tpr B_2L )$-module $\fm{}=\fm{1}\tpr\fm{2}$
(external tensor product of Hilbert modules),
and $\fop{}=\fop{1}\tpr 1 + 1\tpr \fop{2}\in \bdd( \fm{} )$. The
claim is that the restriction $\res{h}{( \fm{}, \fop{} )}$
to the graph of any section $h$ satisfies (aKm1)---(aKm4).
Indeed, due to the inclusions 
$\loc(\fm{1})\tpr\bdd(\fm{2})\subset\loc(\fm{})$
and $\bdd(\fm{1})\tpr\loc(\fm{2})\subset\loc(\fm{})$, 
it is easy to see that $(\fop{}-\fop{}^*) a$, $[\fop{},a]$,
$(g(\fop{})-\fop{})a \in \loc(\fm{})$, for all $a=a_1\tpr a_2 \in A$
(recall Definition \ref{D:J_2dim} for the meaning of $\loc(\fm{})$).
We also have:
$$
\begin{aligned}
(&a_1 \tpr a_2) \big( F^2-1 \big) (a_1\tpr a_2)^* \\
 &= (a_1 \tpr a_2) \big( \fop{1}^2\tpr 1 + 1\tpr \fop{2}^2 - 1
                     \big) (a_1\tpr a_2)^* \\
 &= 
\left\{ 
\begin{array}{ll}
a_1 (\fop{1}^2 -1) a_1^* \tpr a_2 a_2^* + 
 a_1 a_1^*\tpr a_2 \fop{2}^2 a_2^* \; \geq 0, &
  \text{modulo } 
  J_1=\big(\cval(\fm{1})+\loc(\fm{1})\big)\tpr\bdd(\fm{2}), \\
\text{ and } & \\
a_1 \fop{1}^2 a_1^* \tpr a_2 a_2^* + 
 a_1 a_1^*\tpr a_2 (\fop{2}^2 -1) a_2^* \; \geq 0, &
  \text{modulo } 
  J_2=\bdd(\fm{1})\tpr\big(\cval(\fm{2})+\loc(\fm{2})\big).
\end{array}
\right.
\end{aligned}
$$
Apply Lemma \ref{L:ideals}, with $J_1$, $J_2$ ideals in
$\bdd(\fm{1})\tpr\bdd(\fm{2})$, to see that $(a_1\tpr a_2) (F^2-1)
(a_1\tpr a_2)^* \geq 0$, modulo $J_1 J_2 \subseteq \cval(\fm{})
+ \loc(\fm{})$. There is only one thing left: in order
to obtain a right $(B_1\tpr B_2)L$-module (and not an
$ (B_1\tpr B_2)LL$-module as $\fm{}$ is) we restrict $\fm{}$ and
$\fop{}$ to the graph of $h(t)=t$. It is clear that $\fop{h}$
satisfies (aKm1)---(aKm4). The class of $\res{h}{(\fm{}, \fop{})}$
in $KE_G( A_1\tpr A_2, B_1\tpr B_2 )$ is called the {\it external
product} of $( \fm{1}, \fop{1} )$ and $( \fm{2}, \fop{2} )$.
Compare with Example \ref{Ex:exprKK}.
\\
{\it Conclusion.} The external
product of two asymptotic Kasparov $G$-modules
$\{ ( \fm{1,t}, \fop{1,t} ) \}_t $ and
$\{ ( \fm{2,t}, \fop{2,t} ) \}_t $ 
will be the asymptotic Kasparov $G$-module 
$\{ ( \fm{1,t}\tpr \fm{2,t}, 
         \fop{1,t} \tpr 1 + 1\tpr \fop{2,t} ) \}_t .$
\index{external product!in $KE$-theory|)}
\end{example}

In the above example we used:
\noindent
\begin{lemma}\label{L:ideals}
Let $J_1$ and $J_2$ be closed ideals of the \cs -algebra $A$. If
$a\geq 0 \; \text{mod } J_1$, and $a\geq 0 \; \text{mod } J_2$,
then $a\geq 0 \; \text{mod } J_1 J_2 = J_1 \cap J_2$.
\end{lemma}
\begin{proof}
Given a \cs -subalgebra $B$ and a closed ideal $I$ of $A$,
then $(B+I)$ is a \cs -subalgebra of $A$, and the map
$(B+I)/I \ra B/(B\cap I)$ is a \sh isomorphism
(see \cite[1.5.8]{Pdrs}).
Assume that $a\notin J_1$, otherwise interchange the roles of
$J_1$ and $J_2$ in the argument below ($a\in J_1 \cap J_2$
being trivial).
Consider B to be the \cs -subalgebra generated by $J_2$ and $a$,
and $I=J_1$. We obtain the \sh isomorphism:
$(B + J_1)/{J_1} \ra B/(B\cap J_1)=B/(J_1\cap J_2)$,
$ b + J_1 \mapsto b + J_1 J_2$. It sends the positive element
$ a + J_1$ to a positive element, namely $a + J_1 J_2$.
\end{proof}

\subsection{Two-dimensional connections}
As in Kasparov's $KK$-theory, the general product will involve
tensor products of Hilbert modules. Given a Hilbert $DL$-module
$\fm{1}$ and a Hilbert $BL$-module $\fm{2}$, their tensor product
(internal or external) will be a continuous field of modules 
over $[1,\infty) \times [1,\infty)$
(to be precise, it will be a module over the algebra
$BLL$ or $(D\tpr B)LL$). We shall call such modules
over $[1,\infty) \times [1,\infty)$, and corresponding 
families of operators, `two-dimensional.' The ones indexed
by $[1,\infty)$ are `one-dimensional.'
Our construction of the product will be based on an appropriate notion 
of connection, which is going to be a `two-dimensional' operator. 
The original definition of connection, on which ours is
modelled, appears in 
\cite[Def.A.1]{CoSk} and \cite[Def.8]{Sk84} 
(see Definition \ref{D:cnn}).

\begin{definition}\label{D:2dcnn}
Assume that the following elements are given: 
a Hilbert $DL$-module $\fm{1}$, 
a Hilbert $(D, BL)$-module $\fm{2}$,
and $\fop{2}\in\bdd( \fm{2} )$.
Consider the Hilbert $BLL$-module 
$ \fm{} = \fm{1}\itpr{DL} \fm{2}L $, 
with $\fm{2}L = C_0(L)\tpr \fm{2}$.
An operator $\cnn{F}\in \bdd(\fm{})$ is called an
{\it $\fop{2}$-connection for $\fm{1}$} if it has %
\index{connection|textbf}
\index{$F$@$\cnn{F}$|textbf}
the same degree as $\fop{2}$ and if it satisfies, 
for every {\it compactly supported}
$\fv{\xi}$ in $\fm{1}$,
$$
\bigl(\; T_{\fv{\xi}} \; (1\tpr\fop{2}) - 
     (-1)^{\der \fv{\xi} \cdot \der \fop{2}} \cnn{F} \; T_{\fv{\xi}}
   \;\bigr) \in \loc(\fm{2}L, \fm{}),
$$
and
$$
\bigl(\; (1\tpr\fop{2}) \; T^*_{\fv{\xi}} - 
     (-1)^{\der \fv{\xi} \cdot \der \fop{2}} T^*_{\fv{\xi}} \; \cnn{F} 
   \;\bigr) \in \loc(\fm{}, \fm{2}L).
$$
Here $T_{\fv{\xi}}\in \bdd (\fm{2}L, \fm{})$ is defined by 
$T_{\fv{\xi}}(g\tpr\fv{\eta})=\fv{\xi}\itpr{DL} (g\tpr\fv{\eta})$,
for $g\in C_0(L)$, and $\fv{\eta}\in\fm{2}$. Moreover
$\loc(\fm{2}L, \fm{})=\{T\in\bdd(\fm{2}L, \fm{}) \; | \; 
\lim_{t_1,t_2\ra\infty} \| T_{(t_1,t_2)} \| = 0\}$, and
$\loc(\fm{}, \fm{2}L)$ is defined similarly.
\end{definition}

\begin{rmk}
The above two conditions which a connection must satisfy are 
better remembered through the gradedly commutative modulo $\loc$ diagrams
\begin{equation}\label{E:cnn}
\begin{CD}
\fm{2}L           @>{1\tpr\fop{2}}>>       \fm{2}L       \\
@V{T_{\fv{\xi}}}VV                     @VV{T_{\fv{\xi}}}V    \\
\fm{}             @>>{\cnn{F}}>       \fm{}
\end{CD}
\text{ \quad \qquad and \quad \qquad }
\begin{CD}
\fm{2}L           @>{1\tpr\fop{2}}>>       \fm{2}L            \\
@A{T_{\fv{\xi}}^*}AA                     @AA{ T_{\fv{\xi}}^*}A    \\
\fm{}             @>>{\cnn{F}}>       \fm{}
\end{CD}
\qquad \cdot
\end{equation}
\end{rmk}

\begin{prop}\label{P:excnn}
Consider the notation of the previous definition,
with $\fr{2}:D \ra \bdd(\fm{2})$ denoting the left action of $D$ on $\fm{2}$. 
If $\fop{2}$ satisfies, for all $d\in D$,
$[\fop{2}, \fr{2} (d)] \in\loc(\fm{2})$,  then an $\fop{2}$-connection
exists for any countably generated $\fm{1}$.
\end{prop}

\begin{proof}
According with the Stabilization Theorem
\cite[Thm.2]{Kas80}, there exists an element
$V\in \bdd(\fm{1},\hs_{(DL)^{\sim}})$ of degree $0$ such that
$V^* V=1$.  (This follows from the isomorphism
$\fm{1}\oplus \hs_{(DL)^{\sim}} \simeq \hs_{(DL)^{\sim}}$.)
Assume first that the unit of $(DL)^{\sim}$ acts as identity on $\fm{2}L$.
There is then an obvious isomorphism $W: \hs_{(DL)^{\sim}}
\itpr{(DL)^{\sim}} \fm{2}L \ra \hs\tpr \fm{2}L$, given on elementary
tensors by $W( (v\tpr f)\itpr{(DL)^{\sim}} \eta ) = v\tpr f\eta$,
for $v\in\hs$, $f\in (DL)^{\sim}$, $\eta\in\fm{2}L$.
(In $ \hs\tpr \fm{2}L$ the tensor product is an external one.)
We obtain an $\fop{2}$-connection $\cnn{F}$ by imposing the
commutativity of the diagram below:
$$
\begin{CD}
\fm{1}\itpr{DL} \fm{2}L   @>{\qquad\cnn{F}\qquad}>>  \fm{1}\itpr{DL} \fm{2}L  \\
@V{V\itpr{(DL)^{\sim}} 1}VV                   
                                @AA{ V^*\itpr{(DL)^{\sim}} 1 }A     \\
\hs_{(DL)^{\sim}}\itpr{(DL)^{\sim}} \fm{2}L @.         
                  \hs_{(DL)^{\sim}}\itpr{(DL)^{\sim}} \fm{2}L
\quad ,\\
@V{W}VV                                  @AA{W^{-1}}A             \\
\hs\tpr \fm{2}L     @>>{ \quad1\tpr(1\tpr\fop{2})\quad }>  \hs\tpr \fm{2}L
\end{CD}
$$
\ie
\begin{equation}\label{E:defcnn}
\cnn{F} = (V^*\tpr 1)\,W^{-1}\,(1\tpr(1\tpr\fop{2}))\,W\,(V\tpr 1).
\end{equation}
We shall verify only one of the conditions for an $\fop{2}$-connection
(the other one being similar). Let $\fv{\xi}$ be a compactly supported 
homogeneous section of $\fm{1}$, and $V(\fv{\xi}) = 
\sum_{i=1}^{\infty} \, e_i\tpr f_i$, 
where $\{ e_i \}_{i=1}^{\infty}$ is an orthonormal basis in $\hs$, and
$\sum_{i=1}^{\infty} \, f_i^* f_i < \infty$ in $DL$. We have of course
$\der \xi = \der e_i + \der f_i$, and $\text{supp}(f_i)\subseteq 
\text{supp}(\fv{\xi})$. A direct computation gives for any 
$\fv{\eta}\in \fm{2}L$:
$$
\begin{aligned}
W\,(V\tpr 1)\, 
 &\big(\, T_{\fv{\xi}} \, (1\tpr\fop{2}) - 
        (-1)^{\der \fv{\xi} \cdot \der \fop{2}} \cnn{F} \, T_{\fv{\xi}}
        \,\big) ( \fv{\eta} )  \\
 &= W (\,V(\fv{\xi})\itpr{(DL)^{\sim}} (1\tpr\fop{2}) (\fv{\eta})\,) -
        (-1)^{\der \fv{\xi} \cdot \der \fop{2}} \,(1\tpr(1\tpr\fop{2}))\,W\,
        ( V(\fv{\xi})\itpr{(DL)^{\sim}}\fv{\eta}\,) \\
 &= \sum_{i=1}^{\infty} \, e_i\tpr f_i (1\tpr\fop{2}) (\fv{\eta}) -
        (-1)^{ \der f_i \cdot \der \fop{2}}
        \sum_{i=1}^{\infty} \, 
        e_i\tpr (1\tpr\fop{2}) ( f_i \fv{\eta}) \\
 &=  \sum_{i=1}^{\infty} 
        e_i\tpr [ f_i, 1\tpr\fop{2} ] (\fv{\eta}).
\end{aligned}
$$
%
%
Consequently, it remains to show the convergence of the
last infinite sum and that it belongs
to $\loc(\fm{2}L, \fm{})$. This is accomplished by proving
the convergence in {\it operator norm} of the partial sums
$ S_{I} = \sum_{i=1}^{I} 
        e_i\tpr [ f_i,1\tpr\fop{2} ]$,
using the expression given after the second 
equal sign in the above 
computation. The desired result follows because the partial sums 
belong to $\loc(\fm{2}L, \hs\tpr \fm{2}L)$. (The last observation uses the 
hypothesis on $\fop{2}$ and on $\xi$.)

Fix $\eps >0$. We have:
$$
\begin{aligned}
\bigl\| \bigl( S_{I+k} &- S_{I} \bigr) (\fv{\eta}) \bigr\|
 = \bigl\| \sum_{i=I+1}^{I+k} \, e_i\tpr f_i (1\tpr\fop{2}) (\fv{\eta}) -
      (-1)^{ \der f_i \cdot \der \fop{2}} \sum_{i=I+1}^{I+k} e_i\tpr
        (1\tpr\fop{2}) f_i ( \fv{\eta} ) \bigr\| \\
 &\leq \underbrace{\bigl\| 
        \sum_{i=I+1}^{I+k} \, e_i\tpr f_i (1\tpr\fop{2}) (\fv{\eta}) 
\bigr\|}_{\alpha} +
    \underbrace{\bigl\| \sum_{i=I+1}^{I+k} e_i\tpr
        (1\tpr\fop{2}) f_i ( \fv{\eta} ) \bigr\|}_{\beta} .
\end{aligned}
$$
Now:
$$
\alpha^2 
 = \bigl\| \;\langle (1\tpr\fop{2}) (\fv{\eta}),
        \bigl( \sum_{i=I+1}^{I+k} \, f_i^* f_i \bigr) (1\tpr\fop{2}) (\fv{\eta})
        \rangle \;\bigr\|  
 \leq \bigl\| \;\sum_{i=I+1}^{I+k} \, f_i^* f_i \;\bigr\| \cdot
        \bigl\| \fop{2} \bigr\|^2 \cdot \bigl\| \fv{\eta} \bigr\|^2 .
$$
Choose $I$ such that 
$ \bigl\| \,\sum_{i\in\Omega} \, f_i^* f_i \,\bigr\| 
\leq \eps^2/(4 \| F_2 \|^2)$, for every finite set $\Omega$ 
which does not intersect $\{1, 2, \dots, I\}$. Next:
$$
\begin{aligned}
\beta^2 
 &= \bigl\| \;\langle \fv{\eta}, \sum_{i=I+1}^{I+k} \, 
        f_i^* \,(1\tpr\fop{2})^*\, (1\tpr\fop{2})\, f_i \, (\fv{\eta})
        \rangle \;\bigr\| \\
 &\leq \bigl\| \fop{2}^*\, \fop{2} \bigr\| \cdot
\bigl\| \;\sum_{i=I+1}^{I+k} \, f_i^* f_i \;\bigr\| \cdot
\bigl\| \fv{\eta} \bigr\|^2 
 = \bigl\| \fop{2} \bigr\|^2 \cdot 
\bigl\| \;\sum_{i=I+1}^{I+k} \, f_i^* f_i \;\bigr\| \cdot
\bigl\| \fv{\eta} \bigr\|^2.
\end{aligned}
$$
For the chosen $I$, we obtain: $\alpha + \beta \leq (\eps/2 +\eps/2) \,
\bigl\| \fv{\eta} \bigr\|$. Consequently, 
$\bigl\| S_{I+k} - S_{I} \bigr\| \leq \eps$, for all positive 
integers $k$. This proves the norm convergence of the double sum, and
the proposition in the case when 
the unit of $(DL)^{\sim}$ acts as identity on $\fm{2}L$.

In the general case, the equation (\ref{E:defcnn}) needs replaced by:
$$
\cnn{F} = (V^*\tpr 1)\,W^{-1}\,
            \big(1\tpr (1\tpr\fop{2})|_{DL}\big) \,W\,(V\tpr 1),
$$
where $W:\hs_{DL}\itpr{DL} \fm{2}L \ra \hs\tpr(DL \cdot \fm{2}L)$, and
$(1\tpr (1\tpr\fop{2})|_{DL})\in\bdd( \hs\tpr(DL \cdot \fm{2}L))$.
We recall the definition of the restriction operator 
$ (1\tpr\fop{2})\bigr|_{DL} $ of
$1\tpr\fop{2}$ to the closed (but not necessarily complemented)
subspace $DL \cdot \fm{2}L$:
$$
(1\tpr\fop{2})\bigr|_{DL} = \sum_{n=1}^{\infty} \;
  (1\tpr\fr{2})(\delta_n^{1/2}) \; (1\tpr\fop{2}) \;
    (1\tpr\fr{2})(\delta_n^{1/2}),
$$
where $\{ u_n \}_{n=1}^{\infty}$ is an approximate unit for $DL$, 
$\delta_n = u_n -u_{n-1}$, $n=1, 2, \dots $, and $u_0=0$.
The computations are now longer, but there is no new idea
involved in the proof.
\end{proof}

The next result gathers some useful properties of connections
(compare with \cite[Prop.9]{Sk84}). The same notation as in
Definition \ref{D:2dcnn} is used.
\begin{prop}\label{P:propcnn}
(i) Let $\cnn{F}$ be an $\fop{2}$-connection for $\fm{1}$, and
$\cnn{F}'$ be an $\fop{2}'$-connection for $\fm{1}$. Then
$( \cnn{F} + \cnn{F}')$ is an $(\fop{2} + \fop{2}')$-connection 
for $\fm{1}$, and $( \cnn{F}\,\cnn{F}')$ is an 
$(\fop{2}\,\fop{2}')$-connection for $\fm{1}$. \\
(ii) The linear space of $0$-connections for $\fm{1}$ is 
$$
\big\{\; \cnn{F}\in\bdd(\fm{}) \;\big|\; (K\itpr{DL} 1) \cnn{F}, \,
           \cnn{F} (K\itpr{DL} 1) \in\loc(\fm{}), \,
   \text{ for all } K\in\cpct(\fm{1})  \;\big\}.
$$
\end{prop}
\begin{proof}
Both (i) and (ii) follow immediately from the definition of 
connection. 
\end{proof}

\begin{lemma}\label{L:constrh}
Consider the notation of Definition \ref{D:2dcnn} and assume that a 
separable set $K\subset\cval( \fm{1} )$ is given. 
Then there exists a section $h_{00}$ of $[1,\infty)\times [1,\infty)$ 
such that for any other section $h\geq h_{00}$ the following holds:
$$
\res{h}{[\, k\itpr{DL} 1, \cnn{F}\,]} \in 
  \loc\left(  \res{h}{\fm{}} \right), \text{\rm for all } k\in K.
$$
\end{lemma}
 
\begin{proof}
Choose a dense subset $\{\, k_n \,\}_{n=1}^{\infty}$ of $K$. Assume
that one is able to find for each $k_n$ a section $h_n$ such that
$\res{h}{[ k_n\itpr{DL} 1, \cnn{F} ]} \in \loc (  \res{h}{\fm{}} )$,
for any $h\geq h_n$. Apply the diagonalization process described in 
Lemma \ref{L:sec} to obtain a section $h_{00}$ which makes the conclusion 
true for all $k_n$'s. A density argument shows that the result 
holds for all $k\in K$.

Consequently it is enough to construct a section that works for a single
element $ k\in K $. As in the proof of (\ref{E:key})
in the Technical Theorem (Subsection \ref{S:TT}),
 one uses a partition of unity for $L$,
an approximation of $k\itpr{DL} 1$ by finite sums 
$\sum_{i} T_{\xi_i} T_{\eta_i}^*$, with $\xi_i$, $\eta_i \in \fm{1}$, 
and the properties of connections that $\cnn{F}$ satisfies.
\end{proof}

\subsection{Construction of the product}
We are now ready to give the construction of the product (\ref{E:gxmap})
in the case when $B_1 = A_2 = \C$. Before stating the main theorem we 
present an overview of the proof.

\mic
\begin{overview}\label{overview}
\index{product!in $KE$-theory!overview of the construction|textbf}
Consider two asymptotic Kasparov modules
$( \fm{1}, \fop{1} )\in ke_G( A, D )$ and 
$( \fm{2}, \fop{2} )\in ke_G( D, B )$.
Their product, which is an element in $ ke_G( A, B )$, is obtained
by performing the following sequence of steps.

\mic\noindent
{\bf Step 1.} Find a self-adjoint $u \in \cval ( \fm{1} )^{(0)}$ such that:
 \begin{enumerate}
  \item $[\, u, \fop{1} \,] \in \loc(\fm{1})$,
  \item $[\, u, a \,] \in \loc(\fm{1})$, for all $a\in A$,
  \item $(1-u^2) \left( a\,(\,\fop{1}^2 -1\,)\,a^* \right)_{-}
                         \in \loc(\fm{1})$, for all $a\in A$,
  \item $(g(u)-u)\in \loc(\fm{1})$, for all $g\in G$.
 \end{enumerate}

\mic\noindent
{\bf Step 2.} Define $\fm{}=\fm{1}\itpr{DL} \fm{2}L$.
Find $\cnn{F} = \cnn{F}^{\,*}$ an $\fop{2}$-connection for $\fm{1}$,
and define $\fop{} = \fop{1}\itpr{DL} 1 + (u \itpr{DL} 1)\,\cnn{F}$.
(The self-adjointness of $\cnn{F}$ is just a convenience.)

\mic\noindent
{\bf Step 3.} Choose a section $h_{00}$ of 
$[1,\infty)\times [1,\infty)$ such that the restrictions 
of the following operators to the graph of any 
other section $h\geq h_{00}$ are in $\loc ( \res{h}{\fm{}} )$:
\begin{enumerate}\setcounter{enumi}{4}
 \item $[\, u \itpr{DL} 1, \cnn{\fop{}} \,]$,
 \item $[\, u \, \fop{1} \itpr{DL} 1, \cnn{\fop{}} \,]$,
 \item $[\, u \, a \itpr{DL} 1, \cnn{\fop{}} \,]$, for all $a\in A$.
\end{enumerate}

\mic\noindent
{\bf Step 4.} Find $h_0\geq h_{00}$ such that the restriction
to the graph of any $h\geq h_0$ of:
\begin{enumerate}\setcounter{enumi}{7}
 \item $(u\itpr{DL} 1)\,\big( \cnn{F}^2 -1 \big)\,(u\itpr{DL} 1)$
 is positive modulo $\cval( \res{h}{\fm{}} )+\loc ( \res{h}{\fm{}} )$,
 \item $(u \itpr{DL} 1)\big( g(\cnn{F})-\cnn{F} \big)$ 
   is in $\loc ( \res{h}{\fm{}} )$, for all $g\in G$.
 \end{enumerate}

\mic
Once a triple $(u,\cnn{F},h_0)$ satisfying (1)--(9) is constructed,
the conclusion is that the restriction of
$( \fm{}, \fop{} )$ to the graph of any $h\geq h_0$ gives an
asymptotic Kasparov $G$-$(A,B)$-module $(\fm{h}, \fop{h})$,
that we call a {\em product of $( \fm{1}, \fop{1} )$ by
$( \fm{2}, \fop{2} )$}:
\begin{equation}\label{E:x}
\begin{aligned}
\fm{h} &= \res{h}{\fm{}} = \res{h}{\fm{1}\itpr{DL} \fm{2}L}, \\
\fop{h}&= \res{h}{\fop{}} = \res{h}{\fop{1}\itpr{DL} 1 + 
                                        (u \itpr{DL} 1)\,\cnn{F}}
        = \widetilde{\fop{1}\itpr{D,h} 1}+\widetilde{1\itpr{D,h}\fop{2}}.
\end{aligned}
\end{equation}
The notation 
\index{$F1$@$\widetilde{\fop{1}\itpr{D} 1}$}
\index{$F12$@$\widetilde{1\itpr{D}\fop{2}}$}
$\widetilde{\fop{1}\itpr{D,h} 1}=\res{h}{\fop{1}\itpr{DL} 1}$, and
$\widetilde{1\itpr{D,h}\fop{2}}=\res{h}{(u \itpr{DL} 1)\,\cnn{F}}$
is suggested by the form of the product in the external product case.
Note that in terms of families (\ref{E:x}) reads: 
\begin{equation}\label{E:xf}
(\fm{h}, \fop{h}) = 
  \left\{ \:
    \left( \,
      \fm{1,t}\itpr{D}\fm{2,h(t)},
      \fop{1,t}\itpr{D} 1 + (u_t\itpr{D} 1)\, \cnn{F}_{(t,h(t))} \,
    \right) \:
  \right\}_{\ui}.
\end{equation}
\end{overview}

\begin{goodrmk}
We do not have an axiomatic definition of the
product as in \cite[Def.10]{Sk84}, \cite[Thm.A.3]{CoSk}
(see Definition \ref{D:axprod}), so the situation is more like in
$E$-theory.
\end{goodrmk}

The following theorem guarantees that Steps 1-4 of Overview 
\ref{overview} can be performed. Its proof will be given in Subsection
\ref{S:TT}.
\begin{thm}[Technical Theorem]\label{T:TT}
\index{product!in $KE$-theory!Technical Theorem|textbf}
Let $G$ be a locally compact $\sigma$-compact Hausdorff group, and
let $A$, $B$, and $D$ be separable graded \gcs{G}-algebras.
Consider two asymptotic Kasparov modules
$( \fm{1}, \fop{1} )\in ke_G( A, D )$ and 
$( \fm{2}, \fop{2} )\in ke_G( D, B )$.
There exists a triple $(\, u, \cnn{F}, h_0 \,)$, with $u$ 
a self-adjoint element of $\cval^{(0)}( \fm{1} )$,
$\cnn{F}$ an $\fop{2}$-connection for $\fm{1}$, and $h_{0}$ 
a section of $[1,\infty)\times [1,\infty)$, as in Overview
\ref{overview}, such that for any 
other section $h\geq h_{0}$ 
$$
( \fm{h}, \fop{h} ) 
 = (\text{\rm Res}_h)_* 
  \big(\;
     \fm{1}\itpr{DL} \fm{2}L,
     \fop{1}\itpr{DL} 1 + (u \itpr{DL} 1)\,\cnn{F} \;\big)
$$
is an asymptotic Kasparov $G$-$(A,B)$-module.
\end{thm}

\mic
We can now give the definition of the {\it
product map in $KE$-theory} in the form of:
\index{product!in $KE$-theory|textbf}
\begin{thm}\label{T:x}
With the notation of the above theorem, the map
$(\, ( \fm{1}, \fop{1} ), ( \fm{2}, \fop{2} ) \,)
\mapsto ( \fm{h_0}, \fop{h_0} )$ 
 passes to quotients and defines the {\em product map}:
\begin{equation}\label{E:gmap}
KE_G(A,D) \tpr KE_G(D,B) \stack{\quad\sharp_D\quad} KE_G(A,B), \quad
  (x, y) \mapsto \Kprod{x}{D}{y}.
\end{equation}%
\index{$AA.prod$@$\Kprod{x}{D}{y}$|textbf}
\end{thm}

\begin{proof}
The notation is that of Overview \ref{overview}. 
{\em (I) Independence of h.} 
For any two $h_1,h_2 \geq h_0$ we have a homotopy 
between $(\fm{h_1},\fop{h_1})$ and $(\fm{h_2},\fop{h_2})$ given
by the explicit formula:
$$
\big\{  
  \left( \,
    \res{s h_1+ (1-s) h_2}{\fm{}},
    \res{s h_1+ (1-s) h_2}{\fop{}} \,
   \right)
\big\}_{s\in [0,1]}.
$$

\noindent
{\em (II) Independence of the triple $(\, u, \cnn{F}, h_0 \,)$.}
{\em (a)} 
As above, one can construct a homotopy between two asymptotic Kasparov
modules corresponding to different $h_0$'s satisfying 
Step \nolinebreak 4. This proves the independence of $h_0$.
{\em (b)} 
In order to show independence of $\cnn{F}$, consider two
$\fop{2}$-connections $\cnn{F}$ and $\cnn{F}'$ and the same $u$.
Now $(\cnn{F}-\cnn{F}')$ is a $0$-connection, 
and Proposition \ref{P:propcnn}(ii) implies that there exists a section
$h$ such that
$\res{h}{(u\itpr{DL} 1)\,\cnn{F} - (u\itpr{DL} 1)\,\cnn{F}'
         } \in \loc(\fm{h})$. 
Further modify $h$ such that both
$\fop{h}$ 
and
$\fop{h}'$ 
give elements in $ke_G(A,B)$.
Lemma \ref{L:smallp} applies and gives a homotopy
between $\fop{h}$ and $\fop{h}'$.
{\em (c)}
To show independence of $u$, choose two different such elements
$u$ and $u'$, both satisfying the requirements of Step 1, same
$\cnn{F}$, and an $h$ that works for both choices. We obtain a homotopy
by the formula:
$$
\big\{  
   \res{h}{
       \fop{1}\itpr{DL} 1 + \big( s (u\itpr{DL} 1) +
                                  (1-s) (u'\itpr{DL} 1) \big)\,\cnn{F} 
           }
\big\}_{s\in [0,1]}.
$$
Combining (a), (b), and (c) above we get that the homotopy class of
the element $(\fm{h},\fop{h})$ constructed in Theorem \ref{T:TT}
does not depend on the triple $(\, u, \cnn{F}, h_0 \,)$.

\noindent
{\em (III) Passage to quotients.}
Our goal is to show that the homotopy class of the product 
does not depend on the representatives in the class of 
$(\fm{1},\fop{1})$ and $(\fm{2},\fop{2})$, respectively.
Consider $(\boldsymbol{\fm{1}},\boldsymbol{\fop{1}})\in ke_G(A,D[0,1])$
a homotopy between $(\fm{1,0},\fop{1,0})$ and 
$(\fm{1,1},\fop{1,1})$. A product $(\boldsymbol{\fm{}},\boldsymbol{\fop{}})$
of $(\boldsymbol{\fm{1}},\boldsymbol{\fop{1}})$ by 
$\gs_{C[0,1]}((\fm{2},\fop{2}))$ represents a homotopy between the
product of $(\fm{1,0},\fop{1,0})$ by $(\fm{2},\fop{2})$ and a product
of $(\fm{1,1},\fop{1,1})$ by $(\fm{2},\fop{2})$. Consider now 
$(\boldsymbol{\fm{2}},\boldsymbol{\fop{2}})\in ke_G(D,B[0,1])$.
A product $(\boldsymbol{\fm{}},\boldsymbol{\fop{}})$ of
$(\fm{1},\fop{1})$ by $(\boldsymbol{\fm{2}},\boldsymbol{\fop{2}})$
represents a homotopy between the product of $(\fm{1},\fop{1})$ by
$(\fm{2,0},\fop{2,0})$ and a product of $(\fm{1},\fop{1})$ by
$(\fm{2,1},\fop{2,1})$
We obtain that the map from the statement does pass to a well-defined
map at the level of $KE$-theory groups.
\end{proof}

Using Theorem \ref{T:x} and the map $\gs$, we are now in position to
construct the general product (\ref{E:gxmap}) mentioned at the very
beginning of this section (compare with the definition in $KK$-theory
\cite[Def.2.12]{Kas88}).

\begin{definition}\label{D:gx}
Let $G$ be a group, 
and let $A_1$, $A_2$, $B_1$, $B_2$, $D$ be \gcs{G}-algebras.
The {\em general product in $KE$-theory} is the map
\begin{equation}
KE_G(  A_1, B_1\tpr D ) \tpr KE_G( D\tpr A_2, B_2 ) \ra
  KE_G( A_1\tpr A_2, B_1\tpr B_2 ), \tag{\ref{E:gxmap}} 
\end{equation}
defined by:
\begin{equation}
\Kprod{x}{D}{y}=\Kprod{\gs_{A_2}(x)}{B_1\tpr D\tpr A_2}{\gs_{B_1}(y)}.
\end{equation}
The {\em external product} corresponds to $D=\C$.
\index{external product|textbf}
\end{definition}

\mic
This subsection is concluded by showing that, in the case of external 
product, the asymptotic Kasparov module constructed in 
Example \ref{Ex:exprKE} is homotopic with the one given by the general 
product of Definition \ref{D:gx}. This will show that 
Example \ref{Ex:exprKE} really represents the construction of a 
product, and not merely of some other asymptotic Kasparov module. 
Let $x\in KE_G( A_1, B_1 )$ be represented by $( \fm{1}, \fop{1} )$, 
and $y\in KE_G( A_2, B_2 )$ be represented by $( \fm{2}, \fop{2} )$. 
According with Definition \ref{D:gx}, 
$\Kprod{x}{\C}{y} = 
   \Kprod{\gs_{A_2}(x)}{B_1\tpr A_2}{\gs_{B_1}(y)}$.
Now, $\gs_{A_2}(x)$ is represented by $(\fm{1}\tpr A_2, \fop{1}\tpr 1)$,
and $\gs_{B_1}(y)$ is represented by $(B_1\tpr\fm{2}, 1\tpr\fop{2})$.
(Bear in mind the details related to the graded tensor product of
Hilbert modules \cite[14.4.4]{Blck}.) To obtain a module that represents
the product we follow the steps given in Overview \ref{overview}.
The element $u$ of {\em Step 1} can be chosen of the form
$\{ \widetilde{u_t}\tpr \ga_{h(t)}\}_t$, with $\{ \widetilde{u_t} \}_t$
a q.i.q.c.a.u. for $\cpct(\fm{1})$, $\{ \ga_t \}_t$ an a.u. for 
$A_2$, and $h$ an arbitrary section.
In {\em Step 2} we identify $\fm{}$ with 
$\fm{1} \tpr A_2 \fm{2}$, which is a Hilbert
$(B_1 L\tpr B_2 L)$-module, acted on the left by $A_1\tpr A_2$. 
As two-dimensional connection we can take the constant field
$\{ 1\tpr \fop{2,t_2} \}_{(t_1,t_2)\in LL}$. 
With the choices and identifications
made so far, any section $h_{00}$ will do in {\em Step 3}.
In {\em Step 4} choose a section $h$ that makes the restriction to its
graph an asymptotic Kasparov module:
$$
\big( \fm{h},\fop{h} \big) =
  \big\{ \,
    \big( 
      \fm{1,t}\tpr\fm{2,h(t)},
      \fop{1,t}\tpr 1 + \widetilde{u_t}\tpr\ga_{h(t)}\fop{2,h(t)}
   \big) \, 
  \big\}_t 
\in ke_G(A_1\tpr A_2, B_1\tpr B_2).
$$
Lemma \ref{L:htpy} applies and gives a homotopy between
$(\fm{h},\fop{h})$ and
$$
\big( \fm{h}', \fop{h}' \big) =
  \big\{ \,
    \big( 
      \fm{1,t}\tpr\fm{2,h(t)},
      \fop{1,t}\tpr 1 + 1\tpr\fop{2,h(t)}
    \big) \, 
  \big\}_t 
\in ke_G(A_1\tpr A_2, B_1\tpr B_2).
$$
Finally we notice that 
$\{ ( \fm{2,h(t)}, \fop{2,h(t)} ) \}_t $ is just another 
representative of $y$, obtained by `stretching' (Example \ref{Ex:htpy})
the initial representative $(\fm{2},\fop{2})=
\{ ( \fm{2,t}, \fop{2,t} ) \}_t )$. Consequently, using two 
homotopies, we succeeded to show that the product $\sharp_{\C}$ of
Definition \ref{D:gx} is what we called external tensor 
product in Example \ref{Ex:exprKE}.

\subsection{Properties of the product}\label{S:properties}
We study in this subsection some of the properties of the product
in $KE$-theory. They are very similar with the ones that the
Kasparov product satisfies in $KK$-theory. For our first result
compare with \cite[Thm.2.14]{Kas88}.

\begin{thm}
\index{product!in $KE$-theory!properties}
The product $\sharp$ satisfies the following functoriality properties:
\begin{enumerate}
\item[(i)] it is bilinear;
\item[(ii)] it is contravariant in $A$, 
\ie $\Kprod{f^*(x)}{D}{y}=f^*(\Kprod{x}{D}{y})$,
for any \sh homorphism $f:A_1\ra A$,
$x\in KE_G(A, D)$, and $y\in KE_G(D, B)$;
\item[(iii)] it is covariant in $B$, 
\ie $g_*(\Kprod{x}{D}{y})=\Kprod{x}{D}{g_*(y)}$,
for any \sh homorphism $g:B\ra B_1$,
$x\in KE_G(A, D)$, and $y\in KE_G(D, B)$;
\item[(iv)] it is functorial in $D$, 
\ie $\Kprod{f_*( x )}{D_2}{y} = \Kprod {x}{D_1}{f^*(y)}$,
for any \sh homorphism $f: D_1\ra D_2$,
$x\in KE_G(  A, D_1 )$, and $y\in KE_G( D_2, B )$;
\item[(v)] $\sg{D_1}( \Kprod{x}{D}{y} ) = 
  \Kprod{\sg{D_1}( x )}{D\tpr D_1}{\sg{D_1}( y )}$, for 
$x\in KE_G(  A, D )$ and $y\in KE_G( D, B )$.
\end{enumerate}
\end{thm}

\begin{proof}
(i) Let $x=\keclass{(\fm{1},\fop{1})}\in KE_G(A,D)$,
$y_1=\keclass{(\fm{2},\fop{2})}$,
$y_2=\keclass{(\fm{2}',\fop{2}')}\in KE_G(D,B)$. Then:
$
\Kprod{x}{D}{y_1} = \keclass{(\text{\rm Res}_{h_1})_*
   \big((\,\fm{1}\itpr{DL}\fm{2}L,\fop{1}\itpr{DL}1+
     (u\itpr{DL}1)\cnn{F}\,)\big)}$,
$
\Kprod{x}{D}{y_2} = \keclass{(\text{\rm Res}_{h_2})_*
   \big((\,\fm{1}\itpr{DL}\fm{2}'L,\fop{1}\itpr{DL}1+
     (u\itpr{DL}1)\cnn{F}'\,)\big)}$,
$
y_1+y_2 = \keclass{(\fm{2}\oplus\fm{2}',
                     \fop{2}\oplus\fop{2}')}$.
Let $h=\text{sup}\{h_1,h_2\}$. Using $\fm{1}\itpr{DL}
(\fm{2}\oplus \fm{2}')L\simeq (\fm{1}\itpr{DL}\fm{2}L) \oplus 
   (\fm{1}\itpr{DL} \fm{2}'L)$, the definition of connection shows that
$(\cnn{F}\oplus\cnn{F}')$ is an $(\fop{2}\oplus\fop{2}')$-connection
for $\fm{1}$. It is clear that:
$$
\begin{aligned}
\Kprod{x}{D}{y_1} &+ \Kprod{x}{D}{y_2} \\
 &= \keclass{(\text{\rm Res}_{h})_*
   \big((\,\fm{1}\itpr{DL} (\fm{2}\oplus \fm{2}')L,
     \fop{1}\itpr{DL}(1\oplus 1) + (u\itpr{DL}(1\oplus 1))
         (\cnn{F}\oplus\cnn{F}'\,)\big)}    \\
 &= \Kprod{x}{D}(y_1+y_2).
\end{aligned} 
$$
The linearity in the first variable is simpler.

(ii,iii,iv) A proof using the definition of the product can be given
as for (i) above, but these properties are also a direct consequence of
the associativity of the product (Theorem \ref{T:assx}) and of the
following remark: $f^*(x)=\Kprod{\keclass{f}}{A}{x}$ and 
$g_*(y)=\Kprod{y}{B}{\keclass{g}}$.

(v) With $x=\keclass{(\fm{1},\fop{1})}$ and
$y=\keclass{(\fm{2},\fop{2})}$, $\gs_{D_1}(\Kprod{x}{D}{y})$ is 
represented by the restriction of
$
\big( (\fm{1}\itpr{DL}\fm{2}L)\tpr D_1,
      (\fop{1}\itpr{DL} 1)\tpr 1 + ((u\itpr{DL} 1) \cnn{F})\tpr 1
\big)
$
to the graph of a section $h$. Let $\fm{1}'=\fm{1}\tpr D_1$,
$\fm{2}'=\fm{2}\tpr D_1$, $D'=D\tpr D_1$.
The product $\Kprod{\gs_{D_1}(x)}{D\tpr D_1}{\gs_{D_1}(y)}$
is represented by the restriction of
$
\big( \fm{1}'\itpr{D'L}\fm{2}'L,
      (\fop{1}\tpr 1)\itpr{D'L} 1 + (\widetilde{u}\itpr{D'L} 1)\cnn{F}'
\big).
$
Under the identification $\fm{1}'\itpr{D'L}\fm{2}'L\simeq
(\fm{1}\itpr{DL}\fm{2}L)\tpr D_1$, we can take $\cnn{F}'=\cnn{F}\tpr 1$.
Given any quasi-invariant approximate unit $\widetilde{d} =
\{ d_t \}_t$ for $D_1$, we can choose $\widetilde{u}=u\tpr
\widetilde{d} \in\cval^{(0)}(\fm{1}\tpr D_1)$. Finally, after 
considering a common section for both products, Lemma \ref{L:htpy}
applies and gives a homotopy between the two representatives.
\end{proof}


\begin{rmk}
In the proof of the next theorem 
the language of elementary calculus will be used again in order
to `visualize' the construction of a double product in $KE$-theory.
A `3D-cartesian coordinate system' is assumed, with $LLL$ viewed as 
`octant' in this system. The quotations marks required by
such imprecise, but suggestive we hope, terminology will be dropped.
\end{rmk}

\begin{definition}
A {\em 3D-section} is a function $h:L\ra LL$, 
$t\mapsto (h_2(t),h_3(t))$, with $h_2$ and $h_3$ ordinary sections.
\end{definition}

\begin{thm}[Associativity of the product]\label{T:assx}
\index{associativity of the product|textbf}
\index{product!in $KE$-theory!associativity|textbf}
Let $A$, $B$, $D$, and $E$ be \gcs{G}-algebras. Then, for any
$x_1\in KE_G(A,D)$, $x_2\in KE_G(D,E)$, and $x_3\in KE_G(E,B)$,
$$
\Kprod{(\Kprod{x_1}{D}{x_2})}{E}{x_3} =
  \Kprod{x_1}{D}{(\Kprod{x_2}{E}{x_3})}.
$$
\end{thm}

\begin{proof}
Assume that $x_1$, $x_2$, $x_3$ are represented by 
$(\fm{1},\fop{1})\in ke_G(A,D)$,
$(\fm{2},\fop{2})\in ke_G(D,E)$, $(\fm{3},\fop{3})\in ke_G(E,B)$,
respectively. We shall use the notation: $\fm{12}=
\fm{1}\itpr{DL} \fm{2}L$, $\fm{23}=\fm{2}\itpr{EL} \fm{3}L$,
$\fm{}=\fm{1}\itpr{DL} \fm{2}L\itpr{ELL} \fm{3}LL$,
$x_{12,3}=\Kprod{(\Kprod{x_1}{D}{x_2})}{E}{x_3}$,
$x_{1,23}=\Kprod{x_1}{D}{(\Kprod{x_2}{E}{x_3})}$.
An inner product $( \xi\itpr{DL}\eta\itpr{ELL}\gz) \in \fm{}$ 
is abbreviated
as $( \xi\itpr{D}\eta\itpr{E}\gz )$, and similarly for operators on
$\fm{}$. In $LLL$, the first copy of $L$ and the first coordinate $t_1$
correspond to $\fm{1}$, the second copy of $L$ and the second
coordinate $t_2$ correspond to $\fm{2}$, and  the third copy of $L$
and the third coordinate $t_3$ correspond to $\fm{3}$.

We first describe the product $x_{12,3}$. As explained in the previous
subsection, $\Kprod{x_1}{D}{x_2}$ is constructed from a triple
$(u_1, \cnn{F_{12}}, h_{12})$, with $u_1\in \cval(\fm{1})$,
$\cnn{F_{12}}$ an $\fop{2}$-connection for $\fm{1}$, and $h_{12}$
a section in the $(t_1,t_2)$-plane. It is represented by 
$(\fm{12,h_{12}},\fop{12,h_{12}}) =
 \res{h_{12}}{(\fm{12},\fop{12})}$, where
$$
\fop{12} = \fop{1}\itpr{DL} 1 +(u_1\itpr{DL} 1)\cnn{F_{12}}.
$$
The product $x_{12,3}$ is constructed from a triple $(u_{12,h_{12}},
\cnn{F_{12,3}}, h_3)$, with 
$u_{12,h_{12}}\in \cval(\fm{12,h_{12}})$,
$\cnn{F_{12,3}}$ an $\fop{3}$-connection for $\fm{12,h_{12}}$, and
$h_3$ a section in the `surface' $\gS_1=\{ (t_1,t_2,t_3)\in LLL\,|\, 
t_2=h_{12}(t_1) \}$. It is represented by the restriction to the
graph of $h_3$ of $\big(\fm{12,h_{12}}\itpr{EL} \fm{3}L,
\fop{12,h_{12}}\itpr{EL} 1 +
   (u_{12,h_{12}}\itpr{EL} 1) \cnn{F_{12,3}}\big)$.
There is a simpler way of describing a representative. Define
the 3D-section $h(t)=(h_{12}(t), h_3(t))$. Consider the 
three-dimensional objects $\fm{}$ and 
$$
\fop{} = \fop{1}\itpr{D}1\itpr{E}1 + (u_1\itpr{D}1\itpr{E}1)
  (\cnn{F_{12}}\itpr{E}1) + (u_{12}\itpr{E}1) \cnn{F},
$$
with $u_1$, $\cnn{F_{12}}$ as before, $u_{12}\in \cval(\fm{12})$,
and $\cnn{F}$ a three-dimensional $\fop{3}$-connection for $\fm{12}$. 
(Such a three-dimensional connection is a straightforward
generalization of our definition for two-dimensional connection.
See (\ref{E:ass1.5}) for one of the defining, commutative up 
to $\loc$, diagrams.) The product is
represented by the restriction of $(\fm{},\fop{})$ to the graph of $h$.

Similarly, $\Kprod{x_2}{E}{x_3}$ is constructed from a triple
$(u_2, \cnn{F_{23}}, h_{23})$, with $u_2\in \cval(\fm{2})$,
$\cnn{F_{23}}$ an $\fop{3}$-connection for $\fm{2}$, and $h_{23}$
a section in the $(t_2,t_3)$-plane. It is represented by
$(\fm{23,h_{23}},\fop{23,h_{23}}) =
 \res{h_{23}}{(\fm{23},\fop{23})}$, where
$$
\fop{23} = \fop{2}\itpr{EL} 1 +(u_2\itpr{EL} 1)\cnn{F_{23}}.
$$
The product $x_{1,23}$ is constructed from a triple $(u_{1}, 
\cnn{F_{1,23}}, h_3')$, with the same $u_{1}$ as before,
$\cnn{F_{1,23}}$ an $\fop{23,h_{23}}$-connection for $\fm{1}$, and 
$h_3'$ a section in the `surface' $\gS_2=\{ (t_1,t_2,t_3)\in LLL\,|\,
t_3=h_{23}(t_2) \}$. 
Let $h'$ be the 3D-section whose graph is given by the graph of $h_3'$.
We can describe a representative for $x_{1,23}$ by the restriction
to the graph of $h'$ of 
$$
\fop{}' = \fop{1}\itpr{DL} 1 + (u_1\itpr{DL} 1)\cnn{F_{1,23}},
$$
with $\cnn{F_{1,23}}$ an $\fop{23}$-connection for $\fm{1}$.
The properties of connections given in Proposition \ref{P:propcnn}
imply that we can take $\cnn{F_{1,23}} = \cnn{F_{12}}\itpr{E}1+
\cnn{U_2} \,\cnn{F}'$, where $\cnn{U_2}$ is an 
$(u_2\itpr{EL} 1)$-connection for $\fm{1}$, and $\cnn{F}'$
is an $\cnn{F_{23}}$-connection for $\fm{1}$. The best way to see
this choice for $\cnn{F_{1,23}}$
is through the diagram below, which represents
the first of the two diagrams (\ref{E:cnn}) for the connections
under discussion (the other one being constructed in a similar way):
\begin{equation}\label{E:ass1}
\begin{CD}
(\fm{3}L)L @>{1\tpr (1\tpr \fop{3})}>> (\fm{3}L)L @.    \\
@V{f_1\tpr T_{\eta}}VV    @VV{f_1\tpr T_{\eta}}V   @. \\
(\fm{2}\itpr{EL} \fm{3}L)L 
          @>{1\tpr \cnn{F_{23}}}>> (\fm{2}\itpr{EL} \fm{3}L)L 
                   @>{1\tpr(u_2\itpr{\scriptscriptstyle{EL}}1)}>> 
                           (\fm{2}\itpr{EL} \fm{3}L)L  \\
@V{T_{\xi}}VV @VV{T_{\xi}}V @VV{T_{\xi}}V \\
\fm{}
          @>>{\cnn{F}'}> 
                          \fm{}
                   @>>{\cnn{U_2}}>
                                               \fm{}
\end{CD}
\end{equation}
(In the diagram: $f_1\in C_0(L)$, $\eta\in \fm{2}$, $\xi\in\fm{1}$.
We also have made the identification: 
$\fm{1}\itpr{DL} (\fm{2}L\itpr{ELL} \fm{3}LL) \simeq \fm{} \simeq
(\fm{1}\itpr{DL} \fm{2}L) \itpr{ELL} \fm{3}LL$.)
The bottom squares of (\ref{E:ass1}) show that $\cnn{U_2} \,\cnn{F}'$
is indeed a $(u_2\itpr{EL} 1) \cnn{F_{23}}$-connection for $\fm{1}$.
The left squares of (\ref{E:ass1}) 
are nothing but an $\fop{3}$-connection
for $\fm{12}$:
\begin{equation}\label{E:ass1.5}
\begin{CD}
(\fm{3})LL @>{(1\tpr 1)\tpr\fop{3}}>> (\fm{3})LL    \\
@V{T_{\xi\itpr{\scriptscriptstyle{DL}} (f_1\tpr\eta)}}VV 
      @VV{T_{\xi\itpr{\scriptscriptstyle{DL}} (f_1\tpr\eta)}}V  \\
\fm{} @>>{\;\;\cnn{F}'=\cnn{F}\;\;}> \fm{}
\end{CD}
\end{equation}

The outcome of all the above is the following: $x_{12,3}$ and
$x_{1,23}$ can be represented by the restriction of three
dimensional pairs $(\fm{},\fop{})$ and $(\fm{},\fop{}')$, where
\begin{equation}\label{E:ass2}
\begin{aligned}
\fop{} &= \fop{1}\itpr{D}1\itpr{E}1 + (u_1\itpr{D}1\itpr{E}1) 
  (\cnn{F_{12}}\itpr{E}1) + (u_{12}\itpr{E}1) \,\cnn{F},  \\
\fop{}'&= \fop{1}\itpr{D}1\itpr{E}1 + (u_1\itpr{D}1\itpr{E}1)
  (\cnn{F_{12}}\itpr{E}1) + (u_1\itpr{D}1\itpr{E}1)\,\cnn{U_2}\,\cnn{F},
\end{aligned}
\end{equation}
to the graphs of appropriate sections $h$ and $h'$, respectively. 
We complete the 
proof by showing that $h$ and $h'$ can be chosen the same, and that
$\fop{}$ and $\fop{}'$ are homotopic.

The proof of Technical Theorem given in Subsection \ref{S:TT}
(see also the remark that follows that proof) shows that, while
the section $h_0$ that appears in the triple $(u,\cnn{F},h_0)$ used
to define the product of two $KE$-modules 
is an important element, the `right decay
conditions' actually hold true on a two dimensional object, namely
over $\cup_{n=0}^{\infty} [T_{1,n},T_{1,n+1}]\times[T_{2,n},\infty)$,
or over $\{\, (t_1,t_2)\in LL\,|\, t_2\geq h_0(t_1)\,\}$.
(Notation as in the proof of Technical Theorem.)
This implies that in the computation of a product 
the section is important only through the fact that it captures 
the behavior when both $t_1\ra \infty$ and $t_2\ra \infty$. 
This observation is summarized as:

\begin{lemma}\label{L:ass1}
The products $\Kprod{(\Kprod{x_1}{D}{x_2})}{E}{x_3}$ and 
$\Kprod{x_1}{D}{(\Kprod{x_2}{E}{x_3})}$ can be computed by restricting
the operators of (\ref{E:ass2}) to a {\em common} 3D-section $h$.
\end{lemma}

We need one more result:

\begin{lemma}
Define:
$
\loc_0(\fm{}) = \{ \fop{}\in\bdd(\fm{}) \; \big\vert \;
    \lim_{t_1,t_2,t_3\ra \infty} \| F_{(t_1,t_2,t_3)} \| = 0 \}.
$
(Here $t_1,t_2,t_3\ra \infty$ means $t_i\ra \infty$, for $i=1,2,3$.)
Then
$[\, u_1\itpr{D}1\itpr{E}1,\cnn{U_2}\,]\in \loc_0(\fm{})$, and
$u_{12}$ can be chosen such that 
$[\, u_{12}\itpr{E}1, (u_1\itpr{D}1\itpr{E}1)\,\cnn{U_2}\,]\in 
\loc_0(\fm{})$.
\end{lemma}
\begin{proof}
Modulo an element in $\loc(\fm{1})\itpr{D}1\itpr{E}1\subset 
\loc_0(\fm{})$, $(u_1\itpr{D}1\itpr{E}1)$ can be approximated on 
compact intervals in $t_1$-variable by finite sums 
$\sum_i (T_{\xi_i} T_{\eta_i}^*\itpr{E}1)$,
with $\xi_i, \eta_i\in \fm{1}$, compactly supported. (See the 
proof of Technical Theorem in Subsection \ref{S:TT}.) This implies:
$$
\begin{aligned}
&[ u_1\itpr{D} 1\itpr{E} 1, \cnn{U_2} ] \\
  &\qquad\sim \sum_i \big( (T_{\xi_i} T_{\eta_i}^*\itpr{E}1) \cnn{U_2}
                 - \cnn{U_2} (T_{\xi_i} T_{\eta_i}^*\itpr{E}1) \big) 
   && \text{ modulo } \loc_0(\fm{}) \\
  &\qquad\sim (-1)^{\der \eta_i} \sum_i \big( 
       T_{\xi_i} (1\tpr(u_2\itpr{EL}1)) T_{\eta_i}^* 
        - T_{\xi_i} (1\tpr(u_2\itpr{EL}1)) T_{\eta_i}^*\big),
   && \text{ modulo } \loc_0(\fm{}) \\
  &\qquad=0.
\end{aligned}
$$
This proves the first inclusion. For the second one, use the same
approximation for $(u_1\itpr{D}1\itpr{E}1)$ as above to see that,
modulo $\loc_0(\fm{})$, $(u_1\itpr{D}1\itpr{E}1)\,\cnn{U_2}$
is an element of $\bdd(\fm{12})\itpr{E}1$. The claimed
asymptotic-commutativity follows by actually imposing it as an 
{\em extra
requirement} for $u_{12}$ (besides the conditions that appear
in Step 1, Overview \ref{overview}).
\end{proof}

\noindent
This last lemma implies that $a\,[\fop{},\fop{}']\,a^* \geq 0$, modulo
$\loc(\fm{h})$, for any section $h$, and consequently Lemma 
\ref{L:htpy} gives the required homotopy.
We have showed that $x_{12,3}=x_{1,23}$ in $KE_G(A,B)$, and this
completes the proof of Theorem \ref{T:assx}.
\end{proof}

\mic
\begin{rmk}
There is another way to see the homotopy between the operators
from (\ref{E:ass2}). It uses the following result, whose proof
is left to the reader:
\begin{lemma}
$(u_1\itpr{D}1\itpr{E}1)\,\cnn{U_2}$ satisfies the
(properly modified) conditions of
Step 1, Overview \ref{overview}, that $(u_{12}\itpr{E}1)$ satisfies.
\end{lemma}
\noindent
Consequently, the straight line homotopy
$
\{ \, (1-s)(u_{12}\itpr{E}1) + s (u_1\itpr{D}1\itpr{E}1)\,\cnn{U_2}
  \, \}_{s\in[0,1]}
$
can be used to give a homotopy between $\fop{}$ and $\fop{}'$.
\end{rmk}

\mic
Recall from Definition \ref{D:1A} that $1=1_{\C}\in KE_G(\C,\C)$
is the class of the identity homomorphism $\psi=\text{id}:\C\ra\C$.
For the next result compare with \cite[Thm.4.5]{Kas81},
\cite[Prop.17]{Sk84}.

\begin{prop}\label{T:1}
Let $A$ and $B$ be separable \gcs{G}-algebras, then
\index{product!in $KE$-theory!existence of unit|textbf}
$$
\Kprod{1_{\C}}{\C}{x}=\Kprod{x}{\C}{1_{\C}}=x, \text{ for any }
x\in KE_G(A,B).
$$
\end{prop}

\begin{proof}
One equality is easy. We have:
$\Kprod{x}{\C}{1_{\C}}\stackrel{\text{def}}{=}
   \Kprod{x}{B}{\gs_B(1_{\C})} =
    \Kprod{x}{B}{1_B}.$
Let $x$ be represented by $(\fm{},\fop{})$ and $1_B$ be represented by
$(BL,0)$. As $0$-connection for $\fm{}$ we can take the $0$ operator,
and we can restrict to $h(t)=t$ in the construction of the product
to obtain:
$$
\begin{aligned}
 & \res{h}{\fm{}\itpr{BL} BLL} \simeq \fm{}, \text{ via }
     \res{h}{\xi\itpr{BL}(f\tpr g\tpr b)} \mapsto (\xi\cdot (fg\tpr b)),
  \\
 & \res{h}{\fop{}\itpr{BL} 1}=\fop{} 
     \text{ (under the previous isomorphism).}
\end{aligned}
$$ 
Consequently $\Kprod{x}{B}{1_B}=x$.

For the other equality, we start with:
$\Kprod{1_{\C}}{\C}{x}\stackrel{\text{def}}{=}
   \Kprod{\gs_A(1_{\C})}{A}{x} =
    \Kprod{1_A}{A}{x}.$
Let $1_A$ be represented by $(AL,0)$. 
Consider a quasi-invariant approximate unit $\{u_n\}_{n=1}^{\infty}$
for $A$, and construct an element $u=\{u_t\}_{\ui}\in \cval^{(0)}(AL)$ 
by interpolating the $u_n$'s:
$u_t = (1-\{t\}) u_{[t]} + \{t\} u_{[t]+1}$, with $[t]$ denoting the
greatest integer smaller that t, and $\{t\}=t-[t]$. 
We shall exhibit a homotopy between $(\fm{},\fop{})$ and a
representative of the product $\Kprod{1_A}{A}{x}$ constructed 
using $u$.

If $A$ is unital, consider the projection $\fr{}(1)=P\in\bdd(\fm{})$.
With the identification $AL\itpr{AL}\fm{}L \simeq (P\fm{})L$, 
and after choosing 
$u\equiv 1$ and $h(t)=t$ in the definition of the product, we
obtain as representative of $\Kprod{1_A}{A}{x}$ the
asymptotic Kasparov module $(P\fm{},P\fop{})$. There is an operator
homotopy between $(\fm{},\fop{})$ and $(\fm{},P\fop{})=
(P\fm{},P\fop{})\oplus((1-P)\fm{},0)$, with the second summand being
degenerate. This proves that $(\fm{},\fop{})$ represents the product.

Assume now that $A$ is not unital.
Let $A^{\sim}$ be the unitization of $A$, with $1$ acting as
identity on $\fm{}$. 
Following \cite[Prop.17]{Sk84}, let 
$\widehat{A[0,1]}$ be the $G$-$(A,A^{\sim}[0,1])$-module:
$$
\widehat{A[0,1]} =\{\, f:[0,1] \ra A^{\sim} \,|\, f(1)\in A \,\}
  \subseteq A^{\sim}[0,1].
$$
Notice that $A$ acts as multiplication by constant functions.
Let $\widetilde{\fm{}}= \widehat{A[0,1]}L \itpr{A^{\sim}[0,1] L}
(\fm{}[0,1]) L$, and let $\widetilde{\fop{}}$ be an 
$(1\tpr \fop{})$-connection
for $\widehat{A[0,1]}L$. Consider
$\widetilde{u}=\{ (1-s)1+su \}_{s\in [0,1]} \in
\cval^{(0)} (\widehat{A[0,1]}L)$.
Finally, let $\widetilde{h}$ be any section of $LL$ that makes 
$(\widetilde{\ec}_{\widetilde{h}}, \widetilde{F}_{\widetilde{h}})=
\res{\widetilde{h}}{(\widetilde{\ec},(\widetilde{u}\tpr 1)\widetilde{F})}$
an asymptotic Kasparov $G$-$(A,B[0,1])$-module. Then
$(\widetilde{\ec}_{\widetilde{h},0},\widetilde{F}_{\widetilde{h},0})$
is homotopic (via a `stretching') with $(\fm{},\fop{})$, and
$(\widetilde{\ec}_{\widetilde{h},1},\widetilde{F}_{\widetilde{h},1})$
represents $\Kprod{1_A}{A}{x}$.
\end{proof}

\begin{rmk}
Theorem \ref{T:assx} and Proposition \ref{T:1} imply that,
for any \gcs{G}-algebra $A$, $KE_G(A,A)$ is a ring with unit.
\end{rmk}

The following notion is important in further studying the properties
of $KE$-theory and in applications.
 
\begin{definition}
Let $D_1$ and $D_2$ be \gcs{G}-algebras. An element
$\ga\in KE_G(D_1,D_2)$ is called {\em $KE$-equivalence}
\index{KE-equivalence@$KE$-equivalence}
(or {\em invertible}) if there exists an element $\gb \in 
KE_G(D_2,D_1)$ such that $\Kprod{\ga}{D_2}{\gb}=1_{D_1}$
and $\Kprod{\gb}{D_1}{\ga}=1_{D_2}$. 
If such an element $\ga$ exists then
$D_1$ and $D_2$ are called {\em $KE$-equivalent}.
(See \cite[2.17]{Kas88}, \cite[19.1]{Blck}.)
\end{definition}

We use $KE$-equivalence to state a result that bears considerable
theoretical significance:

\begin{thm}[Stability in $KE$-theory]\label{T:stability}
\index{stability of KE@stability in $KE$-theory}
For any \gcs{G}-algebra $A$, $A$ and $A\tpr \cpct(\hs_G)$ are
$KE$-equivalent.
\end{thm}

The proof follows from the corresponding result in $KK$-theory,
as explained in Corollary \ref{C:KE_equiv_via_KK}.
Another proof can be given by rephrasing \cite[2.18]{Kas88}
in terms of $KE$-theory groups.

\begin{cor}
For any separable \gcs{G}-algebras $A$ and $B$, we have
$$
KE_G(A,B)\simeq KE_G(A,B\tpr\cpct(\hs_G))
   \simeq KE_G(A\tpr\cpct(\hs_G),B)\simeq
KE_G(A\tpr\cpct(\hs_G),B\tpr\cpct(\hs_G)).
$$ 
\end{cor}

We end this section by defining in $KE$-theory (as it is the
case in $KK$-theory and $E$-theory) the 
{\em higher order groups}.
We recall that $\pcli{n}$ is the Clifford algebra of
$\R^n$, \ie the universal algebra with odd generators 
$\{e_1,...,e_n\}$ satisfying
$e_i e_j + e_j e_i = + 2 \delta_{ij}$, for $1\leq i,j\leq n$,
${e_i}^* = + e_i$, 
and $\| e_i\|=1$.
(The grading is the standard one, and the
notation coincides with the one from \cite{Kas75}.
The adjoint and the norm refer to the fact that $\pcli{n}$ can be
given the structure of a \cs -algebra.)

\begin{definition}\label{D:higherKE}
$KE^n_G(A,B) = KE_G(A,B\tpr\pcli{n})$, for $n=1,2,\dots$.
\index{KE-theory@$KE$-theory!groups!higher order|textbf}
\end{definition}


\subsection{The proof of the technical theorem}\label{S:TT}

In this subsection the following is proved:

\begin{tthm}[Theorem \ref{T:TT}]
\index{product!in $KE$-theory!Technical Theorem|(}
Let $G$ be a locally compact $\sigma$-compact Hausdorff group, and let
$A$, $B$, and $D$ be separable graded \gcs{G}-algebras.
Consider two asymptotic Kasparov modules
$( \fm{1}, \fop{1} ) \in ke_G( A, D )$ and 
$( \fm{2}, \fop{2} ) \in ke_G( D, B )$.
There exists a triple $(\, u, \cnn{F}, h_0 \,)$, with $u$ 
a self-adjoint element of $\cval^{(0)}( \fm{1} )$,
$\cnn{F}$ an $\fop{2}$-connection for $\fm{1}$, and $h_{0}$ 
a section of $[1,\infty)\times [1,\infty)$, as in Overview
\ref{overview}, such that for any 
other section $h\geq h_{0}$ 
$$
( \fm{h}, \fop{h} ) = (\text{\rm Res}_h)_* 
  \big(\;
     \fm{1}\itpr{DL} \fm{2}L,
     \fop{1}\itpr{DL} 1 + (u \itpr{DL} 1)\,\cnn{F}
  \;\big) 
$$
is an asymptotic Kasparov $G$-$(A,B)$-module.
\end{tthm}
\begin{proof}
We shall justify Steps 1-4 of the Overview \ref{overview}. 
{\it Step 1}, in which $u$ is constructed, is nothing but 
Lemma \ref{L:constru} applied to $( \fm{1}, \fop{1} )$.
The existence of the connection $\cnn{F} = \cnn{F}^{\,*}$ in 
{\it Step 2} follows from Proposition \ref{P:excnn} (and after
choosing $\fop{2}=\fop{2}^*$). As it will become clear from the proof,
the self-adjointness of $\cnn{F}$ is just a convenience. It enables
us to reduce some of the computations to the unified
requirements of Step 3.
So far we succeeded to create the pair of `two-dimensional' objects
$( \fm{}, \fop{} ) = ( \fm{1}\itpr{DL} \fm{2}L,
\fop{1}\itpr{DL} 1 + (u \itpr{DL} 1)\,\cnn{F} )$.
For {\it Step 3}, we obtain $h_{00}$ by applying Lemma \ref{L:constrh}
for the set $K=\{ u, u\fop{1}, u a_1, u a_2, \dots, u a_n, \dots \}$,
where $\{\, a_n \,\}_{n=1}^{\infty}$ is a dense subset of $A$.
The essential {\it Step 4} is concerned with finding an appropriate
section $h_0$ 
such that 
$( \fm{h_0}, \fop{h_0} ) = \res{h_0}{( \fm{}, \fop{} )}$
will be the asymptotic Kasparov $G$-$(A, B)$-module which
represents the product. For this to happen, the axioms (aKm1)--(aKm4)
must be satisfied. The tensor products that appear below are all
inner (over $DL$), but the \cs -algebra will be omitted in order to
simplify the writing.
%
%

\mic
\noindent $\bullet$
The simple computation: ($\fop{} - \fop{}^*)(a\tpr 1) =
\big( \fop{1}\tpr 1 + (u\tpr 1)\, \cnn{F} - \fop{1}^*\tpr 1 -\cnn{F}^*
\, (u\tpr 1)\big) (a\tpr 1) = (\fop{1}-\fop{1}^*) a \tpr 1 + 
[u\tpr 1,\cnn{F}] (a\tpr 1)$, shows
that (aKm1) for $\fop{h}$ is satisfied for {\it any}
$h\geq h_{00}$, due to (aKm1) for $\fop{1}$, and (5) of Step 3.
%
%

\mic
\noindent $\bullet$
Next, given $a\in A$, we have $[\fop, a\tpr 1] = \fop{1}a\tpr 1 +
(u\tpr 1)\, \cnn{F}\, (a\tpr 1) - (-1)^{\der a} a\fop{1}\tpr 1 -
(-1)^{\der a} (au\tpr 1)\, \cnn{F} = [\fop{1}, a]\tpr 1 -
(-1)^{\der a} [ua\tpr 1, \cnn{F}] + (-1)^{\der a} ([u,a]\tpr 1)\,
\cnn{F} + [u\tpr 1, \cnn{F}]\, (a\tpr 1)$. Consequently (aKm2)
for $\fop{h}$ is also satisfied for {\it any}
$h\geq h_{00}$, because of (aKm2) for $\fop{1}$, and (2), (5), and (7).
%
%

\mic
\noindent $\bullet$
For (aKm3), it is noted that:
$$
\begin{aligned}
a \; &\big( \fop{}^2 -1 \big)\; a^* \\
&= (a\tpr 1) \;\Big(\; \fop{1}^2\tpr 1 + (u\tpr 1)\cnn{F}(\fop{1}\tpr 1)
    + (\fop{1}\tpr 1) (u\tpr 1) \cnn{F} 
     - (u\tpr 1) \cnn{F} \,[ \cnn{F},u\tpr 1 ] \\ 
&\qquad\qquad\qquad + (u\tpr 1) \cnn{F}^2 (u\tpr 1) -1 
             \;\Big)\; (a^*\tpr 1) \\
&\sim \big( (au)\, \fop{1}^2\, (au)^* \big)\tpr 1 
      + (1-u^2)\,\big( a (\fop{1}^2-1) a^* \big)\tpr 1 
      - (a\tpr 1)(u\tpr 1) \cnn{F} \,[ \cnn{F},u\tpr 1 ](a^*\tpr 1) \\
&\qquad\qquad + (a\tpr 1) \big(\;
                               ([\fop{1},u]\tpr 1)\cnn{F}
                              +[u\fop{1}\tpr 1, \cnn{F} ] 
                              +[ u\tpr 1,\cnn{F} ](\fop{1}\tpr 1)
                           \;\big) (a^*\tpr 1) \\
&\qquad\qquad + (a\tpr 1)\,(u\tpr 1)\,\big( \cnn{F}^2 -1 \big)\,
                   (u\tpr 1)\,(a\tpr 1)^*, 
  \text{ modulo } \loc(\fm{1})\itpr{DL} 1.
\end{aligned}
$$
(For the second equality $\sim$ above, we used (1) and (2) of Step 1,
and the self-adjointness of $u$.)
The restriction of the first six terms to {\it any} $h\geq h_{00}$
will give a positive element modulo $\loc(\fm{h})$, because of (3), (5)
and (6). So we shall have (aKm3) satisfied provided that
\begin{equation}\label{E:key}
(u\tpr 1)\,\big( \cnn{F}^2 -1 \big)\,(u\tpr 1)
\text{ restricts to a positive element modulo }
\cval(\fm{h})+\loc(\fm{h}).
\end{equation}
Showing (\ref{E:key}) is a critical point in the construction.
Let $\{I_n\}_{n=0}^{\infty}$
be a cover of $[1,\infty)$ by closed intervals of the form $I_n =
[t_n,t_{n+2}]$, with $t_0=1$, and $\{t_n\}_n$ being a strictly 
increasing sequence with $\lim_{n\ra\infty} t_n = \infty$.
Let $T_{1,n}=t_n$, for $n\geq 0$, and $T_{2,0}=1$. 
If $\{ \mu_n \}_{n=0}^{\infty}$ is a partition of unity
subordinated to this cover, then
$u\tpr 1 = \sum_{n=0}^{\infty} (\mu_n u\tpr 1)$. For each $n\geq 1$, we
can approximate $(\mu_n u\tpr 1)$ by a {\it self-adjoint} finite rank
operator
\begin{equation}\label{E:TT1}
K_n = \sum_{i=1}^{N_n} T_{\fv{\xi_i}}\,T_{\fv{\eta_i}}^*
    = \sum_{i=1}^{N_n} T_{\fv{\eta_i}}\,T_{\fv{\xi_i}}^*,
\text{ with } \fv{\xi_i}, \fv{\eta_i}\in\fm{1}|_{I_n},
\text{ for } i=1, 2, \dots, N_n, 
\end{equation}
and such that $\| (\mu_n u\tpr 1) - K_n \| < 
1/(24 n (\|\fop{2}\|^2+1))$. Note that:
\begin{equation}\label{E:TT2}
\begin{aligned}
K_n \;\big( \cnn{F}^2 -1 \big)\; K_n^* 
&= \sum_{i,j=1}^{N_n} T_{\fv{\xi_i}}\,T_{\fv{\eta_i}}^*
     \;\big( \cnn{F}^2 -1 \big)\; T_{\fv{\eta_j}}\,T_{\fv{\xi_j}}^* \\
&\sim \sum_{i,j=1}^{N_n} T_{\fv{\xi_i}}\,T_{\fv{\eta_i}}^*
      T_{\fv{\eta_j}}\;\big( (1\tpr\fop{2}^2) -1 \big)\; T_{\fv{\xi_j}}^* ,
 && \text{ modulo } \loc(\fm{}) \\
&= \sum_{i,j=1}^{N_n} T_{\fv{\xi_i}}\, 
      \langle \fv{\eta_i}, \fv{\eta_j}\rangle
                     \;\big( (1\tpr\fop{2}^2) -1 \big)\; T_{\fv{\xi_j}}^*.
\end{aligned}
\end{equation}
There exists $\tau_{n,1}$ such that $\big\| (\cnn{F}^2 T_{\fv{\eta_j}} -
T_{\fv{\eta_j}} (1\tpr\fop{2}^2) )_{(t_1,t_2)} \big\| <
1/(12 n N_n^2)$, for all $\fv{\eta_j}$, all $t_1\in I_n$, and all
$t_2>\tau_{n,1}$. This implies that the error of the commutation that
was used for the second line of equation (\ref{E:TT2}) is smaller
than $1/(12 n)$, in norm and when restricted to the graph of 
any section $h$ whose values on $I_n$ are bigger than $\tau_{n,1}$.
Using the characterization of positive operators on Hilbert modules
\cite[4.1]{Lan} that generalizes the familiar one from Hilbert space theory,
we see that the matrix 
$P = \left( \langle \fv{\eta_i}, \fv{\eta_j}\rangle \right)
\in M_{N_n}(DL)$ is positive. Consequently $P = Q Q^*$,
with $Q = \left( d_{ij} \right)$, and we get:
\begin{equation}\label{E:TT3}
\begin{aligned}
\sum_{i,j=1}^{N_n} T_{\fv{\xi_i}}\,
      &\langle \fv{\eta_i}, \fv{\eta_j}\rangle
                     \;\big( (1\tpr\fop{2}^2) -1 \big)\; T_{\fv{\xi_j}}^* 
= \sum_{i,j=1}^{N_n} T_{\fv{\xi_i}}\,
     \big( \sum_{k=1}^{N_n} d_{ik}\, d_{jk}^* \big)
        \;\big( (1\tpr\fop{2}^2) -1 \big)\; T_{\fv{\xi_j}}^* \\
&\sim \sum_{k=1}^{N_n} \Big(
    \big(\sum_{i=1}^{N_n} T_{\fv{\xi_i}}\, d_{ik} \big)
     \;\big( (1\tpr\fop{2}^2) -1 \big)\;
      \big(\sum_{j=1}^{N_n} T_{\fv{\xi_j}}\, d_{jk} \big)^* \Big),
            \qquad\text{ modulo } \loc(\fm{}).
\end{aligned}
\end{equation}
There exists $\tau_{n,2}$ such that 
$\big\| [d_{jk},1\tpr\fop{2}^2]_{(t_1,t_2)} \big\| < 1/(12 n N_n^3)$,
for all $d_{jk}$,  all $t_1\in I_n$, and all $t_2>\tau_{n,2}$.
This implies that the error due to asymptotic commutativity 
((aKm2) for $\fop{2}$, used to obtain the second line of equation 
(\ref{E:TT3})) is smaller than $1/(12 n)$,
in norm and when restricted to the graph of
any section $h$ whose values on $I_n$ are bigger than $\tau_{n,2}$.
Let $\{ \delta_m \}_m$ be an approximate unit in $D$.
Because of (aKm3) for $\fop{2}$,
\begin{equation}\label{E:TT4}
\sum_{k=1}^{N_n} \Big(
  \big(\sum_{i=1}^{N_n} T_{\fv{\xi_i}}\, d_{ik} \big)\,
     \big( \: 1\tpr\delta_m ( \fop{2}^2 -1 )\delta_m \: \big)\,
      \big(\sum_{j=1}^{N_n} T_{\fv{\xi_j}}\, d_{jk} \big)^* \Big)
\end{equation}
is positive modulo $\cval(\fm{}|_{I_n})+\loc(\fm{}|_{I_n})$.
Choose $m_0$ such that the entire sum from (\ref{E:TT4}) 
approximates the one from the second line of (\ref{E:TT3})
by $1/(12 n)$. 

Let $T_{2,n}=\max \{ \tau_{n,1}, \tau_{n,2}, T_{2,(n-1)}+1 \}$.
(To be precise, there is also an $\tau_{n,3}$ coming from (aKm4) 
to be taken into account, but we ignore it for the moment.)
Once the sequence $\{T_{2,n}\}_n$ has been constructed, we define
$h_0$ on $[T_{1,n}, T_{1,(n+1)}]$
as the linear function satisfying $h_0(T_{1,n})=T_{2,n}$
and $h_0(T_{1,(n+1)})=T_{2,(n+1)}$.
The estimates above show that the restriction to the graph of
$h_0|_{I_n}$ of 
$(\mu_n u\tpr 1) \, ( \cnn{F}^2 -1 )\, (\mu_n u\tpr 1)^*$ is positive
modulo $\cval(\fm{h_0})$, with an error which is smaller than
$1/(3 n)$, in norm. At most three such terms are non-zero over $I_n$,
this
proves (\ref{E:key}) for any $h\geq h_0$, and consequently $\fop{h}$
satisfies (aKm3).
%
%

\mare
\noindent $\bullet$
Finally, for any $g\in G$, we have:
$$
\begin{aligned}
\big(g ( F ) - F \big) (a\tpr 1)
    &= \big( g\, (\fop{1}\tpr 1) + g\, (u \tpr 1)\; g\,(\cnn{F}) 
          - (\fop{1}\tpr 1) - (u \tpr 1)\,\cnn{F}\big) (a\tpr 1) \\
    &= \big(\, g(\fop{1}) - \fop{1}\big) a \tpr 1 
        + \big(( g(u)-u) \tpr 1\big)\, g(\cnn{F}) (a\tpr 1) \\
    &\qquad\qquad+ (u \tpr 1) 
       \big( g(\cnn{F}) - \cnn{F} \big) (a\tpr 1).
\end{aligned}
$$
Due to (aKm4) for $\fop{1}$ and (4) of Step 1, 
the first two terms put no extra constraints on $h_0$. 
For the third one, $u\tpr 1$ can be approximated,
as in the proof of (aKm3), on each
interval $I_n$, by a finite sum 
$\sum_{i} T_{\xi_i} T_{\eta_i}^*$. A simple computation shows that
$g\, T_{\eta_i}^* = T_{g(\eta_i)}^*$. Consequently:
$$
\begin{aligned}
T_{\xi_i} T_{\eta_i}^*\big(g ( \cnn{F} ) - \cnn{F} \big)  
&= T_{\xi_i} \, g \big(g^{-1}(T_{\eta_i}^*)\,\cnn{F}\big)
       - T_{\xi_i} T_{\eta_i}^* \cnn{F} \\
&\sim (-1)^{\der \eta_i} \,
    T_{\xi_i} \, g\big( \fop{2} T_{g^{-1}(\eta_i)}^*\big)
       - (-1)^{\der \eta_i}\, T_{\xi_i} \fop{2} T_{\eta_i}^*,
 && \text{modulo } \loc(\fm{}) \\
&= (-1)^{\der \eta_i}\,
     T_{\xi_i} \,\big( g(\fop{2})-\fop{2}\big) \,T_{\eta_i}^*.
\end{aligned}
$$
Further modification (increase) of $h_0$, using (aKm4) for $\fop{2}$,
will make the above errors go to zero when restricted to the graph of
$h_0$. (This is the place where the $\tau_{n,3}$ mentioned 
when we defined $T_{2,n}$
makes its appearance.) This shows that (aKm4) holds for $\fop{h}$,
for any $h\geq (\text{new }h_0)$, and the proof of the Technical 
Theorem is complete.
\index{product!in $KE$-theory!Technical Theorem|)}
\end{proof}

\begin{rmk}
The only important fact that $h_0$ encodes in the construction of 
the product is a certain behavior that occurs when $t_1\ra \infty$ 
and $t_2\ra\infty$, with $h_0$ correlating $t_1$ and $t_2$.
We have noticed that certain decay properties hold true on entire
`stripes' $[T_{1,n},T_{1,n+1}]\times[T_{2,n},\infty)$, and not
only on the graph of $h_0$. This observation is used in the proof
of the associativity of the product (see Lemma \ref{L:ass1}),
where it allows us to focus on the analysis of the operators that appear
in the construction rather than on the sections.
\end{rmk}


\section{$KE$-theory: some examples}

We further investigate, by means of examples, the significance 
the axioms (aKm1)--(aKm4) that lie at the foundation of 
$KE$-theory. A consequence of our discussion is
that non-equivariant $KE$-theory groups recover
the ordinary $K$-theory for trivially graded \cs -algebras.

\subsection{A non-equivariant example: $K$-theory}
\label{ss:ktheory}
In this subsection the \cs -algebras are trivially graded (ungraded),
and there is no group action. We consider that there is some
merit in the proof of the next result:

\begin{prop}\label{P:ktheory}
\index{KE-theory@$KE$-theory!in relation to!$K$-theory}
Let $B$ be an ungraded separable \cs -algebra, then 
$$
KE^*(\C,B) \simeq  KK^*(\C,B) \simeq K_*(B), \; \text{ for } *=0,1.     
$$
\end{prop}

\begin{proof}
The second group isomorphism, $KK^*(\C,B) = K_*(B)$, is well-known
(see for example \cite[17.5.6, 17.5.7]{Blck}). 
Consequently the main point behind this proposition is the 
following: we shall show that the axioms of Kasparov modules
(Definition \ref{D:Kasmod}, (\ref{E:axkk})) can be successively 
modified, in the case when $A=\C$, to give the 
axioms (aKm1--3) of asymptotic Kasparov modules 
(Definition \ref{D:akm}). In this way we obtain intermediate
abelian groups $\widetilde{KK}(\C,B)$, $\widetilde{KE}(\C,B)$,
and group homomorphisms $\alpha$, $\beta$, $\gamma$ between the 
four groups under consideration, that can be depicted
in the diagram:
\begin{equation}\label{E:kth1}
\begin{CD}
KK(\C,B) @>{\ga}>> \widetilde{KK}(\C,B) 
         @>{\gb}>> \widetilde{KE}(\C,B)
         @>{\gamma}>> KE(\C,B) \; .
\end{CD}
\end{equation}
(Note that $\widetilde{KK}(\C,B)$ has nothing to do with
the group denoted by same symbol in \cite[Def.2(8)]{Sk84},
which is the quotient of $kk(\C,B)$ by the equivalence relation 
generated by addition of degenerate elements and operatorial 
homotopy.)  The claimed isomorphism between the $KK$-theory group
and the $KE$-theory group is deduced from the fact that 
$\alpha$, $\beta$, and $\gamma$ are proven to be isomorphisms.

$\widetilde{KK}(\C,B)$ is the abelian group (under direct sum) of 
homotopy classes of pairs $(\ec,F)$, where $\ec$ is a 
Hilbert $B$-module, admitting an action of $\C$
via a \sh homomorphism $\gvf:\C \ra \bdd(\E)$,
and $F\in\bdd(\E)$ is an odd operator such that:
\begin{equation}\label{E:kth2}
\gvf(1)=\text{id}, \, F=F^*, \, \text{ and }  
        (F^2-1/2) \geq 0, \,\text{ modulo } \cpct(\ec).
\end{equation}
(See Remark \ref{R:kk-ke}.) To construct the group homomorphism
$\ga : KK(\C,B) \ra \widetilde{KK}(\C,B)$ we recall some of the
standard simplifications of the axioms that a Kasparov module has
to satisfy \cite[17.4]{Blck}. Let $(\ec,F)\in kk(\C,B)$ be an 
arbitrary Kasparov module. By replacing $F$ with $F' =
(F+F^*)/2$ we find a homotopic module $(\ec,F')$ with the
operator self-adjoint. Next, consider the projection
$\gvf(1)=P\in \bdd(\ec)$. The pair $(\ec,F')$ is operator homotopic to
$(\ec,PF'P) = (P\ec,PF'P) + ((1-P)\ec,0)$, with the second
summand being degenerate. Consequently, in the homotopy class of the 
initial Kasparov module we find a representative
$(\widetilde{\ec},\widetilde{F}) = (P\ec,PF'P)$, with
$1\in \C$ acting as identity, $\widetilde{F}$ self-adjoint,
and $\widetilde{F}^2=1\geq 1/2$, modulo $\cpct(\widetilde{\ec})$.
This defines the group homomorphism $\ga$ (all the changes above
preserve homotopies and direct sums):
\begin{equation}\label{E:kth3}
\ga : KK(\C,B) \ra \widetilde{KK}(\C,B), \;
   (\ec,F) \mapsto (\widetilde{\ec},\widetilde{F}).
\end{equation}
For the inverse map, let $\psi:\R\ra\R$, be $\psi(x)=-1$, for
$x\leq -1/{\sqrt{2}}$, $\psi(x)=\sqrt{2}\,x$, for $x\in (-1/{\sqrt{2}},
1/{\sqrt{2}})$, and $\psi(x)=1$, for $x\geq 1/{\sqrt{2}}$. Define
$$
\ga': (\widetilde{\ec},\widetilde{F}) \mapsto (\widetilde{\ec},
    \psi(\widetilde{F})).
$$
The only non-trivial checking is $\psi(\widetilde{F})^2-1 = 
2 \widetilde{F}^2-1 \geq 0$, modulo $\cpct(\widetilde{\ec})$.
We observe that $[\psi(\widetilde{F}), \widetilde{F}]\geq 0$, and
consequently both compositions $\ga'\circ\ga$ and $\ga\circ\ga'$
give results homotopic with the initial module. It follows that
$\ga$ is an isomorphism, with $\ga^{-1}=\ga'$.

Define $\widetilde{KE}(\C,B)$ to be the abelian group (under direct
sum) of homotopy classes of
asymptotic Kasparov $(\C,B)$-modules 
$(\ofm{},\ofop{})$ satisfying the {\em extra conditions}:
\begin{equation}\label{E:kth4}
\gvf(1)=\text{id}, \, \ofop{}=\ofop{}^*, \, \text{ and }
        (\ofop{}^2-1/2) \geq 0, \,\text{ modulo } \cval(\ofm{}).
\end{equation}
The map $\gamma: \widetilde{KE}(\C,B)\ra KE(\C,B)$ is the forgetting
map at the level of asymptotic Kasparov modules. To define the inverse
$\gamma '$, let $(\ofm{},\ofop{})$ be an arbitrary asymptotic
Kasparov module. We can make the action of $\C$ unital as in 
$KK$-theory: there is a homotopy followed by a `small perturbation'
connecting $(\ofm{},\ofop{})$ with $(\ofm{\:}',\ofop{\:}'') = 
(P\ofm{},P\ofop{}P)$, where $P=\fr{}(1)$. 
As we have already observed in Corollary \ref{C:self-adj}, 
there is a homotopy from this last pair to another one 
$(\ofm{\:}',\ofop{\:}')$, with $\ofop{\:}'$ self-adjoint.
Finally, (aKm3) implies that $(\ofop{t}')^2 - 1 \geq U_t + V_t$,
with $U=\{U_t\}_t \in \cval(\ofm{\:}')$ and $V=\{V_t\}_t\in\loc(\ofm{\:}')$.
Let $T$ be such that $\|V_t\|<1/2$, for all $t>T$. It follows that
$(\ofop{t}')^2 - 1/2\geq U_t$, for $t>T$. We define $\gamma '$ via
a `translation' (see Example \ref{Ex:htpy}):
$$
\gamma ' : \{\, (\ofm{t},\ofop{t}) \,\}_t \mapsto
   \{\, (\widehat{\ec}_{t+T},\widehat{F}_{t+T}) \,\}_t.
$$
All the operations used to define $\gamma '$ preserve homotopies and 
direct sums, and consequently both $\gamma '\circ \gamma $ and
$\gamma \circ \gamma '$ are identity, and $\gamma^{-1}=\gamma '$.

Finally, define 
\begin{equation}\label{E:kth5}
\gb: \widetilde{KK}(\C,B)\ra  \widetilde{KE}(\C,B),
\; (\widetilde{\ec},\widetilde{F}) \mapsto 
    \{\, (\widetilde{\ec},\widetilde{F}) \,\}_t \;
    \text{ (constant family)}, 
\end{equation}
and
$$
\gb': \widetilde{KE}(\C,B) \ra \widetilde{KK}(\C,B),
\; (\ofm{},\ofop{}) = \{\, (\ofm{t},\ofop{t}) \,\}_t 
    \mapsto (\ofm{1},\ofop{1}) \;
    \text{ (the `fiber' at $t=1$)}.
$$
The composition $\gb'\circ \gb = \text{id}$ is obvious. 
Let now $(\ofm{},\ofop{}) = \{\, (\ofm{t},\ofop{t}) \,\}_t$
be an element of $\widetilde{KE}(\C,B)$. There exists a homotopy
$(\boldsymbol{\fm{}},\boldsymbol{\fop{}})$ between 
$(\ofm{},\ofop{})$ and $(\gb\circ\gb')( (\ofm{},\ofop{}) ) =
\{ \, (\ofm{1},\ofop{1}) \,\}_t$ given by explicit formulas:
$$
\boldsymbol{\ec}_{t,s}=\widehat{\ec}_{s+(1-s)t}, \;
\boldsymbol{F}_{t,s}=\widehat{F}_{s+(1-s)t}, \;
\text{ for } s\in [0,1], t\in [1,\infty).
$$
This proves that $\gb$ is also an isomorphism, with $\gb^{-1}=\gb '$.

The claimed isomorphism is $\gamma\circ\gb\circ\ga: KK(\C,B) \ra
KE(\C,B)$. Finally we get:
$KK^1(\C,B) \stackrel{\scriptstyle\text{def}}{\simeq} KK(\C,B\tpr \pcli{1})
 \stackrel{\scriptstyle\text{as above}}{\simeq} KE(\C,B\tpr \pcli{1}) 
  \stackrel{\scriptstyle\text{def}}{\simeq} KE^1(\C,B)$.
\end{proof}

\subsection{An equivariant example: $KE_{\Gamma}(\C,\C)$,
for $\Gamma$ discrete}
\label{ss:keGamma}

The next result is similar with Remark 2, after \cite[2.15]{Kas88},
namely the dual of Green-Julg theorem in $KK$-theory:

\begin{prop}
Let $\gG$ be a discrete group and $A$ a separable \gcs{\gG}-algebra, then 
$KE_{\gG}(A,\C) = KE(\mcs{\gG,A},\C).$
\end{prop}

\begin{proof}
The presentation of crossed products contained in \cite[Ch.8]{Dvds} 
should be enough to follow the argument below. We start by choosing
$(\fm{},\fop{})\in ke_{\gG}(A,\C)$. Using the Stability Theorem \ref{T:stability}
and the Stabilization Theorem \cite[Thm.2]{Kas80}, we can assume that
$\fm{}=HL$, for a fixed Hilbert space $H$ (see \ref{P:KEstabilization}).
The field of Hilbert spaces $\fm{}$ is endowed with a unitary action 
$U:\gG \ra \bdd(\fm{})$, and an equivariant \sh representation 
$\fr{}:A \ra \bdd(\fm{})$. In terms of families $(\fm{},\fop{})$ gives
bounded and \sh strong continuous families 
$\{ \, U_t:\gG \ra \uni(H) \,\}_t$ and $\{\, \fr{t}:A \ra \bdd(H) \,\}_t$. 
We denote $U_t(g)$ by $g_t\in \uni(H)$, for each $\ui$. 
The equivariance of each $\fr{t}$ implies that we actually 
have a covariant representation of the dynamical system $(A,\gG)$. 
Consequently we can construct actions 
$\widetilde{\fr{t}}: C_c(\gG,A) \ra \bdd(H)$ 
in the usual way:
$$
\widetilde{\fr{t}}(f)=\sum_{g\in\gG} \fr{t}(a_g) \, g_t,
 \text{ for } f=\sum_{g\in\gG} a_g \, \gd_g \in C_c(\gG,A).
$$
Note that $\{\, \widetilde{\fr{t}}(f) \,\}_t$ is bounded and \sh strong
continuous for each $f\in C_c(\gG,A)$. By the norm density of $C_c(\gG,A)$
in $\mcs{\gG,A}$ we obtain representations 
$\widetilde{\fr{t}}: \mcs{\gG,A} \ra \bdd(H)$. We claim that with this
representation $\widetilde{\fr{}}: \mcs{\gG,A} \ra \bdd(\fm{})$
the asymptotic Kasparov module $(\fm{},\fop{})$ gives an element
$(\widetilde{\fm{}},\widetilde{\fop{}})= (\fm{},\widetilde{\fr{}},\fop{})$
in $ke(\mcs{\gG,A},\C)$. It is enough to check the axioms for 
$f=a_g \, \gd_g \in C_c(\gG,A)$.
\newline $\bullet$ {\it $\widetilde{\fop{}}$ satisfies} (aKm1).
$(\widetilde{\fop{}} - \widetilde{\fop{}}^*) \widetilde{\fr{}}(f)
 = (\fop{}-\fop{}^*) \fr{}(a_g) g \sim 0$, by (aKm1) for $\fop{}$.
\newline $\bullet$ {\it $\widetilde{\fop{}}$ satisfies} (aKm2). Indeed:
$ 
[ \widetilde{\fop{}} , \widetilde{\fr{}}(f)]
 = \fop{}  \fr{}(a_g) g - (-1)^{\der a_g} \fr{}(a_g) g \fop{}
 = [\fop{}, \fr{}(a_g)] g + \fr{}(a_g) (\fop{} - g (\fop{})) g
 \sim 0, \text{ by (aKm2) and (aKm4) for } \fop{}.
$
\newline $\bullet$ {\it $\widetilde{\fop{}}$ satisfies} (aKm3).
$$
\begin{aligned}  
\widetilde{\fr{}}(f) \,(\widetilde{\fop{}}^2 -1)\, \widetilde{\fr{}}(f)^*
 &= \fr{}(a_g) \, g \, (\fop{}^2-1) \, g^{-1} \, \fr{}(a_g) \\
 &\sim \fr{}(a_g) \, (\fop{}^2-1) \,  \fr{}(a_g)
  && \text{ by (aKm4) for } \fop{}    \\
 &\geq 0, \text{ modulo } \cval(\fm{})+\loc(\fm{}),
  &&  \text{ by (aKm3) for } \fop{}.  
\end{aligned}   
$$
The computation above shows also that a homotopy in $ke_{\gG}(A,C([0,1]))$
is sent to a homotopy in $ke(\mcs{\gG,A},C([0,1]))$. We obtain in this
way a group homorphism $KE_{\gG}(A,\C)\ra KE(\mcs{\gG,A},\C)$,
$\keclass{(\fm{},\fop{})}\mapsto \keclass{(\widetilde{\fm{}},\widetilde{\fop{}})}$.
Now for the inverse group homomorphism, consider and asymptotic Kasparov module 
$(\widetilde{\fm{}},\widetilde{\fop{}})\in ke(\mcs{\gG,A},\C)$,
where $\widetilde{\fm{}}=HL$, for a fixed Hilbert space $H$, and
$\widetilde{\fr{}}:\mcs{\gG,A}\ra \bdd(\widetilde{\fm{}})$.
Let $\{ \alpha_n\}_n$ be an approximate unit for $A$. We obtain
by a standard construction bounded and \sh strong continuous families of
representations ($\ui$):
$$
\fr{t}:A\ra \bdd(H), \; \fr{t}(a)=\widetilde{\fr{t}}(a \gd_e),
\text{ with $e$ the unit of }\gG,
$$
and
$$
U_t:\gG \ra \uni(H), \; U_t(g)=\lim_{n\ra \infty}\widetilde{\fr{t}}(\alpha_n \gd_g),
\text{ for } g\in\gG.
$$
Denote by $(\fm{},\fop{})$ the pair $(\widetilde{\fm{}},\widetilde{\fop{}})$
with the actions of $\gG$ and $A$ obtained in this way. The claim is that
$(\fm{},\fop{})$ belongs to $ke_{\gG}(A,\C)$, and the only non-trivial
axiom to be checked is (aKm4). We have:
$$
\begin{aligned}  
\left( \, g(\fop{})-\fop{} \,\right)\, \fr{}(a)
 &= \lim_{n\ra \infty} \left( \widetilde{\fr{}}(\alpha_n \gd_g) \,\widetilde{\fop{}}
       \,\widetilde{\fr{}}(\alpha_n \gd_{g^{-1}}) - \widetilde{\fop{}} \right)
        \widetilde{\fr{}}(a \gd_e) \\
 &= \lim_{n\ra \infty} \left( \widetilde{\fr{}}(\alpha_n \gd_g) \,
       [\widetilde{\fop{}},\widetilde{\fr{}}(\alpha_n a \gd_{g^{-1}})] \right.\\
 &\qquad\qquad + \left. (-1)^{\der a} 
       \left( \widetilde{\fr{}}(\alpha_n^2 a \gd_e)- 
                 \widetilde{\fr{}}(a \gd_e) \right) \widetilde{\fop{}}
        + (-1)^{\der a} [\widetilde{\fr{}}(a \gd_e),\widetilde{\fop{}}\,] \right)\\
 &\sim 0, \text{ by (aKm2) for } \widetilde{\fop{}}.
\end{aligned}   
$$
Again homotopies in $ke(\mcs{\gG,A},C([0,1]))$ are sent to homotopies in
$ke_{\gG}(A,C([0,1]))$, and it is clear that we obtain in this way the
inverse group homomorphism.
\end{proof}

\noindent
{\em Note.} Using Bott periodicity we also have:
$KE^1_{\gG}(A,\C) \stackrel{\scriptstyle\text{Bott}}{\simeq} 
  KE_{\gG}(A\tpr \pcli{1},\C) \stackrel{\scriptstyle\text{above}}{\simeq} 
   KE(\mcs{\gG,A\tpr \pcli{1}},\C) \simeq
    KE(\mcs{\gG,A}\tpr \pcli{1},\C) \stackrel{\scriptstyle\text{Bott}}{\simeq} 
     KE^1(\mcs{\gG,A},\C)$.


\section{Final remarks: the relation with $KK$-theory and $E$-theory}
\label{S:kke}

Assume that a group G (locally compact, $\gs$-compact, Hausdorff)
is given.
In this final section we construct two functors: 
$\kke: \mathbf{KK_G} \ra \mathbf{KE_G}$, and
$\kee: \mathbf{KE_G} \ra \mathbf{E_G}$.
The three categories have all the same objects:
the separable and graded \gcs{G}-algebras.
The morphisms of $\mathbf{KK_G}$ 
\index{$KK$xcategory@$\mathbf{KK_G}$}
\index{KK-theory@$KK$-theory!a category theory perspective,$\mathbf{KK_G}$}
(\cite{Hg87a}, \cite{Hg90b}, \cite{Thms98}) 
are the $KK$-theory groups, with
composition given by the Kasparov product (see Theorem \ref{T:xkk}).
The morphisms of $\mathbf{KE_G}$ are the $KE$-theory groups,
\index{$KE$xcategory@$\mathbf{KE_G}$}
\index{KE-theory@$KE$-theory!a category theory perspective,$\mathbf{KE_G}$}
with composition given by the product defined in Section \ref{S:prod}.
The morphisms of $\mathbf{E_G}$ (\cite{GHT}, \cite{HgKas97})
\index{$E$xcategory@$\mathbf{E_G}$}
\index{E-theory@$E$-theory!a category theory perspective,$\mathbf{E_G}$}
are the $E$-theory groups, 
with the corresponding composition product.
Both functors are the identity on objects.

One consequence of the existence of these two functors 
is the construction of an {\em explicit}
natural transformation, namely the composition $\kee\,\circ\,\kke$,
between $KK$-theory and $E$-theory. This transformation 
{\em preserves the product structures} of the two theories. 
This connecting functor is roughly:
\begin{equation}\label{Eq:compositionKKE}
\begin{array}{ccccc}
KK_G(A,B) & \stack{\,\;\kke\,\;} & KE_G(A,B) &
 \stack{\;\,\kee\,\;} & E_G(A,B) \\ [2ex]
(\ec,F) & \mapsto &
  \big\{ \big( \ec,(1-u_t) F (1-u_t) \big) \big\}_{t} & \mapsto &
\big\{ f\tpr a \stackrel{\gvf_t}{\mapsto} 
         f\big( (1-u_t) F (1-u_t) \big) \, a \big\}_{t}
\end{array}
\end{equation}
\nopagebreak
Here $(\ec,F)$ is a Kasparov module, $\{ u_t \}_t$
is a quasi-invariant quasi-central approximate unit for $\cpct(\ec)$,
and $\{\gvf_t\}_t:\CI \tpr A\asy B\tpr\cpct$ is an asymptotic family.
The suggestive but somehow imprecise (see subsection \ref{sS:kee})
formula of the composition in (\ref{Eq:compositionKKE}), 
namely $\kee\,\circ\,\kke : KK_G(A,B)\ra E_G(A,B)$,
$(\ec,F)\mapsto 
   \{ f\tpr a \stackrel{\gvf_t}{\mapsto} 
         f\left( (1-u_t) F (1-u_t) \right) \, a \}_{t}$
appears also in \cite[4.5.1]{Pop00}, in the context of groupoid
actions.

\subsection{The map $KK_G \ra KE_G$}

Let $G$ be a group, $A$ and $B$ be \gcs{G}-algebras. 
Consider $(\fm{}, \fop{})\in kk_G(A, B)$.
This means that $\ec$ is a graded Hilbert $G$-$B$-module,
acted on by $A$, and $F\in\bdd(\ec)$ is an odd operator such that
$(F-F^*) a$, $[F,a]$, $(F^2-1) a$, $(g(F)-F) a\,$ belong to
$\cpct(\ec)$, for all $a\in A$, $g\in G$.
Denote by $C^*(\cpct(\ec),A,F)$ the smallest \cs -subalgebra of 
$\bdd(\ec)$ that contains $\cpct(\ec)$, $\gvf(A)$, and $F$,
and let $u=\{ u_t \}_{\ui}$ be a quasi-invariant quasi-central 
approximate unit for 
$\cpct(\ec)\subset C^*(\cpct(\ec),A,F)\subset \bdd(\fm{})$.
It will be convenient, at least for notational purposes, to 
regard $u$ as an element of $\cval(\ec L)$.
We make the notation:
$\ofm{} = \ec L$ (constant family of modules), and
$\ofop{} = \{ (1-u_t) F (1-u_t) \}_t = (1-u) F (1-u)$.

\begin{goodclaim}\label{Cl:defkke}
$\big\{ \big( \ec , (1-u_t) F (1-u_t) \big) \big\}_t 
= \big( \ofm{},\ofop{} \big)$ 
is an asymptotic Kasparov $G$-$(A,B)$-module.
\end{goodclaim}

With this result at our disposal, the connection between the 
$KK$-theory and $KE$-theory groups is given by the following
two results:

\begin{thm}\label{T:kkemap}
Let $G$ be a group, $A$ and $B$
be \gcs{G}-algebras. Consider $(\fm{}, \fop{})\in kk_G(A, B)$,
and let $u=\{ u_t \}_{\ui}$ be a quasi-invariant quasi-central 
approximate unit for 
$\cpct(\fm{})\subset 
\bdd(\fm{})$. The map
\begin{equation}\label{E:kke}
\kke : kk_G(A,B) \ra ke_G(A,B) , \, 
   \big(\,\fm{}, \fop{}\,\big)\mapsto
      \big\{\,\big(\fm{},(1-u_t) F (1-u_t)\big)\,\big\}_{\ui} ,
\end{equation}
passes to quotients and gives a group homomorphism
$\kke : KK_G(A,B) \ra KE_G(A,B)$.%
\index{connecting map!between $KK$-theory and $KE$-theory, $\kke$|textbf}
\index{$xXHeta$@$\kke$, map $KK_G\ra KE_G$|textbf}
\end{thm}

\begin{thm}\label{T:kkeX}
$\kke: \mathbf{KK_G} \longrightarrow \mathbf{KE_G}$ is a functor.
\end{thm}

One consequence is worth noticing:

\begin{cor}\label{C:KE_equiv_via_KK}
A $KK$-equivalence is sent by $\kke$ into a $KE$-equivalence.
In particular we obtain that $A$ and $A\tpr \cpct(\hs_G)$ are 
$KE$-theory equivalent, for any \gcs{G}-algebra $A$, and that the 
$KE$-theory groups satisfy Bott periodicity.
\end{cor}

\subsection{The map $KE_G \ra E_G$}\label{sS:kee}

The $E$-theory groups were introduced and studied in 
\cite{CoHg89}, \cite{CoHg90}, 
the equivariant ones under the action of a group in \cite{GHT}, 
and under the action of a groupoid in \cite{Pop}. 
We use here the approach taken in \cite[Sec.2]{HgKas97}.
Let $\s$ be the \cs-algebra $C_0(\R)$ graded by even and odd functions.

\begin{definition}
(\cite[Def.2.2]{HgKas97}) We denote by $E_G(A,B)$%
\index{$E$theory groups@$E_G(A,B)$}
\index{E-theory@$E$-theory!groups}
the set of all homotopy equivalence classes of
asymptotic families from $\s A\tpr\cpct(\hs_G)=\s\tpr A\tpr\cpct(\hs_G)$
to $B\tpr\cpct(\hs_G)$:
$ E_G(A,B) = \lam \, \s A\tpr\cpct(\hs_G) , B\tpr\cpct(\hs_G) \, \ram.$
\end{definition}

Our construction of the connecting map 
between $KE$-theory and $E$-theory is performed via a description of the 
$E$-theory groups which involves $\CI$ instead of $\s$. 
Such a modification seems more appropriate when working with bounded
operators. As with $\s=C_0(\R)$, the \cs -algebra $C_0((-1,1))$ will be
graded by even and odd functions.

Let $G$ be a group, $A$ and $B$
be \gcs{G}-algebras. We consider first a particular case of
asymptotic Kasparov $(A,B)$-modules: $(\fm{}, \fop{}) =
\left\{ ( \hm , \fop{t} ) \right\}_{t} \in ke_G(A,B)$, where $\hm$ is
a {\em fixed} Hilbert $G$-$B$-module acted upon by $A$ through the \sh
homomorphism $\gvf:A\ra\bdd(\hm)$ (or through a family of 
\sh homomorphisms $\gvf_t:A\ra\bdd(\hm)$, but the argument remains
unchanged). This means that $F_t = F_t^*\in \bdd(\hm)$
is an odd self-adjoint operator, for every~$t$, 
such that $[F_t,a]$, $(g(F_t)-F_t)a$
converge in norm to 0 as $t\ra \infty$, for all $a\in A$, $g\in G$,
and that $a (F_t^2 - 1) a^* \geq 0$, modulo compacts, with an error
that converges in norm to 0 as $t\ra \infty$.

\begin{goodclaim}\label{Cl:defkee}
The family of maps 
\begin{equation}\label{E:defkee}
\gf_F=\{ \gf_{F,t} \}_{\ui} \: : \: \CI \tpr A \ra \cpct(\hm),\;
    f \tpr a \stackrel{\textstyle{\gf_{F,t}}}{\mapsto} f(F_t)\, a,
\end{equation}
for $f\in C_0((-1,1))$, $a\in A$, is an asymptotic family, 
in the sense of $E$-theory \cite[Def.1.3]{GHT}.
\end{goodclaim}

The asymptotic family constructed above indicates that a `$\CI$-picture'
of $E$-theory is in order. The next lemma is the first step towards
such a characterization.

\begin{lemma}\label{L:CIpicture1}
Let $A$ and $D$ be \gcs{G}-algebras, and consider an equivariant
asymptotic family $\gf_F=\{ \gf_{F,t} \}_t : \CI\tpr A \asy D$.
Then there exists a unique, up to homotopy, equivariant 
asymptotic family
$\psi_F=\{ \psi_{F,t} \}_t : \s A \asy D$ such that the
diagram
$$
\begin{CD}
\CI\tpr A @>{\gf_F}>> D \\
@V{\text{\rm inclusion}}V{\gi}V @| \\
\s A @>>{\psi_F}> D
\end{CD}
$$
commutes up to homotopy.
\end{lemma}

To discuss the general case we mention the following possible simplification
in the definition of asymptotic Kasparov modules:

\begin{prop}\label{P:KEstabilization}
Given two \gcs{G}-algebras $A'$ and $B'$, let $A=A'\tpr \cpct(\hs_G)$
and $B=B'\tpr \cpct(\hs_G)$. Then, in 
the definition of $KE_G(A,B)$ it is enough to consider modules
of the form $(\hs_{BL},\fop{})$.
\end{prop}

The proposition implies that the previous construction of the
asymptotic morphism associated to an asymptotic Kasparov module
with constant `fibers' can be carried over the general case.
Consider an arbitrary Kasparov module
$(\fm{},\fop{})\in ke_G(A,B)$. We can construct
an asymptotic morphism $\gf : \CI\tpr A \ra \cval(\fm{})/\cpct(\fm{})$.
This in turn gives an asymptotic morphism:
\begin{equation}\label{E:1defkee}
\gf\tpr 1: \CI\tpr A\tpr \cpct(L^2(G))\ra
    \cval(\fm{}\tpr L^2(G))/\cpct(\fm{}\tpr L^2(G)).
\end{equation}
By ignoring the action of $G$, apply the Stabilization Theorem 
(\cite[Thm.2]{Kas80}, with $G$=\{e\}) to get a non-equivariant isometry
$V:\fm{} \ra \hs_{BL}$. 
Apply next the Fell's trick 
to construct an equivariant $BL$-linear isometry 
$W:\fm{}\tpr L^2(G)\ra \hs_{BL} \tpr L^2(G)$.
Use it, and the fact that now we have a {\em constant field} $\hs_{BL}$
of modules, to transform the asymptotic morphism $\gf\tpr 1$ of 
(\ref{E:1defkee}) into an asymptotic family:
\begin{equation}\label{E:2defkee} 
\gf_F: \CI\tpr A\tpr \cpct(L^2(G)) \asy
   \cpct(\hs_{B})\tpr\cpct(L^2(G)).
\end{equation}
After tensoring with $\cpct$, we can use Lemma \ref{L:CIpicture1}
to  obtain an asymptotic morphism
$
\psi_F:\s A\tpr \cpct \asy B\tpr \cpct.
$
The connection between $KE$-theory and $E$-theory is given by:

\begin{thm}\label{T:keemap}
For any group $G$, and any two \gcs{G}-algebras $A$ and $B$,
the map $\kee: (\fm{},\fop{})\mapsto \psi_F$, from asymptotic Kasparov
$G$-$(A,B)$-modules to asymptotic families from $\s A\tpr \cpct$ to
$B\tpr \cpct$, passes to quotients and gives a
natural group homomorphism
\begin{equation}\label{E:3defkee}
\kee : KE_G(A,B) \ra E_G(A,B), \;
\kee ((\fm{},\fop{})) = \eclass{\psi_F}.
\end{equation}
\index{connecting map!between $KE$-theory and $E$-theory, $\kee$|textbf}\index{$xXMXi$@$\kee$, map $KE_G\ra E_G$|textbf}
\end{thm}

We next have:

\begin{thm}\label{T:keeX}
$\kee: \mathbf{KE_G} \longrightarrow \mathbf{E_G}$ is a functor.
\end{thm}

\subsection{Puppe sequences and long exact sequences in $KE$-theory}

The last subsection of the paper discusses some partial results
that we have obtained related to the characterization of the excision
in $KE$-theory. Further study is necessary, but at least we can mention
two facts about the non-equivariant groups: the existence of Puppe
sequences \cite{Rsn82} and the split-exactness \cite{CuSk}.

\begin{thm}\label{T:Puppe_seq}
Let $A$, $B$, $D$ be graded separable \cs -algebras and $\varphi:A\ra B$
a \sh morphism, then we have exact sequences
$$
KE(D,A(0,1)) \xrightarrow{S\gvf_*} KE(D,B(0,1)) \xrightarrow{i_*}
   KE(D,C_{\gvf}) \xrightarrow{p_*} KE(D,A) \xrightarrow{\gvf_*} KE(D,B)
$$
and
$$
KE(A(0,1),D) \xleftarrow{S\gvf^*} KE(B(0,1),D) \xleftarrow{i^*}
   KE(C_{\gvf},D) \xleftarrow{p^*} KE(A,D) \xleftarrow{\gvf^*} KE(B,D).
$$
\end{thm}
\begin{proof}
In the statement  $C_{\gvf}$ is the mapping cone of $\gvf$,
$C_{\gvf}=\{ (a,f) \,|\, a\in A, f\in B[0,1), \gvf(a)=f(0) \}$,
$p: C_{\gvf}\ra A$ is the projection onto the first factor,
and $i:B(0,1) \ra C_{\gvf}$ is the inclusion $i(f)=(0,f)$.
The justification of the theorem is a minor modification of 
\cite[Thm.1.1]{CuSk}.
\end{proof}

\begin{thm}
Consider a short exact sequence of separable graded \cs -algebras
$$
\begin{CD}
0 @>>> J  @>{\;j\;}>>A  @>{\;q\;}>> B @>>> 0,
\end{CD}
\leqno{( \star)}
$$
such that $q$ admits a completely positive (grading preserving
and norm decreasing) cross-section. Then six-term exact sequences 
exist in $KE$-theory:
$$
\begin{CD}
KE(D,J)    @>{\;j_*\;}>>   KE(D,A)    @>{\;q_*\;}>>   KE(D,B)    \\
@A{\delta}AA             @.                       @VV{\delta}V \\
KE^1(D,B)  @<<{\;q_*\;}<   KE^1(D,A)  @<<{\;j_*\;}<   KE^1(D,J)    
\end{CD}
$$
and
$$
\begin{CD}
KE(J,D)    @<{\;j^*\;}<<   KE(A,D)    @<{\;q^*\;}<<   KE(B,D)    \\
@V{\delta}VV             @.                       @AA{\delta}A \\
KE^1(B,D)  @>>{\;q^*\;}>   KE^1(A,D)  @>>{\;j^*\;}>   KE^1(J,D). 
\end{CD}
$$
\end{thm}
\begin{proof}
Use Theorem \ref{T:Puppe_seq} and Corollary \ref{C:KE_equiv_via_KK} 
to get a $KE$-equivalence from the $KK$-equivalence \cite[Thm.2.1]{CuSk}
given by $e : J \ra C_q$, $e(x)=(j(x),0)$.
\end{proof}

\mic
Full details and proofs for the results contained in this section
will form the substance of a second paper.


\end{document}